\documentclass[12pt]{article}
\usepackage{array,epsfig,latexsym,tabularx}

\catcode`\@=12

\font\sectionfont=cmbx12 at 12pt

 \def\lskipamount{12pt}
 \def\lskip{\vskip\lskipamount plus3pt minus2pt}
 \def\lbreak{\par \ifdim\lastskip<\lskipamount
  \removelastskip \penalty-200 \lskip \fi}

 \def\lnobreak{\par \ifdim\lastskip<\lskipamount
  \removelastskip \penalty200 \lskip \fi}

\def\section#1{\refstepcounter{section}%
  \vskip 1.5truepc\centerline{%
  \hbox {{\sectionfont \thesection.  #1}}}
\vskip 0.5truepc
\par}

\def\subsection#1{\refstepcounter{subsection}%
  \vskip 1.5truepc
  \hbox {{\sectionfont \thesubsection.  #1}}
\vskip 0.5truepc
\par}

\def\thebibliography#1{\vskip 1.5pc{\centerline {\bf References}}\vskip 4pt
\list
 {[\arabic{enumi}]}{\settowidth\labelwidth{[#1]}\leftmargin\labelwidth
 \advance\leftmargin\labelsep
 \usecounter{enumi}}
 \def\newblock{\hskip .11em plus .33em minus .07em}
 \sloppy\clubpenalty4000\widowpenalty4000
 \sfcode`\.=1000\relax}

\def\@makeschapterhead#1{ \vspace*{50pt} { \parindent 0pt \raggedright
 \Huge \bf #1\par
 \nobreak \vskip 40pt } }

\def\chapter{\cleardoublepage
   \thispagestyle{plain}%
   \global\@topnum\z@

   \@afterindentfalse
   \secdef\@chapter\@schapter}

\def\@chapter[#1]#2{\ifnum \c@secnumdepth >\m@ne
        \refstepcounter{chapter}%
        \typeout{\@chapapp\space\thechapter.}%
        \addcontentsline{toc}{chapter}{\protect
        \numberline{\thechapter}#1}\else
      \addcontentsline{toc}{chapter}{#1}\fi
   \chaptermark{#1}%
   \addtocontents{lof}%
       {\protect\addvspace{10\p@}}
   \addtocontents{lot}%
       {\protect\addvspace{10\p@}}
   \if@twocolumn
           \@topnewpage[\@makechapterhead{#2}]%
     \else \@makechapterhead{#2}%
           \@afterheading
     \fi}

\def\@schapter#1{\if@twocolumn \@topnewpage[\@makeschapterhead{#1}]%
        \else \@makeschapterhead{#1}%
              \@afterheading\fi}

\begin{document}
 

\newtheorem{example}{Example}[section]
\newtheorem{note}[example]{Note}
\newtheorem{theorem}[example]{Theorem}
\newtheorem{corollary}[example]{Corollary}
\newtheorem{definition}[example]{Definition}
\newtheorem{proposition}[example]{Proposition}
\newtheorem{algorithm}[example]{Algorithm}
\newtheorem{lemma}[example]{Lemma}
\newtheorem{problem}[example]{Problem}
\newtheorem{conjecture}[example]{Conjecture}

 
\font\twelvesym=msbm10 at 12pt
\font\tensym=msbm10
\font\sevensym=msbm7
\font\fivesym=msbm5
\newfam\ssymfam
\textfont\ssymfam=\tensym
\scriptfont\ssymfam=\sevensym
\scriptscriptfont\ssymfam=\fivesym
\def\ssym{\fam\ssymfam\tensym}


\newcommand{\Z}{{\ssym Z}}
\newcommand{\N}{{\ssym N}}
\newcommand{\B}{{\mathcal B}}
\newcommand{\D}{{\mathcal D}}
\newcommand{\M}{{\mathcal M}}
\renewcommand{\P}{{\mathcal P}}

\newcommand{\infinity}{\infty}
\newcommand{\implies}{\Rightarrow}
\newcommand{\mod}{\mbox{mod}\,}
\newcommand{\wt}{{\rm wt\,}}
\newcommand{\owt}{{wt}}
\newcommand{\mwt}{\tilde{wt}}
\newcommand{\ochi}{{\chi}}
\newcommand{\mchi}{\tilde{\chi}}
\newcommand{\tkappa}{\tilde{\kappa}}
\newcommand{\Proof}{\medskip\noindent {\it Proof: }}
\newcommand{\cqfd}{\hfill $\Box$ \medskip}
\newcommand{\boldm}{\mbox{\boldmath$m$}}
\newcommand{\boldn}{\mbox{\boldmath$n$}}
\newcommand{\bolde}{\mbox{\boldmath$e$}}
\newcommand{\boldu}{\mbox{\boldmath$u$}}
\newcommand{\boldv}{\mbox{\boldmath$v$}}
\newcommand{\boldC}{\mbox{\boldmath$C$}}
\newcommand{\hboldC}{\hat{\mbox{\boldmath$C$}}}
\newcommand{\boldQ}{\mbox{\boldmath$Q$}}
\newcommand{\sboldm}{\mbox{\boldmath$\scriptstyle m$}}
\newcommand{\sboldn}{\mbox{\boldmath$\scriptstyle n$}}
\newcommand{\sbolde}{\mbox{\boldmath$\scriptstyle e$}}
\newcommand{\sboldu}{\mbox{\boldmath$\scriptstyle u$}}
\newcommand{\sboldC}{\mbox{\boldmath$\scriptstyle C$}}
\newcommand{\sboldQ}{\mbox{\boldmath$\scriptstyle Q$}}
\newcommand{\boldlambda}{\mbox{\boldmath$\lambda$}}
\newcommand{\boldDelta}{\mbox{\boldmath$\Delta$}}
\newcommand{\sboldDelta}{\mbox{\boldmath$\scriptstyle \Delta$}}
\newcommand{\wombat}{\rule[-6pt]{0pt}{46pt}}


\hyphenation{boson-ic 
             ferm-ion-ic 
	     para-ferm-ion-ic
             two-dim-ension-al
	     two-dim-ension-al}


\title{On the combinatorics of Forrester-Baxter models%
\thanks{Research supported by the Australian Research Council (ARC)}}
\author{Omar~Foda\thanks{foda@maths.mu.oz.au} \relax\
        and Trevor~A.~Welsh\thanks{trevor@maths.mu.oz.au}\\
        Department of Mathematics and Statistics,\\
        The University of Melbourne,
        Victoria, Australia.}

\maketitle
                 		      
\begin{abstract}
		       
We provide further boson-fermion $q$-polynomial identities 
for the `finitised' Virasoro characters $\chi^{p, p'}_{r,s}$
of the Forrester-Baxter minimal models $M(p,p')$, for certain 
values of $r$ and $s$. The construction is based on a detailed 
analysis of the combinatorics of the set $\P^{p, p'}_{a, b, c}(L)$ 
of $q$-weighted, length-$L$ Forrester-Baxter paths, whose
generating function $\chi^{p, p'}_{a, b, c}(L)$ provides 
a finitisation of $\chi^{p, p'}_{r,s}$. In this paper, 
we restrict our attention to the case where the startpoint 
$a$ and endpoint $b$ of each path both belong to the set of 
{\it \lq Takahashi lengths\rq}. In the limit $L \to \infty$, 
these polynomial identities reduce to $q$-series identities 
for the corresponding characters.

We obtain two closely related fermionic polynomial forms for 
each (finitised) character. The first of these forms uses the 
classical definition of the Gaussian polynomials, and includes 
a term that is a (finitised) character of a certain 
$M(\hat p,\hat p')$ where $\hat p'<p'$. We provide a combinatorial 
interpretation for this form using the concept of 
{\it \lq particles\rq}. The second form, which was first obtained 
using different methods by the Stony-Brook group, requires 
a modified definition of the Gaussian polynomials, and its 
combinatorial interpretation requires not only the concept of 
particles, but also the additional concept of {\it \lq particle
annihilation\rq}.

\end{abstract}

\vfill
\newpage

\setcounter{secnumdepth}{10}
\setcounter{section}{-1}

\section{Introduction}

\subsection{Motivation}

The physical spectrum of exactly-solvable lattice models 
can be described in the language of highest-weight infinite 
dim\-en\-sion\-al representations of affine and Virasoro 
algebras \cite{jimbo-miwa-book}. The characters of these 
representations are $q$-series that contain detailed 
information on the structure and symmetries of the 
corresponding models. In the following discussion, we wish 
to restrict attention to the characters of Virasoro 
highest-weight representations. 

The earliest known expressions for these characters are due 
to Feigen and Fuchs \cite{feigen-fuchs} and Rocha-Caridi 
\cite{rocha-caridi}. These expressions have alternating-signs. 
A number of years ago, the Stony Brook group 
discovered completely new expressions for the character 
formulae\footnote{For references to the original Stony Brook 
papers, please refer to \cite{stony-brook-review}.}. These 
expressions have constant-signs\footnote{The characterisation
of the different expressions of the characters as `alternating-sign'
and `constant-sign' $q$-series is valid only for Virasoro but
not for affine characters.}. 

For physical reasons that are beyond the scope of this work, 
the original alternating-sign expressions are also known 
as {\it `bosonic characters'}. Correspondingly, the con\-st\-ant-sign 
expressions are also known as {\it \lq fermionic characters\rq} 
\footnote{For a complete discussion of the physical motivation 
of the terms `bosonic' and `fermionic', please refer to the 
original literature on the subject as cited in 
\cite{stony-brook-review}.}. 

The structure of these new character formulae hints at the 
presence of a completely new formulation of exactly-solvable 
models\footnote{Analogous developments in the context of 
highest-weight representations of affine algebras also took 
place. They are outside the scope of this work.}. This possibility 
has attracted attention for a number of reasons. One of these 
reasons is the fact that certain physical problems, such as 
the long-distance asymptotics of the correlation functions, 
are too difficult to handle in the current formulation. Further, 
there are reasons to believe that the new formulation could be 
the right starting point to tackle them (see \cite{stony-brook-review}
and references therein). At a more technical level, the availability 
of two distinct formulations is mathematically enriching, as we 
can use one to learn about the other. 

However, although the bosonic characters are technically simple
to write down, and are completely known for all Virasoro 
representations, the structure of the fermionic characters 
is strictly-speaking known explicitly only in special cases, 
and generally only conceptually. In particular, the characters 
of the {\it `non-unitary'} Virasoro representations 
have turned out to be rather resistant to a complete formulation 
in fermionic form\footnote{The reason for that may of course 
eventually turn out to be the fact that we are not using the 
most efficient tools to tackle this problem.}.   

This work is part of a series of papers that aim at a complete 
and explicit derivation of the fermionic characters of a certain 
class of models first discussed by Forrester and Baxter 
\cite{forrester-baxter}. The characters of the Forrester-Baxter 
models correspond to the complete set of Virasoro characters of 
the discrete, though not necessarily unitary, Virasoro algebras 
with central charge $c < 1$, first discussed in \cite{bpz}. As 
such, they form the largest class of Virasoro characters with 
no $W$-symmetries.

As in previous works, our approach is purely combinatorial. 
Further, the exposition is self-contained, in the sense that 
we have included all concepts required in the derivations. 
Our main result is a combinatorial derivation of two related 
finitised fermionic forms for the characters of a certain class 
of Forrester-Baxter models. The first of these requires the use 
of the classical form of Gaussian polynomials and can be 
interpreted combinatorially using the concept of {\it particles}. 
The second has already appeared in the works of Berkovich, McCoy 
and Schilling \cite{bms}, requires the use of a modified form of 
Gaussian polynomials, and has a combinatorial interpretation in 
terms of particles and {\it particle annihilation}. 

In a forthcoming paper, we further extend and refine the techniques 
of this work to obtain a complete and explicit derivation of the 
fermionic characters of the complete set of Forrester-Baxter models 
\cite{foda-welsh-big}. 

\subsection{Overview of content of paper}

The aim of this paper is to obtain fermionic 
expressions for $\ochi^{p, p'}_{a, b, c}(L)$, the generating function 
for the set $\P^{p, p'}_{a, b, c}(L)$ of restricted length-$L$ paths
that have startpoint $a$ and endpoint $b$.

These functions\footnote{To be precise, a certain renormalisation 
thereof.} first arose in the calculation of one-point functions of 
the Forrester-Baxter models \cite{forrester-baxter}. The weighting 
originally assigned in \cite{forrester-baxter} to the paths is 
significantly different from that used here. The weighting described 
in the current paper arose by obtaining a \lq weight-preserving\rq\ 
bijection between partitions with prescribed hook-differences that 
were considered in \cite{abbbfv}, and the paths of 
\cite{forrester-baxter}. This bijection is described in \cite{flpw}.

The paths in $\P^{p, p'}_{a, b, c}(L)$ may be depicted on a $(p'-2)\times L$
grid that we refer to as the $(p,p')$-model,
as described in Section \ref{BandSec}. Of particular importance
is the shading of the $(p,p')$-model, which determines the
weights $\owt(h)$ that we assign to the paths $h$.

A bosonic expression for $\ochi^{p, p'}_{a, b, c}(L)$
is given in Section \ref{BosonicGenSec}.
This expression is readily proved using $L$-recurrence
relations \cite{forrester-baxter}, or by using the
generating function for partitions with prescribed hook-differences
given in \cite{abbbfv}, and the bijection of \cite{flpw}.
The polynomial $\ochi^{p, p'}_{a, b, c}(L)$ is seen to
be a finitisation of a Virasoro character.

In this paper, we tackle the particular cases where $a$ and $b$ are
each one of the Takahashi lengths $\mathcal T$, or one of
$\mathcal T'=\{p'-s:s\in\mathcal T\}$.
These values depend on $p$ and $p'$, and are defined in Section
\ref{ContFSec}.
Our methods and results are a common generalisation of those of
\cite{flpw,foda-welsh}.

On equating the bosonic expression for $\ochi^{p, p'}_{a, b, c}(L)$
with either of the fermionic expressions, we obtain
boson-fermion polynomial identities.
Taking the $L\to\infinity$ limit (using, for example, the variable change
employed in \cite{flpw,flw}), these become $q$-series identities.
Amongst them, in particular, are the Rogers-Ramanujan identities,
and their generalisations by Andrews and Gordon \cite{andrews-red-book}.
In fact, the techniques employed by Agarwal and Bressoud
\cite{agarwal-bressoud,bressoud} in their combinatorial proof
of the Andrews-Gordon identities provided the genesis of the
techniques employed here.

Before we develop a generalisation of Agarwal and Bressoud's
\lq Volcanic activity\rq, we define in Section \ref{WingSec},
a slightly different set $\P^{p, p'}_{a, b, e, f}(L)$
of paths, which have assigned pre-segments and post-segments
that are determined by $e,f\in\{0,1\}$.
Their generating function $\mchi^{p, p'}_{a, b, e, f}(L)$
is defined in terms of a path weighting that differs slightly
from that defined earlier.

The $\B$-transform, which is described in Section \ref{BTranSec},
enables $\mchi^{p, p'+p}_{a', b', e, f}(L')$, for certain $a',b'$
to be expressed in terms of $\mchi^{p, p'+p}_{a, b, e, f}(L)$.
We derive this transform combinatorially in three steps.
The first step is known as the $\B_1$-transform and
enlarges the features of a path, so that the resultant path
resides in a larger model.
The second step, referred to as a $\B_2(k)$-transform,
lengthens a path by appending $k$ pairs of segments to the path.
Each of these pairs is known as a particle.
The third step, the $\B_3(\lambda)$-transform deforms the
path in a particular way. This process may be viewed as the
particles {\em moving} through the path.
The resulting transformation of generating functions is
given in Corollary \ref{BijCor}.

In Section \ref{DTranSec}, we see that $\mchi^{p'-p, p'}_{a, b, 1-e, 1-f}(L)$
may be obtained from $\mchi^{p, p'}_{a, b, e, f}(L)$ in
a combinatorially trivial way.
This process is referred to as a $\D$-transform.
In fact, it is more convenient to use the $\D$-transform
combined with the $\B$-transform.
The resulting transformation of generating functions is
given in Corollary \ref{DijCor}.

To obtain a particular generating function
$\mchi^{p,p'}_{a,b,e,f}(L)$, where $p$ and $p'$ are co-prime, we begin
with one of the trivial generating
functions $\mchi^{1,3}_{a',b',e',f'}(L)$ given in Lemma
\ref{SeedLem}, and perform a sequence of $\B$- and $\B\D$-transforms.
This sequence is determined by the continued fraction
of $p'/p$ which is described in Section \ref{ContFSec}.

In fact, a basic application of the transforms does not
generate all elements of $\P^{p,p'}_{a,b,e,f}(L)$ in some cases.
In these instances, the set generated is deficient in the full set
of paths that do not rise above (or below) a certain height.
Various results obtained in Section \ref{SegmentSec} enable us
to keep track of this height.
Lemma \ref{SegmentLem} shows that this height bounds a portion
of the $(p,p')$-model which is identical to a smaller
$(\hat p,\hat p')$-model. This property enables (in one case),
the final generating function to be expressed using the
generating function for paths in the $(\hat p,\hat p')$-model.

Section \ref{ExtAttSec} provides one further ingredient
for the final construction.
There, it is shown how appending or removing the first segment
of the path affects the generating function.

Everything is now in place to carry out the proof of the main results.
These results are stated in Section \ref{ResultsSec}.
We provide two similar expressions for $\ochi^{p,p'}_{a,b,c}(L)$.
These are Theorems \ref{Ferm1Thrm} and \ref{Ferm2Thrm}.
The first of these makes use of the classical definition of the
Gaussian polynomial:
\begin{equation}\label{Gaussian}
\left[ {A \atop B} \right]_q =
\left\{
  \begin{array}{cl}
     \frac{\displaystyle (q)_A}{\displaystyle (q)_{A-B}(q)_B} &
          \quad\mbox{if } 0\le B\le A;\\[3mm]
     0 &
          \quad\mbox{otherwise},
  \end{array} \right.
\end{equation}
where $(q)_0=1$ and $(q)_n=\prod_{i=1}^n(1-q^i)$ for $n>0$.
In some cases, the expression also includes a term
$\ochi^{\hat p,\hat p'}_{a,b,c}(L)$ for $\hat p'<\hat p$.
Thus this expression may be viewed as a recursive
fermionic expression for $\ochi^{p,p'}_{a,b,c}(L)$.
In the cases where this additional term is not present
(for $a$ and $b$ further restricted in a certain way),
the expressions were first stated in \cite{berkovich-mccoy}.

The expression of Theorem \ref{Ferm2Thrm} makes use of a modified
definition of the Gaussian polynomial (\cite{gasper-rahman}):
\begin{equation}\label{ModGaussian}
\left[ {A \atop B} \right]^\prime_q =
\left\{
  \begin{array}{cl}
     \frac{\displaystyle (q^{A-B+1})_B}{\displaystyle (q)_B} &
          \quad\mbox{if } 0\le B;\\[3mm]
     0 &
          \quad\mbox{otherwise},
  \end{array} \right.
\end{equation}
where $(z)_0=1$ and $(z)_n=\prod_{i=0}^{n-1}(1-zq^i)$ for $n>0$.
These expressions were first presented and proved in \cite{bms}.
In fact, invoking the definition (\ref{ModGaussian}) is somewhat overkill,
since the only value of $\left[ {A \atop B} \right]^\prime$ that we
require that differs from $\left[ {A \atop B} \right]$ is
$\left[ {-1\atop 0}\hbox to0pt{$\phantom{A \atop B}$\hss}\right]^\prime=1$.

In \cite{bms}, expressions for $\ochi^{p,p'}_{a,b,c}(L)$ are presented,
where $b$ is now any value with $1\le b\le p'-1$.
However, only $a\in\mathcal T\cup\mathcal T'$ is still permitted.
In \cite{foda-welsh-big}, we show that it is Theorem \ref{Ferm1Thrm},
and not Theorem \ref{Ferm2Thrm},
that generalises to provide fermionic expressions for
the most general $\ochi^{p,p'}_{a,b,c}(L)$.

The remainder of Section \ref{FermSec} is concerned with the detailed
derivation of the expression for first 
$\mchi^{p,p'}_{a,b,e,f}(L)$, and then converting it to
$\ochi^{p,p'}_{a,b,c}(L)$.
Section \ref{MNsysSec} describes the $\boldm\boldn$-system which
aids the actual evaluation of the fermionic expressions obtained.
Section \ref{BadFermSec} describes how the proof for Theorem \ref{Ferm1Thrm}
modifies to provide a proof for Theorem \ref{Ferm2Thrm}.
Here we see that the appearance of
$\left[ {-1\atop 0}\hbox to0pt{$\phantom{A \atop B}$\hss}\right]^\prime$
may be viewed in terms of \lq particle annihilation\rq.

\section{Paths}\label{PathSec}

\subsection{Paths and the $(p,p')$-model}\label{BandSec}

Let $p$ and $p'$ be positive co-prime integers for which $0<p<p'$.
Then, given $a,b,c,L\in\Z_{\ge0}$ such that $1\le a,b,c\le p'-1$,
$b=c\pm1$, $L + a - b \equiv 0$ ($\mod2$),
a path $h \in \P^{p, p'}_{a, b, c}(L)$ is a sequence
$h_0, h_1, h_2, \ldots, h_L,$ of integers such that:

\begin{enumerate}
\item $1 \le h_i \le p'-1$  for $0 \le i \le L$, 
\item $h_{i+1} = h_i \pm 1$ for $0 \le i <   L$, 
\item $h_0 = a, 
       h_L = b.$
\end{enumerate}

\noindent Note that the values of $p$ and $c$ do not feature in 
the above restrictions. As described below, they specify how the 
elements of $\P^{p, p'}_{a, b, c}(L)$ are weighted.

The integers $h_0,h_1,h_2,\ldots,h_L,$ are readily depicted as a sequence 
of {\it heights} on a two-dimensional $L \times (p'-2)$ grid. Adjacent 
heights are connected by {\it line segments} passing from $(i,h_i)$ to 
$(i+1,h_{i+1})$ for $0 \le i < L$. 

\begin{figure}[ht]
\centerline{\epsfig{file=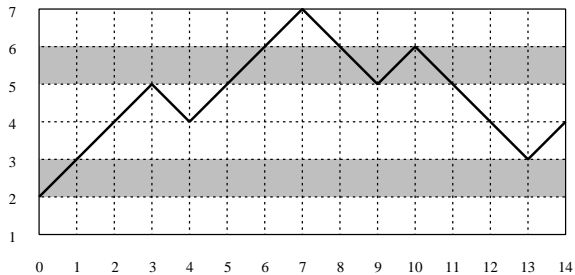}}
\caption{Typical path.}
\label{TypicalFig}
\medskip
\end{figure}

Scanning the path from left to right, each of these line segments points 
either in the NE direction or in the SE direction. Fig.~\ref{TypicalFig}
shows a typical path in the set $\P^{3,8}_{2, 4, 3}(14)$.
The shadings in Fig.~\ref{TypicalFig} are explained below.

In the grid introduced above, the horizontal strip between adjacent
heights is referred to as a {\em band}.  There are $p'-2$ bands.
The $h$th band lies between heights $h$ and $h+1$.

We now assign a parity to each band: the $h$th band is said to be an {\it even}
band if $\lfloor hp/p'\rfloor=\lfloor (h+1)p/p'\rfloor$; and an {\it odd}
band if $\lfloor hp/p'\rfloor\ne\lfloor (h+1)p/p'\rfloor$.
The array of odd and even bands so obtained will be referred to
as the $(p,p')$-model.
It may immediately be deduced that the $(p,p')$-model has $p'-p-1$ even
bands and $p-1$ odd bands.
In addition, it is easily shown that for $1\le r<p$, the band lying
between heights $\lfloor rp'/p\rfloor$ and $\lfloor rp'/p\rfloor+1$
is odd: it will be referred to as the $r$th odd band.

When drawing the $(p,p')$-model, we distinguish the bands by shading
the odd bands. This was done in Fig.~\ref{TypicalFig} for the
$(3,8)$-model.

We note that the band structure of the $(p,p')$-model is up-down
symmetrical, and that if $p'>2p$ then the 1st band and the $(p'-2)$th
band are both even, and there are no two adjacent odd bands.

For $2\le a\le p'-2$, we say that $a$ is {\em interfacial\/} if
$\lfloor (a+1)p/p'\rfloor=\lfloor(a-1)p/p'\rfloor+1$.
Thus $a$ is interfacial if and only if $a$ lies between an
odd and even band in the $(p,p')$-model.
Thus for the case of the $(3,8)$-model depicted in Fig.~\ref{TypicalFig},
$a$ is interfacial for $a=2,3,5,6$.
Note that if $a$ is interfacial, the odd band that it borders
is the $\lfloor (a+1)p/p'\rfloor$th.

As is easily seen, the $(p'-p,p')$-model  
differs from the $(p,p')$-model in that each band has changed parity.
It follows that if $a$ is interfacial in the $(p,p')$-model
then $a$ is also interfacial in the $(p'-p,p')$-model.

\subsection{Weighting the paths}

\begin{figure}[ht]
\centerline{\epsfig{file=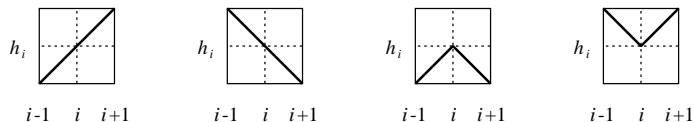}}
\caption{Vertex shapes.}
\label{VertexFig}
\medskip
\end{figure}

Given a path $h$ of length $L$, for $1\le i<L$, the values
of $h_{i-1}$, $h_i$ and $h_{i+1}$
determine the shape of the vertex at the point $i$.
The four possible shapes are given in Fig.~\ref{VertexFig}.

\noindent
The four types of vertices shown in Fig.~\ref{VertexFig} are referred
to as a {\em straight-up} vertex, a {\em straight-down} vertex,
a {\em peak-up} vertex and a {\em peak-down} vertex respectively.
Each vertex is also assigned a parity: this is the parity of the band
in which the segment between $(i,h_i)$ and $(i+1,h_{i+1})$ lies.
Thus, there are eight types of paritied vertex.

For paths $h\in{\P}^{p,p'}_{a,b,c}(L)$, we define $h_{L+1}=c$,
whereupon the shape and parity of the vertex at $i=L$ is well-defined.

The weight function for the paths is best specified in terms of a
$(x,y)$-coordinate system which is inclined at $45^o$ to the original
$(i,h)$-coordinate
system and whose origin is at the path's initial point at $(i=0,h=a)$.
Specifically,
$$
x={{i-(h-a)} \over {2}},\qquad y={{i+(h-a)} \over {2}}.
$$
Note that at each step in the path, either $x$ or $y$ is incremented
and the other is constant. In this system, the path depicted
in Fig.~\ref{TypicalFig} has its first few coordinates at
$(0,0)$, $(0,1)$, $(0,2)$, $(0,3)$, $(1,3)$, $(1,4)$, $(1,5)$,
$(1,6)$, $(2,6)$, $\ldots$

Now, for $1\le i\le L$, we define the weight 
$c_i= c(h_{i-1},h_i,h_{i+1})$ of the $i$th vertex
according to its shape, its parity and its $(x,y)$-coordinate,
as specified in Table~\ref{WtsTable}.

\begin{table}[ht]
\begin{center}
\begin{tabular}{|c|@{\hspace{3mm}}c@{\hspace{3mm}}|c|@{\hspace{3mm}}
c@{\hspace{3mm}}|}
\hline
Vertex&
${c}_i$&
Vertex&
${c}_i$\\
\hline\hline\wombat
{\epsfig{file=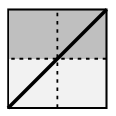}}&
\raisebox{14pt}[0pt]{$x$}&
{\epsfig{file=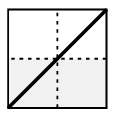}}&
\raisebox{14pt}[0pt]{$0$}\\
\hline\wombat
{\epsfig{file=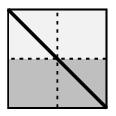}}&
\raisebox{14pt}[0pt]{$y$}&
{\epsfig{file=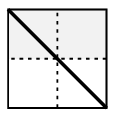}}&
\raisebox{14pt}[0pt]{$0$}\\
\hline\wombat
{\epsfig{file=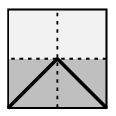}}&
\raisebox{14pt}[0pt]{$0$}&
{\epsfig{file=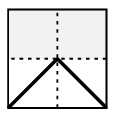}}&
\raisebox{14pt}[0pt]{$x$}\\
\hline\wombat
{\epsfig{file=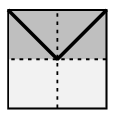}}&
\raisebox{14pt}[0pt]{$0$}&
{\epsfig{file=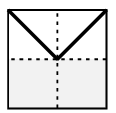}}&
\raisebox{14pt}[0pt]{$y$}\\
\hline
\end{tabular}
\end{center}
\medskip
\caption{Vertex weights.}
\label{WtsTable}
\end{table}

\noindent
In Table~\ref{WtsTable}, the lightly shaded bands can be either even or odd
bands (or when $h_i=p'-1$ or $h_i=1$ in the lowermost four cases,
not a band in the model at all).
Note that for each vertex shape, only one parity case has non-zero
weight in general.
We shall refer to those four vertices, with assigned parity, for
which in general, the weight is non-zero, as {\em scoring} vertices.
The other four vertices will be termed {\em non-scoring}.

We now define:

\begin{equation}\label{WtDef}
\owt(h)=\sum_{i=1}^L c_i.
\end{equation}

To illustrate this procedure, consider again the path $h$ depicted in
Fig.~\ref{TypicalFig}. The above table indicates that there are scoring
vertices at $i=3$, $4$, $5$, $7$, $8$, $13$ and $14$.
This leads to
$$
\owt(h)=0+3+1+1+6+7+6=24.
$$

The generating function $\ochi^{p,p'}_{a,b,c}(L)$ for
the set of paths ${\P}^{p,p'}_{a,b,c}(L)$ is defined to be:
\begin{equation}\label{PathGenDef1}
\ochi^{p,p'}_{a,b,c}(L;q)
=\sum_{h\in{\P}^{p,p'}_{a,b,c}(L)} q^{\owt(h)}.
\end{equation}

\noindent
Often, we drop the base $q$ from the notation so that
$\ochi^{p,p'}_{a,b,c}(L)=\ochi^{p,p'}_{a,b,c}(L;q)$.
The same will be done for other functions without comment.

\subsection{Bosonic generating function}\label{BosonicGenSec}

By setting up recurrence relations for $\ochi^{p,p'}_{a,b,c}(L)$,
it may be readily verified that:


\begin{eqnarray}\label{FinRochaEq}
\chi^{p,p'}_{a,b,c}(L)&=&
\sum_{\lambda=-{\infinity}}^{\infinity}
q^{\lambda^{2} p p'+ \lambda (p'r-pa)}
\left[ {L \atop {\frac{L+a-b}{2}}-p'\lambda} \right]_q\\[0.5mm]
&&\qquad\qquad\qquad
-\sum_{\lambda=-\infinity}^\infinity
q^{(\lambda p+r)(\lambda p'+a)}
\left[ {L \atop {\frac{L+a-b}{2}}-p'\lambda-a} \right]_q,\nonumber
\end{eqnarray}

\noindent where 

\begin{equation}\label{groundstatelabel}
r=\lfloor pc/p'\rfloor+(b-c+1)/2.
\end{equation}

In the limit $L\mapsto\infinity$, we obtain

\begin{equation}\label{ChiLimEq}
\lim_{L\to\infinity}
\chi^{p,p'}_{a,b,c}(L) = \chi^{p,p'}_{r,a},
\end{equation}

\noindent where $r$ is defined in (\ref{groundstatelabel}) and

\begin{equation}\label{RochaEq}
\chi^{p, p'}_{r, s}=
{\frac{1}{(q)_\infty}}\sum_{\lambda=-\infty}^\infty
(q^{\lambda^2pp'+\lambda(p'r-ps)}-q^{(\lambda p+r)(\lambda p'+s)})
\end{equation}

\noindent is, up to a normalisation, the Rocha-Caridi expression 
\cite{rocha-caridi} for the Virasoro character of central charge 
$c = 1 - {6(p' - p)^2}/{p p'}$
and conformal dimension $\Delta^{p, p'}_{r, s} = 
{\left( (p' r - p s)^2 - (p' - p)^2 \right)}/{4 p p'} $.
Therefore, $\chi^{p,p'}_{a,b,c}(L)$ provides a finite analogue of 
the character $\chi^{p,p'}_{r, a}$.

\goodbreak
\section{Winged generating functions}\label{WingSec}
\nobreak

For $h\in\P^{p,p'}_{a,b,c}(L)$, the values of $b$ and $c$ serve to specify
a path {\em post-segment} that extends between $(L,b)$ and $(L+1,c)$.
We now define another set of paths which specifies both the direction
of a post-segment and a {\em pre-segment}.

Let $p$ and $p'$ be positive co-prime integers for which $0<p<p'$.
Then, given $a,b,L\in\Z_{\ge0}$ such that $1\le a,b\le p'-1$,
$L + a - b \equiv 0$ ($\mod2$), and $e,f\in\{0,1\}$,
a path $h \in \P^{p, p'}_{a, b, e, f}(L)$ is a sequence
$h_0, h_1, h_2, \ldots, h_L,$ of integers such that:

\begin{enumerate}
\item $1 \le h_i \le p'-1$  for $0 \le i \le L$, 
\item $h_{i+1} = h_i \pm 1$ for $0 \le i <   L$, 
\item $h_0 = a, 
       h_L = b.$
\end{enumerate}

\noindent If $f=0$ (resp.\ $f=1$) then the post-segment of each
$h \in \P^{p, p'}_{a, b, e, f}(L)$ is defined to be in the NE (resp.\ SE)
direction.
If $e=0$ (resp.\ $e=1$) then the pre-segment of each
$h \in \P^{p, p'}_{a, b, e, f}(L)$ is defined to be in the SE (resp.\ NE)
direction.
This enables a shape and a parity to be assigned to both the
zeroth and the $L$th vertices of $h$.
For $h \in \P^{p, p'}_{a, b, e, f}(L)$, we define $e(h)=e$ and $f(h)=f$.

We now define a weight $\mwt(h)$, for
$h \in \P^{p, p'}_{a, b, e, f}(L)$.
For $1\le i<L$, set $\tilde c_i=c(h_{i-1},h_i,h_{i+1})$ as above.
Then, set
\begin{displaymath}
\tilde c_L=
\left\{
  \begin{array}{ll}
x \qquad &\mbox{if } h_L-h_{L-1}=\phantom{-}1 \mbox{ and } f(h)=1;\\[1.5mm]
y \qquad &\mbox{if } h_L-h_{L-1}=-1 \mbox{ and } f(h)=0;\\[1.5mm]
0 \qquad &\mbox{otherwise,}
  \end{array} \right.
\end{displaymath}
where $(x,y)$ is the coordinate of the $L$th vertex of $h$.
We then designate this vertex as scoring if it is a peak
vertex ($h_L=h_{L-1}-(-1)^{f(h)}$), and as non-scoring otherwise.

We define:

\begin{equation}\label{WtDef2}
\mwt(h)=\sum_{i=1}^L \tilde c_i.
\end{equation}

Consider the corresponding path $h'\in\P^{p,p'}_{a,b,c}(L)$
with $c=b+(-1)^f$, defined by $h_i'=h_i$ for $0\le i\le L$.
{}From Table \ref{WtsTable}, we see that $\mwt(h)=\owt(h')$
if the post-segment of $h$ lies in an even band.

In what follows, we work entirely in terms of $\mwt(h)$,
and the generating functions that we derive from it.
Only at the end of our work, do we revert back to $\owt(h)$
to obtain fermionic expressions for $\ochi^{p,p'}_{a,b,c}(L)$.

Define the generating function
\begin{equation}\label{PathGenDef2}
\mchi^{p,p'}_{a,b,e,f}(L;q)
=\sum_{h\in{\P}^{p,p'}_{a,b,e,f}(L)} q^{\mwt(h)},
\end{equation}
\noindent
where $\mwt(h)$ is given by (\ref{WtDef2}).
Of course, $\mchi^{p,p'}_{a,b,0,f}(L)=\mchi^{p,p'}_{a,b,1,f}(L)$.

\subsection{Striking sequence of a path}\label{StrikeSec}

For each path $h$, define $\pi(h)\in\{0,1\}$ to be the
parity of the band between heights $h_0$ and $h_1$
(if $L(h)=0$, we set $h_1=h_0+(-1)^{f(h)}$).
Thus, for the path $h$ shown in Fig.~\ref{TypicalFig}, we have $\pi(h)=1$.
In addition, define $d(h)=0$ when $h_1-h_0=1$ and $d(h)=1$
when $h_1-h_0=-1$.
We then see that if $e(h)+d(h)+\pi(h)\equiv0\,(\mod2)$ then
the $0$th vertex is a scoring vertex,
and if $e(h)+d(h)+\pi(h)\equiv1\,(\mod2)$ then it is
a non-scoring vertex.

Now consider each path $h\in{\P}^{p,p'}_{a,b,e,f}(L)$ as a sequence of
straight lines, alternating in direction between NE and SE.
Then, reading from the left, let the lines be of lengths
$w_1$, $w_2$, $w_3,\ldots,w_l,$ for some $l$, with $w_i>0$ for $1\le i\le l$.
Thence $w_1+w_2+\cdots+w_l=L(h)$, where $L(h)=L$ is the length of $h$.

For each of these lines, the last vertex will be considered to be
part of the line but the first will not. Then, the $i$th of these 
lines contains $w_i$ vertices, the first $w_i-1$ of which are
straight vertices. Then write $w_i=a_i+b_i$ so that $b_i$ is the 
number of scoring vertices in the $i$th line. The striking sequence 
of $h$ is then the array:
\begin{displaymath}\label{HseqDef}
\left(\begin{array}{ccccc}
  a_1&a_2&a_3&\cdots&a_l\\ b_1&b_2&b_3&\cdots&b_l
 \end{array}\right)^{(e(h),f(h),d(h))}.
\end{displaymath}

With $\pi=\pi(h)$, $e=e(h)$ and $d=d(h)$, we define
\begin{displaymath}
m(h)=
\left\{
  \begin{array}{ll}
       (e+d+\pi)\,\mod2+\sum_{i=1}^l a_i
          \qquad
          &\mbox{if } L>0;\\[1.5mm]
       \vert f-e\vert
          \qquad
          &\mbox{if } L=0,
  \end{array} \right.
\end{displaymath}
whence $m(h)$ is the number of non-scoring vertices possessed
by $h$ (altogether, $h$ has $L(h)+1$ vertices).
We also define $\alpha(h)=(-1)^d((w_1+w_3+\cdots)-(w_2+w_4+\cdots))$
and for $L>0$,
\begin{displaymath}
\beta(h)=
\left\{
  \begin{array}{l}
(-1)^d((b_1+b_3+\cdots)-(b_2+b_4+\cdots))\\
          \hskip40mm \mbox{if } e+d+\pi\equiv0\,(\mod2);\\[1.5mm]
(-1)^d((b_1+b_3+\cdots)-(b_2+b_4+\cdots))+(-1)^e\\
          \hskip40mm \mbox{otherwise.}
  \end{array} \right.
\end{displaymath}
For $L=0$, we set $\beta(h)=f-e$.

For example, for the path shown in Fig.~\ref{TypicalFig} for which
$d(h)=0$ and $\pi(h)=1$, the striking sequence is:
$$
\def\qua{\hskip5pt}
\left(
{2\qua0\qua1\qua1\qua1\qua2\qua0\atop
 1\qua1\qua2\qua1\qua0\qua1\qua1}
\right)^{(e,1,0)}.
$$
In this case, $m(h)=8-e$, $\alpha(h)=2$, and $\beta(h)=2-e$.

We note that given the startpoint $h_0=a$ of the path, the path
can be reconstructed from its striking sequence\footnote{We only
need $w_1,w_2,\ldots,w_l$ together with $d$.}.
In particular, $h_L=b=a+\alpha(h)$. In addition, the nature of the
final vertex may be deduced from $a_l$ and $b_l$%
\footnote{Thus the value of $f$ in the striking sequence
is redundant --- we retain it for convenience.}

\begin{lemma}\label{WtHashLem}
Let the path $h$ have the striking sequence
$\left({a_1 \atop b_1}\,{a_2 \atop b_2}\,{a_3 \atop b_3}\,
 {\cdots\atop\cdots}\,{a_l\atop b_l} \right)^{(e,f,d)}\!,$
with $w_i=a_i+b_i$ for $1\le i\le l$.
Then
$$
\mwt(h)=\sum_{i=1}^l b_i(w_{i-1}+w_{i-3}+\cdots+w_{1+i\bmod2}).
$$
\end{lemma}

\Proof For $L=0$, both sides are clearly $0$.
So assume $L>0$.
First consider $d=0$.
For $i$ odd, the $i$th line is in the NE direction and
its $x$-coordinate is $w_2+w_4+\cdots+w_{i-1}$. By the prescription
of the previous section, and the definition of $b_i$, this line
thus contributes $b_i(w_2+w_4+\cdots+w_{i-1})$ to the weight
$\mwt(h)$ of $h$. Similarly, for $i$ even, the $i$th line is in
the SE direction and contributes $b_i(w_1+w_3+\cdots+w_{i-1})$
to $\mwt(h)$.
The lemma then follows for $d=0$.
The case $d=1$ is similar.
\cqfd
\medskip

\subsection{Path parameters}

We make the following definitions:
\begin{displaymath}
\begin{array}{ll}
\alpha^{p,p'}_{a,b}&=\: b-a;\\[0.5mm]
\beta^{p,p'}_{a,b,e,f} &=\:
  \left\lfloor\frac{bp}{p'}\right\rfloor 
       - \left\lfloor\frac{ap}{p'}\right\rfloor + f-e;\\[2.5mm]
\delta^{p,p'}_{a,e} &=\:
  \left\{
    \begin{array}{ll}
       0 \quad &
            \mbox{if }
              \left\lfloor\frac{(a+(-1)^e)p}{p'}\right\rfloor
              =\left\lfloor\frac{ap}{p'}\right\rfloor;\\[3mm]
       1 \quad &
            \mbox{if }
              \left\lfloor\frac{(a+(-1)^e)p}{p'}\right\rfloor
              \ne\left\lfloor\frac{ap}{p'}\right\rfloor.
    \end{array} \right.
\\
\end{array}
\end{displaymath}

\noindent
(The superscripts of $\alpha^{p,p'}_{a,b}$ are superfluous, of course.)
It may be seen that the value of $\delta^{p,p'}_{a,e}$
gives the parity of the band in which the path pre-segment resides.

\begin{lemma}\label{BetaConstLem}
Let $h\in{\P}^{p,p'}_{a,b,e,f}(L)$.
Then $\alpha(h)=\alpha^{p,p'}_{a,b}$ and $\beta(h)=\beta^{p,p'}_{a,b,e,f}$.
\end{lemma}

\Proof That $\alpha(h)=\alpha^{p,p'}_{a,b}$ follows immediately
from the definitions.

The second result is proved by induction on $L$.
If $h\in{\P}^{p,p'}_{a,b,e,f}(0)$ then $a=b$,
whence $\beta^{p,p'}_{a,b,e,f}=f-e=\beta(h)$,
immediately from the definitions.

For $L>0$, let $h\in{\P}^{p,p'}_{a,b,e,f}(L)$ and assume that
the result holds for all $h'\in{\P}^{p,p'}_{a,b',e,f'}(L-1)$.
We consider a particular $h'$ by setting $h_i'=h_i$ for
$0\le i<L$, $b'=h_{L-1}$ and choosing $f'\in\{0,1\}$ so that
$f'=0$ if either $b-b'=1$ and the $L$th segment of $h$ lies in an
even band, or $b-b'=-1$ and the $L$th segment of $h$ lies in an
odd band; and $f'=1$ otherwise.
It may easily be checked that the $(L-1)$th vertex of $h'$ is
scoring if and only if the $(L-1)$th vertex of $h$ is scoring.
Then, from the definition of $\beta(h)$, we see that:
\begin{displaymath}
\beta(h)=
  \left\{
    \begin{array}{ll}
       \beta(h')+1 \quad &
            \mbox{if } b-b'=\phantom{-}1 \mbox{ and } f=1;\\
       \beta(h')-1 \quad &
            \mbox{if } b-b'=-1 \mbox{ and } f=0;\\
       \beta(h') \quad &
            \mbox{otherwise.}
    \end{array} \right.
\end{displaymath}
The induction hypothesis gives
$\beta(h')= \lfloor b'p/p'\rfloor - \lfloor ap/p'\rfloor +f'-e$.
Then when the $L$th segment of $h$ lies in an even band so
that $\lfloor bp/p'\rfloor=\lfloor b'p/p'\rfloor$, consideration
of the four cases of $b-b'=\pm1$ and $f\in\{0,1\}$ shows that
$\beta(h)= \lfloor bp/p'\rfloor - \lfloor ap/p'\rfloor +f-e$.
When the $L$th segment of $h$ lies in an odd band so
that $\lfloor bp/p'\rfloor=\lfloor b'p/p'\rfloor+b-b'$, consideration
of the four cases of $b-b'=\pm1$ and $f\in\{0,1\}$ again shows that
$\beta(h)= \lfloor bp/p'\rfloor - \lfloor ap/p'\rfloor +f-e$.
The result follows by induction.
\cqfd
\medskip

\subsection{Scoring generating functions}

We now define a generating function for paths that have a
particular number of non-scoring vertices.
First define ${\P}^{p,p'}_{a,b,e,f}(L,m)$ to be the subset of
${\P}^{p,p'}_{a,b,e,f}(L)$ comprising those paths $h$ for which
$m(h)=m$. Then define:
\begin{equation}\label{ResPathGenDef}
\ochi^{p,p'}_{a,b,e,f}(L,m;q)
=\sum_{h\in{\P}^{p,p'}_{a,b,e,f}(L,m)} q^{\mwt(h)}.
\end{equation}

\begin{lemma}\label{ResPathGenLem}
Let $\beta=\beta^{p,p'}_{a,b,e,f}$. Then
\begin{displaymath}
\ochi^{p,p'}_{a,b,e,f}(L)
=\sum_{\scriptstyle m\equiv L+\beta\atop
  \scriptstyle\strut(\mbox{\scriptsize\rm mod}\,2)}
\ochi^{p,p'}_{a,b,e,f}(L,m).
\end{displaymath}
\end{lemma}

\Proof Let $h\in\P^{p,p'}_{a,b,e,f}(L)$.
We claim that $m(h)+L(h)+\beta(h)\equiv0\,(\mod2)$.
This will follow from showing that $L(h)-m(h)+(-1)^{d(h)}\beta(h)$
is even.
If $h$ has striking sequence
$\left({a_1 \atop b_1}\:{a_2 \atop b_2}\:{a_3 \atop b_3}\:
 {\cdots\atop\cdots}\:{a_l\atop b_l} \right)^{(e,f,d)},$
then $L(h)-m(h)=(b_1+b_2+\cdots+b_l)-(e+d+\pi)\,\mod2$,
where $\pi=\pi(h)$.
For $e+d+\pi\equiv0\,(\mod2)$, we immediately obtain
$L(h)-m(h)+(-1)^{d}\beta(h)=2(b_1+b_3+\ldots)$.
For $e+d+\pi\not\equiv0\,(\mod2)$, we obtain
$L(h)-m(h)+(-1)^{d}\beta(h)=2(b_1+b_3+\ldots)-1+(-1)^{d+e}$,
whence the claim is proved in all cases.
The lemma then follows, once it is noted, via Lemma \ref{BetaConstLem},
that $\beta(h)=\beta^{p,p'}_{a,b,e,f}$.
\cqfd
\medskip

\begin{note}
Since each element of ${\P}^{p,p'}_{a,b,e,f}(L,m)$ has $L+1$
vertices, it follows that
$\ochi^{p,p'}_{a,b,e,f}(L,m)$ is non-zero only if $0\le m\le L+1$.
Therefore the sum in Lemma \ref{ResPathGenLem} may be further
restricted to $0\le m\le L+1$.
\end{note}


\subsection{A seed}

The following result provides a seed on which the results of
later sections will act.

\begin{lemma}\label{SeedLem}
If $L\ge0$ is even then:
\begin{displaymath}
\ochi^{1,3}_{1,1,0,0}(L,m)=
\ochi^{1,3}_{2,2,1,1}(L,m)=
\delta_{m,0}q^{\frac14 L^2}.
\end{displaymath}
If $L>0$ is odd then:
\begin{displaymath}
\ochi^{1,3}_{1,2,0,1}(L,m)=
\ochi^{1,3}_{2,1,1,0}(L,m)=
\delta_{m,0}q^{\frac14 (L^2-1)}.
\end{displaymath}
\end{lemma}

\Proof
The $(1,3)$-model comprises one even band.
Thus when $L$ is even, there is precisely one $h\in\P^{1,3}_{1,1,0,0}(L)$.
It has $h_i=1$ for $i$ even, and $h_i=2$ for $i$ odd.
We see that $h$ has striking sequence 
$\left({0 \atop 1}\:{0 \atop 1}\:{0 \atop 1}\:
 {\cdots\atop\cdots}\:{0\atop 1} \right)^{(0,0,0)}$
and $m(h)=0$.
Lemma \ref{WtHashLem} then yields
$\mwt(h)=0+1+1+2+2+3+\cdots+(\frac12L-1)+\frac12L=(L/2)^2$, as required.

The other expressions follow in a similar way.
\cqfd
\medskip

\subsection{Partitions}

A partition $\lambda=(\lambda_1,\lambda_2,\ldots,\lambda_k)$ is a
sequence of $k$
integer parts $\lambda_1,\lambda_2,\ldots,\lambda_k,$ satisfying
$\lambda_1\ge\lambda_2\ge\cdots\ge\lambda_k>0$.
It is to be understood that $\lambda_i=0$ for $i>k$.
The weight $\wt(\lambda)$ of $\lambda$ is given
by $\wt(\lambda)=\sum_{i=1}^k\lambda_i$.

We define ${\mathcal Y}(k,m)$ to be the set of all partitions $\lambda$
with at most $k$ parts, and for which $\lambda_1\le m$.
A proof of the following well known result may be found in
\cite{andrews-red-book}.

\begin{lemma}\label{PartitionGenLem}
The generating function,
$$
\sum_{\lambda\in{\mathcal Y}(k,m)} q^{\wt(\lambda)}=
\left[{m+k\atop m}\right]_q.
$$
\end{lemma}


\section{The $\B$-transform}\label{BTranSec}

In this section, we introduce the $\B$-transform which maps
paths ${\P}^{p,p'}_{a,b,e,f}(L)$ into
${\P}^{p,p'+p}_{a',b',e,f}(L')$ for certain $a'$, $b'$
and various $L'$.

The band structure of the $(p,p'+p)$-model is easily obtained from
that of the $(p,p')$-model. Indeed, according to Section \ref{BandSec},
for $1\le r<p$, the $r$th odd band of the $(p,p'+p)$-model lies
between heights $\lfloor r(p'+p)/p\rfloor=\lfloor rp'/p\rfloor+r$ and
$\lfloor r(p'+p)/p\rfloor+1=\lfloor rp'/p\rfloor+r+1$.
Thus the height of the $r$th odd band in the $(p,p'+p)$-model
is $r$ greater than that in the $(p,p')$-model.
Therefore, the $(p,p'+p)$-model may be obtained from the $(p,p')$-model
by increasing the distance between neighbouring odd bands by one
unit and appending an extra even band to both the top and the bottom
of the grid. For example, compare the $(3,8)$-model of
Fig.~\ref{TypicalFig} with the $(3,11)$-model of Fig.~\ref{B_1aFig}.

The $\B$-transform has three components, which we refer to as
{\em path-dilation}, {\em particle-insertion}, and
{\em particle-motion}. These three components will also
be known as the $\B_1$-, $\B_2$- and $\B_3$-transforms
respectively.
In fact, particle-insertion is dependent on a parameter $k\in\Z_{\ge0}$,
and particle-motion is dependent on a partition $\lambda$
that has certain restrictions.
Consequently, we sometimes refer to particle-insertion
and particle-motion as $\B_2(k)$- and $\B_3(\lambda)$-transforms
respectively.
Then, combining the $\B_1$-, $\B_2(k)$- and $\B_3(\lambda)$-transforms
produces the $\B(k,\lambda)$-transform.

\subsection{Path-dilation} The $\B_1$-transform acts on a
path $h\in{\P}^{p,p'}_{a,b,e,f}(L)$ to yield a path
$h^{(0)}\in{\P}^{p,p'+p}_{a',b',e,f}(L^{(0)})$,
for certain $a'$, $b'$ and $L^{(0)}$.
First, the starting point $a'$ of the new path $h^{(0)}$ is specified
to be:
$$
a'=a+\left\lfloor\frac{ap}{p'}\right\rfloor + e.
$$
If $r=\lfloor ap/p'\rfloor$ then $r$ is the number of odd bands below
$h=a$ in the $(p,p')$-model.
Since the height of the $r$th odd band in the $(p,p'+p)$-model
is $r$ greater than that in the $(p,p')$-model,
we thus see that under path-dilation, the height of the
startpoint above the next lowermost odd band
(or if there isn't one, the bottom of the grid)
has either increased by one or remained constant.

We define $d(h^{(0)})=d(h)$.
The above definition specifies that $e(h^{(0)})=e(h)$ and
$f(h^{(0)})=f(h)$.

In the case that $L=0$ and $e=f$, we specify $h^{(0)}$ by
setting $L^{(0)}=L(h^{(0)})=0$. When $L=0$ and $e\ne f$, we
leave the action of the $\B_1$-transform on $h$ undefined
(it will not be used in this case).
Thus in Lemmas \ref{BHashLem}, \ref{BWtLem}, \ref{WtShiftLem},
\ref{BresLem}, \ref{BijLem}, \ref{BDresLem}, \ref{DijLem} and
and Corollary \ref{EndPtCor}, we implicitly exclude consideration of the
case $L=0$ and $e\ne f$.
However, it must be considered in the proofs of
Corollaries \ref{BijCor} and \ref{DijCor}.

In the case $L>0$ consider, as in Section \ref{StrikeSec}, $h$ to
comprise $l$ straight lines that alternate in direction, the
$i$th of which is of length $w_i$ and possesses $b_i$ scoring
vertices. $h^{(0)}$ is then defined to comprise $l$ straight
lines that alternate in direction (since $d(h^{(0)})=d(h)$, the direction
of the first line in $h^{(0)}$ is the same as that in $h$),
the $i$th of which has length
\begin{displaymath}
w'_i=
  \left\{
    \begin{array}{l}
         \hbox to 25mm{$w_i+b_i$\hss}
         \mbox{if } i\ge2 \mbox{ or } e(h)+d(h)+\pi(h)\equiv0\,(\mod2);\\[1mm]
       w_1+b_1+2\pi(h)-1\\
         \hskip25mm
         \mbox{if } i=1 \mbox{ and } e(h)+d(h)+\pi(h)\not\equiv0\,(\mod2).
    \end{array} \right.
\end{displaymath}
In particular, this determines $L^{(0)}=L(h^{(0)})$ and
$b'=h^{(0)}_{L^{(0)}}$.

As an example, consider the path $h$ shown in Fig.~\ref{TypicalFig}
as an element of $\P^{3,8}_{2, 4, e, 1}(14)$.
Here $d(h)=0$, $\pi(h)=1$ and $\lfloor ap/p'\rfloor=0$.

\begin{figure}[ht]
\centerline{\epsfig{file=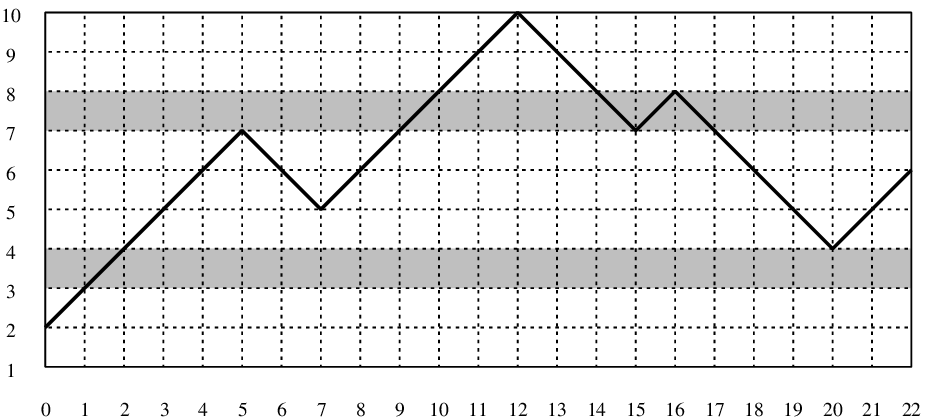}}
\caption{}
\label{B_1aFig}
\medskip
\end{figure}

Thus when $e=0$, the action of path-dilation on $h$ produces the
path given in Fig.~\ref{B_1aFig}.

\noindent
This path is an element of $\P^{3,11}_{2, 6, e, 1}(22)$.

\begin{figure}[ht]
\centerline{\epsfig{file=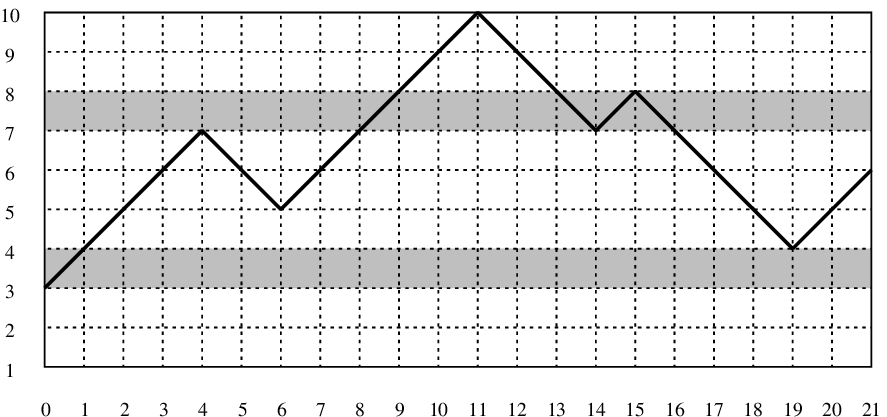}}
\caption{}
\label{B_1bFig}
\medskip
\end{figure}

When $e=1$, the action of path-dilation on $h$ produces the 
element of $\P^{3,11}_{3, 6, e, 1}(21)$ given in Fig.~\ref{B_1bFig}.

The situation at the start point may be considered as
falling into one of eight cases, corresponding to
$e(h),d(h),\pi(h)\in\{0,1\}$.\footnote{Theses cases may be seen
to correspond to the eight cases of vertex type as listed
in Table~\ref{WtsTable}.}
In Table~\ref{BatStartTable}, we illustrate the four cases that arise
when $d(h)=0$ (the four cases for $d(h)=1$ may be obtained
from these by an up-down reflection and changing the value
of $e(h)$).\footnote{The examples here
are such that $w_1\ge3$.}

\begin{table}[ht]
\bigskip\noindent
\raisebox{20pt}[0pt]{$e(h)=0;\atop\pi(h)=0:$}
$\quad\epsfig{file=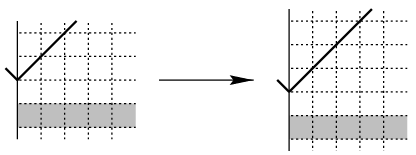}$
\qquad
\raisebox{20pt}[0pt]{$e(h)=1;\atop\pi(h)=0:$}
$\quad\epsfig{file=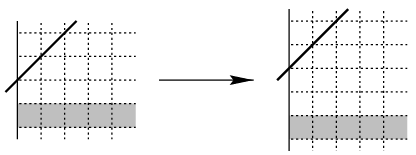}$

\bigskip\noindent
\raisebox{20pt}[0pt]{$e(h)=0;\atop\pi(h)=1:$}
$\quad\epsfig{file=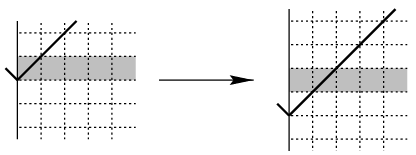}$
\qquad
\raisebox{20pt}[0pt]{$e(h)=1;\atop\pi(h)=1:$}
$\quad\epsfig{file=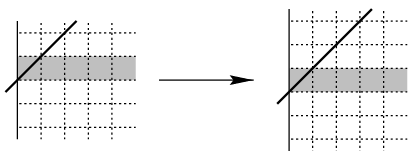}$
\medskip

\caption{$\B_1$-transforms at the startpoint.}
\label{BatStartTable}
\end{table}

\begin{lemma}\label{StartPtLem}
Let $1\le p<p'$, $1\le a<p'$,
$e\in\{0,1\}$ and $a'=a+\lfloor ap/p'\rfloor+e$.
Then $\lfloor a'p/(p'+p)\rfloor=\lfloor ap/p'\rfloor$ and
$\delta^{p,p'+p}_{a',e}=0$.
%
\end{lemma}

\Proof Let $r=\lfloor ap/p'\rfloor$ whence $p'r\le pa<p'(r+1)$.
Then, for $x\in\{0,1\}$, we have
$(p'+p)r\le p(a+r+x)<(p'+p)r+p'+xp$, so that
$\lfloor (a+r+x)p/(p'+p)\rfloor=r$.
In particular, $\lfloor a'p/(p'+p)\rfloor=r$,
and $\lfloor (a+r+e+(-1)^e)p/(p'+p)\rfloor=r$.
Thus $r=\lfloor a'p/(p'+p)\rfloor=\lfloor (a'+(-1)^e)p/(p'+p)\rfloor$
which gives the required results.
%
\cqfd
\medskip

\noindent This result asserts, amongst other things, that the pre-segment
of $h^{(0)}$ always lies in an even band. This is also evident from
Table~2.

\begin{note}\label{Start1Note}
The action of path-dilation on $h\in\P^{p,p'}_{a,b,e,f}(L)$ yields a
path $h^{(0)}\in\P^{p,p'+p}_{a',b',e,f}(L^{(0)})$
that has, including the vertex at $i=0$, no adjacent scoring vertices,
except in the case where $\pi(h)=1$ {\em and} $e(h)=d(h)$, when a single
pair of scoring vertices occurs in $h^{(0)}$ at $i=0$ and $i=1$.

Also note that $\pi(h^{(0)})=\pi(h)$ unless $\pi(h)=1$ {\em and}
$e(h)=d(h)$, in which case $\pi(h^{(0)})=0$.
\end{note}

Now compare the $i$th line of $h^{(0)}$ (which has length $w'_i$)
with the $i$th line of $h$ (which has length $w_i$).
Now for the sake of the following argument, assume that there are odd
bands immediately below (i.e.~between heights $0$ and $1$), and 
immediately above (i.e.~between heights $p'-1$ and $p'$)
the $(p,p')$-model and do likewise for the $(p,p'+p)$-model.

If the lines in question are in the NE direction, we claim that
the height of the final vertex of that in $h^{(0)}$ above
the next lower odd band is one greater than that in $h$.
If the lines in question are in the SE direction, we claim that
the height of the final vertex of that in $h^{(0)}$ below
the next higher odd band is one greater than that in $h$.
In particular, if either the first or last segment
of the $i$th line is in an odd band, then the corresponding segment
of $h^{(0)}$ lies in the same odd band.

We also claim that if that of $h$ has a straight vertex that passes
into the $k$th odd band in the $(p,p')$-model then
that of $h^{(0)}$ has a straight vertex that passes into the
$k$th odd band in the $(p,p'+p)$-model.

These claims follow because in passing from the $(p,p')$-model
to the $(p,p'+p)$-model, the distance between neighbouring odd bands
has increased by one, and because the length of each line has
increased by one for every scoring vertex and possibly a small
adjustment made to the length of the first line.
In effect, a new straight vertex has been inserted
immediately prior to each scoring vertex
and, if $e(h)+d(h)+\pi(h)\not\equiv0\,(\mod2)$,
adjusting the length of the resulting first line by $2\pi(h)-1$.

\begin{lemma}\label{BHashLem}

Let $h \in \P^{p, p'}_{a,b,e,f}(L)$
have striking sequence
$\left({a_1\atop b_1}{a_2\atop b_2}{a_3\atop b_3}
 {\cdots\atop\cdots}{a_l\atop b_l} \right)^{(e,f,d)}\!,$
and let
$h^{(0)} \in \P^{p,p'+p}_{a',b',e,f}(L^{(0)})$ be obtained from
the action of the $\B_1$-transform on $h$.
If $e(h)+d(h)+\pi(h)\equiv0\,(\mod2)$ then $h^{(0)}$ has
striking sequence:
\begin{displaymath}
\left(\begin{array}{ccccc}
  a_1+b_1&a_2+b_2&a_3+b_3&\cdots&a_l+b_l\\
  b_1&b_2&b_3&\cdots&b_l
 \end{array}\right)^{(e,f,d)},
\end{displaymath}
and if $e(h)+d(h)+\pi(h)\not\equiv0\,(\mod2)$ then $h^{(0)}$ has
striking sequence:
\begin{displaymath}
\left(\begin{array}{ccccc}
  a_1+b_1+\pi-1&a_2+b_2&a_3+b_3&\cdots&a_l+b_l\\
  b_1+\pi&b_2&b_3&\cdots&b_l
 \end{array}\right)^{(e,f,d)}.
\end{displaymath}
Moreover, if $m=m(h)$:
\begin{itemize}
\item $m(h^{(0)})=L$;
\item $L^{(0)}=
\left\{
  \begin{array}{ll}
2L-m+2 \quad &
          \mbox{if } \pi=1 \mbox{ and } e=d,\\[1.5mm]
2L-m \quad & \mbox{otherwise};
  \end{array} \right. $
\item $\alpha(h^{(0)})=\alpha(h)+\beta(h)$;
\item $\beta(h^{(0)})=\beta(h)$.
\end{itemize}
\end{lemma}

\Proof The form of the striking sequence for $h^{(0)}$ follows
because, for $i>1$, every scoring vertex in the $i$th line of $h$
accounts for an extra non-scoring vertex in that line. The same is true
when $i=1$, except in the case $(e(h)+d(h)+\pi(h))\equiv1$
(throughout this paper, in proofs, we take all equivalences, modulo 2.)
when the length of the new $1$st line
becomes $a_1+2b_1+2\pi-1$. That there are $b_1+\pi$
scoring vertices in this case, follows from examining Table~2.

Let $e=e(h)$, $d=d(h)$, $\pi=\pi(h)$ and $\pi'=\pi(h^{(0)})$.
Then $e(h^{(0)})=e$ and $d(h^{(0)})=d$.

If $(e+d+\pi)\equiv0$ then $(e+d+\pi')\equiv0$ by Note \ref{Start1Note}.
Thereupon $m^{(0)}=\sum_{i=1}^l(a_i+b_i)=L$.
Additionally, $L^{(0)}=\sum_{i=1}^l(a_i+2b_i)=2L-\sum_{i=1}^la_i=2L-m$.
That $\beta(h^{(0)})=\beta(h)$ and
$\alpha(h^{(0)})=\alpha(h)+\beta(h)$ both follow immediately
in this case.

On the other hand, if $(e+d+\pi)\not\equiv0$ then
$\pi=0\implies e\ne d$ and $\pi=1\implies e=d$.
In each instance, Note \ref{Start1Note} implies that $\pi'=0$.
Thereupon,
$m^{(0)}=(e+d+\pi')\,\mod2+\pi-1+\sum_{i=1}^l(a_i+b_i)
=\sum_{i=1}^l(a_i+b_i)=L$.
Additionally,
$L^{(0)}=2\pi-1+\sum_{i=1}^l(a_i+2b_i)=2L-(1+\sum_{i=1}^la_i)+2\pi=2L-m+2\pi$.
This is the required value.
Now in this case,
$\beta(h)=(-1)^d((b_1+b_3+\cdots)-(b_2+b_4+\cdots))+(-1)^e$.
When $\pi=0$ so that $(e+d+\pi')\equiv1$ then
$\beta(h^{(0)})=\beta(h)$ follows immediately. 
When $\pi=1$, we have
$\beta(h^{(0)})=(-1)^d((b_1+1+b_3+\cdots)-(b_2+b_4+\cdots))$.
$\beta(h^{(0)})=\beta(h)$ now follows in this case because
$(e+d+\pi)\not\equiv0$ implies that $e=d$.
Finally,
$\alpha(h^{(0)})=\alpha(h)+(-1)^d((b_1+b_3+\cdots)-(b_2+b_4+\cdots))
+(-1)^d(2\pi-1)$.
Since $(-1)^d(2\pi-1)=-(-1)^d(-1)^\pi=(-1)^e$,
the lemma then follows.
\cqfd
\medskip

\begin{corollary}\label{EndPtCor}
Let $h \in{\P}^{p,p'}_{a,b,e,f}(L)$ 
and
$h^{(0)} \in \P^{p,p'+p}_{a',b',e,f}(L^{(0)})$ be the path obtained 
by the action of the $\B_1$-transform on $h$.
Then $a'=a+\lfloor ap/p'\rfloor+e$ and $b'=b+\lfloor bp/p'\rfloor+f$.
\end{corollary}

\Proof $a'=a+\lfloor ap/p'\rfloor+e$ is by definition.
Lemma \ref{BHashLem} gives $\alpha(h^{(0)})=\alpha(h)+\beta(h)$,
whence Lemma \ref{BetaConstLem} implies that
$\alpha^{p,p'+p}_{a',b'}=\alpha^{p,p'}_{a,b}+\beta^{p,p'}_{a,b,e,f}$.
Expanding this gives
$b'-a'=b-a+\lfloor bp/p'\rfloor-\lfloor ap/p'\rfloor+f-e$,
whence $b'=b+\lfloor bp/p'\rfloor+f$.
%
\cqfd
\medskip

\noindent
The above result implies that the $\B_1$-transform maps
${\P}^{p,p'}_{a,b,e,f}(L)$ into a set of paths that have the same
startpoint as one another and the same endpoint as one another.
However, the lengths of these paths are not necessarily equal.
We also see that the transformation of the endpoint is analogous
to that which occurs at the startpoint.
In particular, Lemma \ref{StartPtLem} implies that
$\delta^{p,p'+p}_{b',f}=0$ so that the path post-segment of $h^{(0)}$
always resides in an even band.
For the four cases where $h_L=h_{L-1}-1$, the $\B_1$-transform
affects the endpoint as in Table \ref{BatEndTable}
(the value $\pi'(h)$ is the parity of the band in which the
$L$th segment of $h$ lies).

\begin{table}[ht]
\bigskip\noindent
\hbox{
\raisebox{20pt}[0pt]{$f(h)=0;\atop\pi'(h)=0:$}
$\quad\epsfig{file=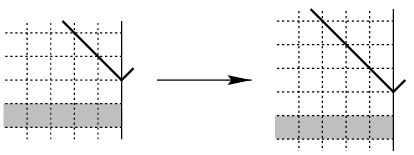}$
\qquad
\raisebox{20pt}[0pt]{$f(h)=1;\atop\pi'(h)=0:$}
$\quad\epsfig{file=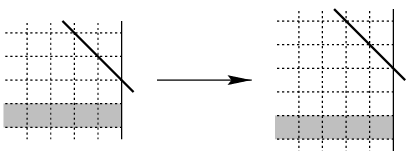}$}

\bigskip\noindent
\hbox{
\raisebox{20pt}[0pt]{$f(h)=0;\atop\pi'(h)=1:$}
$\quad\epsfig{file=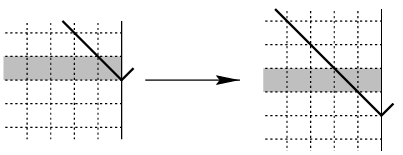}$
\qquad
\raisebox{20pt}[0pt]{$f(h)=1;\atop\pi'(h)=1:$}
$\quad\epsfig{file=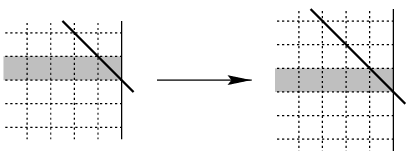}$}
\medskip

\caption{$\B_1$-transforms at the endpoint.}
\label{BatEndTable}
\end{table}

\begin{lemma}\label{ParamLem}
Let $1\le p<p'$, $1\le a,b<p'$, $e,f\in\{0,1\}$,
$a'=a+\lfloor ap/p'\rfloor+e$, and $b'=b+\lfloor bp/p'\rfloor+f$.
Then
$\alpha^{p,p'+p}_{a',b'}=\alpha^{p,p'}_{a,b}+\beta^{p,p'}_{a,b,e,f}$
and
$\beta^{p,p'+p}_{a',b',e,f}=\beta^{p,p'}_{a,b,e,f}$.
\end{lemma}

\Proof
Lemma \ref{StartPtLem} implies that
$\lfloor a'p/(p'+p)\rfloor=\lfloor ap/p'\rfloor$,
$\lfloor b'p/(p'+p)\rfloor=\lfloor bp/p'\rfloor$.
The results then follow immediately from the definitions.
\cqfd
\medskip

\begin{lemma}\label{BWtLem}
Let $h \in{\P}^{p,p'}_{a,b,e,f}(L)$ 
and
$h^{(0)} \in \P^{p,p'+p}_{a',b',e,f}(L^{(0)})$ be the path obtained 
by the action of the $\B_1$-transform on $h$.
Then
$$
\mwt(h^{(0)})=\mwt(h)+{\frac{1}{4}}\left((L^{(0)}-m^{(0)})^2-\beta^2\right),
$$
where $m^{(0)}=m(h^{(0)})$ and $\beta=\beta^{p,p'}_{a,b,e,f}$.
\end{lemma}

\Proof Let $h$ have striking sequence
$\left({a_1\atop b_1}\:{a_2\atop b_2}\:{a_3\atop b_3}\:
 {\cdots\atop\cdots}\:{a_l\atop b_l} \right)^{(e,f,d)}$,
and let $\pi=\pi(h)$.
If $(e+d+\pi)\equiv0\,(\mod2)$, then
Lemmas \ref{BHashLem} and \ref{WtHashLem} show that
$$
\mwt(h^{(0)})-\mwt(h)=(b_1+b_3+b_5+\cdots)(b_2+b_4+b_6+\cdots).
$$
Via Lemma \ref{BHashLem}, we obtain
$L^{(0)}-m^{(0)}=L-m(h)=b_1+b_2+\cdots+b_l$.
Then since $\beta(h)=\pm((b_1+b_3+b_5+\cdots)-(b_2+b_4+b_6+\cdots))$,
it follows that
$$
\mwt(h^{(0)})-\mwt(h)={{1} \over {4}}((L^{(0)}-m^{(0)})^2-\beta(h)^2).
$$

If $(e+d+\pi)\not\equiv0\,(\mod2)$, then
Lemmas \ref{BHashLem} and \ref{WtHashLem} show that
\begin{eqnarray*}
\mwt(h^{(0)})-\mwt(h)&=&(2\pi-1+b_1+b_3+b_5+\cdots)(b_2+b_4+b_6+\cdots)\\
&=&{{1} \over {4}}((L^{(0)}-m^{(0)})^2-\beta(h)^2),
\end{eqnarray*}
the second equality resulting because
$L^{(0)}-m^{(0)}=L-m(h)+2\pi=b_1+b_2+\cdots+b_l+2\pi-1$ and
\begin{eqnarray*}
\beta(h)&=&(-1)^d((b_1+b_3+b_5+\cdots)-(b_2+b_4+b_6+\cdots))+(-1)^e\\
&=&\pm((2\pi-1+b_1+b_3+b_5+\cdots)-(b_2+b_4+b_6+\cdots)),
\end{eqnarray*}
on using $(-1)^{e+d}=-(-1)^\pi=2\pi-1$.

Finally, Lemma \ref{BetaConstLem} gives $\beta(h)=\beta^{p,p'}_{a,b,e,f}=
\beta$.
\cqfd
\medskip

\subsection{Particle insertion}\label{InsertSec}

Let $p'>2p$ so that the $(p,p')$-model has no two neighbouring odd bands,
and let $\delta^{p,p'}_{a',e}=0$.
Then if $h^{(0)}\in\P^{p,p'}_{a',b',e,f}(L^{(0)})$,
the pre-segment of $h^{(0)}$ lies in an even band.
By {\em inserting a particle} 
into $h^{(0)}$, we mean displacing $h^{(0)}$ two positions to
the right and inserting two segments: the leftmost of these
is in the NE (resp.~SE) direction if $e=0$ (resp.~$e=1$),
and the rightmost is in the opposite direction, which is thus
the direction of the pre-segment of $h^{(0)}$.
In this way, we obtain a path $h^{(1)}$ of length $L^{(0)}+2$.
We assign $e(h^{(1)})=e$ and $f(h^{(1)})=f$.
Note also that $d(h^{(1)})=e$ and $\pi(h^{(1)})=0$.

Thereupon, we may repeat this process of particle insertion.
After inserting $k$ particles into $h^{(0)}$,
we obtain a path $h^{(k)}\in\P^{p,p'}_{a',b',e,f}(L^{(0)}+2k)$.
We say that $h^{(k)}$ has been obtained 
by the action of a $\B_2(k)$-transform on $h^{(0)}$.

\begin{figure}[ht]
\centerline{\epsfig{file=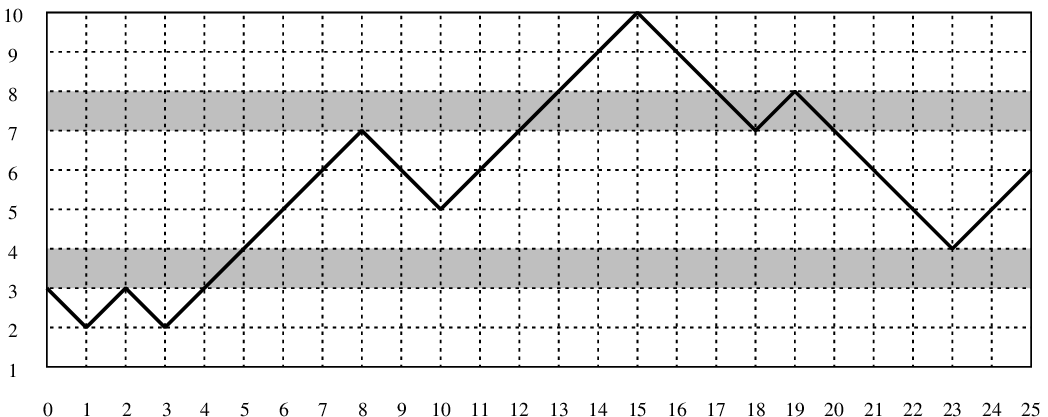}}
\caption{}
\label{B_2fig}
\medskip
\end{figure}

In the case of the element of  $\P^{3,11}_{3, 6, 1, 1}(21)$ shown
in Fig.~\ref{B_1bFig}, the insertion of two particles produces the
element of $\P^{3,11}_{3, 6, 1, 1}(25)$ shown in Fig.~\ref{B_2fig}.

\begin{lemma}\label{WtShiftLem}
Let $h \in \P^{p,p'}_{a,b,e,f}(L)$.
Apply a $\B_1$-transform to $h$ to obtain the path
$h^{(0)} \in \P^{p,p'+p}_{a',b',e,f}(L^{(0)})$. 
Then obtain $h^{(k)}\in\P^{p,p'+p}_{a',b',e,f}(L^{(k)})$ by
applying a $\B_2(k)$-transform to $h^{(0)}$.
If $m^{(k)}=m(h^{(k)})$,
then $L^{(k)}=L^{(0)}+2k$, $m^{(k)}=m^{(0)}$ and
\begin{equation}\label{Proof1Eq}
\mwt(h^{(k)})=\mwt(h)+{ {1} \over {4}}((L^{(k)}-m^{(k)})^2-\beta^2),
\end{equation}
where $\beta=\beta^{p,p'}_{a,b,e,f}$.
\end{lemma}
\Proof
That $L^{(k)}=L^{(0)}+2k$ follows immediately from the
definition of a $\B_2$-transform.
Lemma \ref{BWtLem} yields:
$$
\mwt(h^{(0)})=\mwt(h)
   +{{1} \over {4}}\left((L^{(0)}-m(h^{(0)}))^2-\beta^2\right).
$$
Let the striking sequence of $h^{(0)}$ be
$\left({a_1\atop b_1}\:{a_2\atop b_2}\:
 {\cdots\atop\cdots}\:{a_l\atop b_l} \right)^{(e,f,d)},$
and let $\pi=\pi(h^{(0)})$.

If $e=d$, we are restricted to the case $\pi=0$,
since $\delta^{p,p'+p}_{a',e}=0$ by Lemma \ref{StartPtLem}.
The striking sequence of $h^{(1)}$ is then
$\left({0\atop1}\:{0\atop1}\:{a_1\atop b_1}\:{a_2\atop b_2}\:
 {\cdots\atop\cdots}\:{a_l\atop b_l} \right)^{(e,f,e)}$.
Thereupon $m(h^{(1)})=\sum_{i=1}^l a_i=m(h^{(0)})$.
In this case, Lemma \ref{WtHashLem} shows that
$\mwt(h^{(1)})-\mwt(h^{(0)})=1+b_1+b_2+\cdots+b_l=L^{(0)}-m^{(0)}+1$.

If $e\ne d$, the striking sequence of $h^{(1)}$ is
$\left({0\atop1}\:{a_1+1-\pi\atop b_1+\pi}\:{a_2\atop b_2}\:
 {\cdots\atop\cdots}\:{a_l\atop b_l} \right)^{(e,f,e)}$.
Then $m(h^{(1)})=1-\pi+\sum_{i=1}^la_i$
which equals $m(h^{(0)})=(e+d+\pi)\,\mod2+\sum_{i=1}^la_i$
for both $\pi=0$ and $\pi=1$.
Here, Lemma \ref{WtHashLem} shows that
$\mwt(h^{(1)})-\mwt(h^{(0)})=\pi+b_1+b_2+\cdots+b_l$.
Since $L^{(0)}-m^{(0)}=-(e+d+\pi)\,\mod2+b_1+b_2+\cdots+b_l$,
we once more have
$\mwt(h^{(1)})-\mwt(h^{(0)})=L^{(0)}-m^{(0)}+1$.

Repeated application of these results, yields
$m(h^{(k)})=m(h^{(0)})$ and
$$
\mwt(h^{(k)})=\mwt(h^{(0)})+k\left(L^{(0)}-m(h^{(0)})\right)+k^2.
$$
Then, on using (\ref{Proof1Eq}) and $L^{(k)}=L^{(0)}+2k$, the lemma follows.
\cqfd
\medskip

\subsection{Particle moves}\label{PartMovesSec}

In this section, we once more restrict to the case
$p'>2p$ so that the $(p,p')$-model has no two neighbouring odd bands,
and consider only paths $h\in\P^{p,p'}_{a',b',e,f}(L')$, where
$\delta^{p,p'}_{a',e}=\delta^{p,p'}_{b',f}=0$.

We specify six types of local deformations of a path.
These deformations will be known as {\em particle moves}.
In each of the six cases, 
a particular sequence of four segments of a path is changed to a 
different sequence, the remainder of the path being unchanged. The 
moves are as follows --- the path portion to the left of the arrow is 
changed to that on the right:

\bigskip
\centerline{\epsfig{file=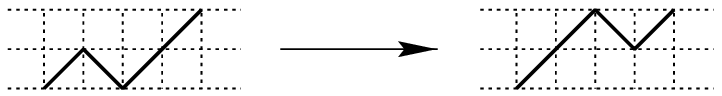}}
\nobreak
\centerline{Move 1.}
\bigskip

\centerline{\epsfig{file=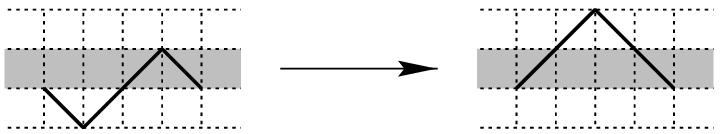}}
\nobreak
\centerline{Move 2.}
\bigskip

\centerline{\epsfig{file=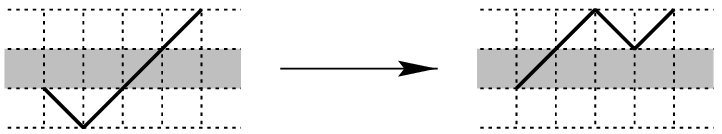}}
\nobreak
\centerline{Move 3.}
\bigskip

\centerline{\epsfig{file=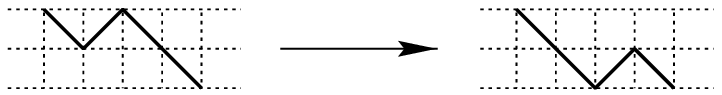}}
\nobreak
\centerline{Move 4.}
\bigskip

\centerline{\epsfig{file=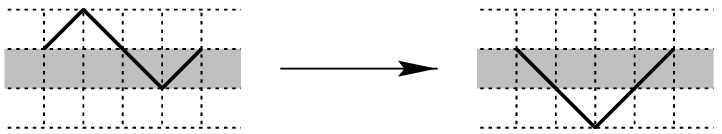}}
\nobreak
\centerline{Move 5.}
\bigskip

\centerline{\epsfig{file=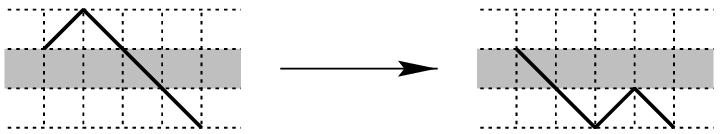}}
\nobreak
\centerline{Move 6.}
\bigskip

\noindent Since $p'>2p$, each odd band is straddled by a pair
of even bands. Thus, there is no impediment to enacting moves
2 and 5 for paths in $\P^{p,p'}_{a,b,e,f}(L)$.

Note that moves 4--6 are inversions of moves 1--3. Also note 
that moves 2 and 3 (likewise moves 5 and 6) may be considered 
to be the same move since in the two cases, the same sequence 
of three edges is changed.

In addition to the six moves described above, we permit certain
deformations of a path close to its left and right extremities in
certain circumstances.
Each of these moves will be referred to as an {\em edge-move}.
They, together with their validity, are as follows:
 
\bigskip
\hbox to 0mm{\qquad If $e=1$:\hss}\hskip28mm
\raisebox{0pt}[0pt]{\epsfig{file=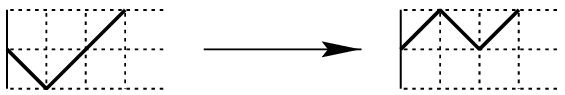}}
\par\medskip\nobreak
\centerline{Edge-move 1.}
\bigskip
 
\bigskip
\hbox to 0mm{\qquad If $e=0$:\hss}\hskip28mm
\raisebox{0pt}[0pt]{\epsfig{file=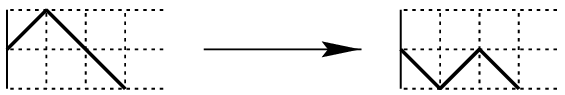}}
\par\medskip\nobreak
\centerline{Edge-move 2.}
\bigskip
 
\bigskip
\hbox to 0mm{\qquad If $f=0$:\hss}\hskip28mm
\raisebox{0pt}[0pt]{\epsfig{file=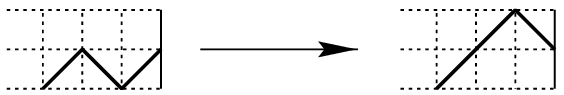}}
\par\medskip\nobreak
\centerline{Edge-move 3.}
\bigskip
 
\bigbreak
\hbox to 0mm{\qquad If $f=1$:\hss}\hskip28mm
\nobreak
\raisebox{0pt}[0pt]{\epsfig{file=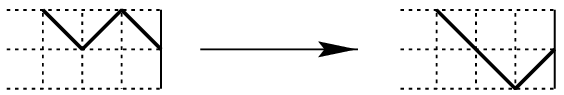}}
\par\medskip\nobreak
\centerline{Edge-move 4.}
\bigskip

In fact, the above four edge-moves may be considered as instances of
moves 1 and 4 described beforehand, if for edge-moves 1 and 2, we append
the appropriate pre-segment to the path,
and for edge-moves 3 and 4, we append the appropriate post-segment
to the path.

\begin{lemma}\label{WtChngLem}
Let the path $\hat h$ differ from the path $h$ in that
four consecutive segments have changed according to one of the
six moves described above, or in that three consecutive
segments have changed according to one of the four edge-moves
described above (subject to their restrictions).
Then
$$
\mwt(\hat h)=\mwt(h)+1.
$$
Additionally, $L(\hat h)=L(h)$ and $m(\hat h)=m(h)$.
\end{lemma}
\Proof For each of the six moves and four edge-moves, take the
$(x,y)$-coordinate of the leftmost point of the depicted portion
of $h$ to be $(x_0,y_0)$. Now consider 
the contribution to the weight of the three vertices in question 
before and after the move (although the vertex at $(x_0,y_0)$
may change, its contribution doesn't). In each 
of the ten cases, the contribution is $x_0+y_0+1$ before the move 
and $x_0+y_0+2$ afterwards. Thus $\mwt(\hat h)=\mwt(h)+1$.
The other statements are immediate on inspecting all ten moves.
\cqfd
\medskip

Now observe that for each of the ten moves specified above, the
sequence of path segments before the move consists of an adjacent
pair of scoring vertices followed by a non-scoring vertex.
The specified move replaces this combination with a non-scoring
vertex followed by two scoring vertices.
As anticipated above, the pair of adjacent scoring vertices is
viewed as a particle. Thus each of the above ten moves describes
a particle moving to the right by one step.

\begin{figure}[ht]
\centerline{\epsfig{file=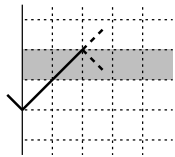}}
\caption{Not a particle}
\label{NotPFig}
\medskip
\end{figure}

When $p'>2p$, so that there are no two adjacent odd bands in the
$(p,p')$-model, and noting that $\delta^{p,p'}_{b',f}=0$, we
see that each sequence comprising two scoring vertices followed
by a non-scoring vertex is present amongst the ten configurations
prior to a move, except for the case depicted in Fig.~\ref{NotPFig}
and its up-down reflection.

\noindent
Only in these cases, where the 0th and 1st segments are scoring and
the first two segments are in the same direction, do we {\em not}
refer to the adjacent pair of scoring vertices as a particle.

Also note that when $p'>2p$ and $\delta^{p,p'}_{a',e}=
\delta^{p,p'}_{b',f}=0$, each
sequence of a non-scoring vertex followed by two scoring vertices
appears amongst the ten configurations that result from a move.
In such cases, the move may thus be reversed.

\subsection{The $\B_3$-transform}

Since in each of the moves described in Section \ref{PartMovesSec},
a pair of scoring vertices shifts to the right by one step,
we see that a succession of such 
moves is possible until the pair is followed by another scoring 
vertex. If this itself is followed by yet another scoring vertex,
we forbid further movement. However, if it is followed by a non-scoring 
vertex, further movement is allowed after considering the latter two 
of the three consecutive scoring vertices to be the particle (instead 
of the first two).

As in Section \ref{InsertSec}, let $h^{(k)}$ be a path resulting from
a $\B_2(k)$-transform
acting on a path that itself is the image of a $\B_1$ transform.
We now consider moving the $k$ particles that have been inserted.

\begin{lemma}\label{GaussLem} Let $\delta^{p,p'}_{b',f}=0$.
There is a bijection between the set of paths obtained by moving the
particles in $h^{(k)}$ and ${\mathcal Y}(k,m)$, where $m=m(h^{(k)})$.
This bijection is such that if $\lambda\in{\mathcal Y}(k,m)$ is the
bijective image of a particular $h$ then
$$
\mwt(h)=\mwt(h^{(k)})+\wt(\lambda).
$$
Additionally, $L(h)=L(h^{(k)})$ and $m(h)=m(h^{(k)})$.
\end{lemma}

\Proof
Since each particle moves by traversing a non-scoring vertex,
and there are $m$ of these to the right of the rightmost particle 
in $h^{(k)}$, and there are no consecutive scoring vertices to its 
right, this particle can make $\lambda_1$ moves to the right, with 
$0\le\lambda_1\le m$. Similarly, the next rightmost particle can 
make $\lambda_2$ moves to the right with $0\le\lambda_2\le\lambda_1$.
Here, the upper restriction arises because the two scoring vertices 
would then be adjacent to those of the first particle. Continuing in 
this way, we obtain that all possible final positions of the particles 
are indexed by $\lambda=(\lambda_1,\lambda_2,\ldots,\lambda_k)$ with 
$m\ge\lambda_1\ge\lambda_2\ge\cdots\ge\lambda_k\ge0$, that is, by 
partitions of at most $k$ parts with no part exceeding $m$. Moreover, 
since by Lemma \ref{WtChngLem} the weight increases by one for each 
move, the weight increase after the sequence of moves specified by 
a particular $\lambda$ is equal to $\wt(\lambda)$.
The final statement also follows from Lemma \ref{WtChngLem}.
\cqfd
\medskip

\noindent We say that a path obtained by moving the particles
in $h^{(k)}$ according to the partition $\lambda$ has been obtained
by the action of a $\B_3(\lambda)$-transform.

Having defined $\B_1$, $\B_2(k)$ for $k\ge0$ and $\B_3(\lambda)$
for $\lambda$ a partition with at most $k$ parts, we now define a 
$\B(k,\lambda)$-transform as the composition 
$\B(k,\lambda)=\B_3(\lambda)\circ\B_2(k)\circ\B_1$.

\begin{lemma}\label{BresLem}
Let $h' \in \P^{p,p'+p}_{a',b',e,f}(L')$ be obtained from
$h \in \P^{p,p'}_{a,b,e,f}(L)$ by the action of the $\B(k,\lambda)$-transform.
If $\pi=\pi(h)$ and $m=m(h)$ then:
\begin{displaymath}
\begin{array}{l}
\bullet\quad L'=
  \left\{
    \begin{array}{ll}
  2L-m+2k+2 \quad &
            \mbox{if } \pi=1 \mbox{ and } e=d,\\[1.5mm]
  2L-m+2k \quad & \mbox{otherwise};
    \end{array} \right.
\\[1mm]
\bullet\quad m(h')=L;\\[1mm]
\bullet\quad \mwt(h')=\mwt(h) + \frac{1}{4}\left( (L'-L)^2
 - \beta^2 \right) + \wt(\lambda),
\end{array}
\end{displaymath}
where $\beta=\beta^{p,p'}_{a,b,e,f}$.
\end{lemma}

\Proof These results follow immediately from
Lemmas \ref{BHashLem}, \ref{WtShiftLem} and \ref{GaussLem}.
\cqfd
\medskip

\begin{note}\label{Start2Note}
Since particle insertion and the particle moves don't change the
startpoint, endpoint or value $e(h)$ or $f(h)$ of a path $h$,
then in view of Lemma \ref{StartPtLem} and Corollary \ref{EndPtCor},
we see that the action of a $\B$-transform on
$h\in\P^{p,p'}_{a,b,e,f}(L)$
yields a path $h'\in\P^{p,p'+p}_{a',b',e,f}(L')$,
where 
$a'=a+\lfloor ap/p'\rfloor+e$, $b'=b+\lfloor bp/p'\rfloor+f$,
and $\delta^{p,p'+p}_{a',e}=\delta^{p,p'+p}_{b',f}=0$.
\end{note}

\goodbreak
\subsection{Particle content of a path}\label{ParticleSec}
\nobreak

Again restrict to the case
$p'>2p$ so that the $(p,p')$-model has no two neighbouring odd bands,
and let $h'\in\P^{p,p'}_{a',b',e,f}(L')$.
In the following lemma, we once more restrict to the cases for which
$\delta^{p,p'}_{a,e}=\delta^{p,p'}_{b,f}=0$,
and thus only consider the cases for which
the pre-segment and the post-segment of $h'$ lie in even bands.

\begin{lemma}\label{UniqueLem}
For $1\le p<p'$ with $p'>2p$, let $1\le a',b'<p'$ and $e,f\in\{0,1\}$,
with $\delta^{p,p'}_{a',e}=\delta^{p,p'}_{b',f}=0$.
If $h'\in\P^{p,p'}_{a',b',e,f}(L')$,
then there is a unique triple $(h,k,\lambda)$
where $h\in\P^{p,p'-p}_{a,b,e,f}(L)$ for some $a,b,L$, such that
the action of a $\B(k,\lambda)$-transform on $h$ results in $h'$.
\end{lemma}

\Proof This is proved by reversing the constructions described
in the previous sections. Locate the leftmost pair of consecutive 
scoring vertices in $h'$, and move them leftward by reversing the 
particle moves, until they occupy the $0$th and $1$st positions.
This is possible in all cases when
$\delta^{p,p'}_{a',e}=\delta^{p,p'}_{b',f}=0$.
Now ignoring these two vertices, do the same with the next leftmost 
pair of consecutive scoring vertices, moving them leftward until 
they occupy the third and fourth positions. Continue in this way 
until all consecutive scoring vertices occupy the leftmost positions 
of the path.
Denote this path by $h^{(\cdot)}$.
At the leftmost end of $h^{(\cdot)}$, there will be a number of even
segments (possibly zero) alternating in direction.
Let this number be $2k$ or $2k+1$ according to whether is it even or odd.
Clearly $h'$ results from $h^{(\cdot)}$ 
by a $\B_3(\lambda)$-transform for a particular $\lambda$ with at most 
$k$ parts.

Removing the first $2k$ segments of $h^{(\cdot)}$ yields a path
$h^{(0)}\in\P^{p,p'}_{a',b',e,f}$. This path thus has no
two consecutive scoring vertices, except possibly at the $0$th
and $1$st positions, and then only if the first vertex is a straight
vertex (as in Fig.~\ref{NotPFig}).
Moreover, $h^{(k)}$ arises by the action of
a $\B_2(k)$-transform on $h^{(0)}$.

Ignoring for the moment the case where there are scoring vertices at the
$0$th and $1$st positions, $h^{(0)}$ has by construction no pair of
consecutive scoring vertices.
Therefore, beyond the $0$th vertex,
we may remove a non-scoring vertex before every scoring
vertex to obtain a path $h\in\P^{p,p'-p}_{a,b,e,f}(L)$ for
some $a,b,L$, from which $h^{(0)}$ arises by the action of a
$\B_1$-transform.

On examining the third case depicted in Table~\ref{BatStartTable},
we see that the case where $h^{(0)}$ has a pair of scoring vertices
at the $0$th and $1$st positions, arises similarly from a particular
$h\in\P^{p,p'-p}_{a,b,e,f}(L)$ for some $a,b,L$.
The lemma is then proved.
\cqfd
\medskip

\noindent
The value of $k$ obtained above will be referred to as the
particle content of $h'$.

\begin{lemma}\label{BijLem}
For $1\le p<p'$, let $1\le a,b<p'$ and $e,f\in\{0,1\}$,
with $\delta^{p,p'}_{a,e}=0$.
Set $a'=a+e+\lfloor ap/p'\rfloor$ and $b'=b+f+\lfloor bp/p'\rfloor$.
Fix $m_0,m_1\ge0$.
Then the map $(h,k,\lambda)\mapsto h'$ effected by the action of
a $\B(k,\lambda)$-transform on $h$, is a bijection between
$\bigcup_{k}
\P^{p,p'}_{a,b,e,f}(m_1,2k+2m_1-m_0)\times{\mathcal Y}(k,m_1)$
and $\P^{p,p'+p}_{a',b',e,f}(m_0,m_1)$.
Moreover,
\begin{displaymath}
\mwt(h')=\mwt(h)+\frac{1}{4}\left( (m_0-m_1)^2 - \beta^2 \right)+\wt(\lambda),
\end{displaymath}
where $\beta=\beta^{p,p'}_{a,b,e,f}$.
\end{lemma}

\Proof Given $h\in\P^{p,p'}_{a,b,e,f}(m_1,m)$, let
$h'$ be the result of a $\B(k,\lambda)$-transform on $h$.

Since $\delta^{p,p'}_{a,e}=0$ so that
$\lfloor (a+(-1)^e)p/p')\rfloor=\lfloor ap/p'\rfloor$,
it follows that if $\pi(h)=1$ then $e(h)\ne d(h)$.
Then, with $m=2m_1+2k-m_0$,
we obtain
$h'\in\P^{p,p'+p}_{a',b',e,f}(m_0,m_1)$
via Lemma \ref{BresLem}. 

Lemma \ref{StartPtLem} shows that
$\delta^{p,p'+p}_{a',e}=\delta^{p,p'+p}_{b',f}=0$.
Thereupon, Lemma \ref{UniqueLem} shows that each
$h'\in\P^{p,p'+p}_{a',b',e,f}(m_0,m_1)$ arises from a unique triple
$(h,k,\lambda)$, with $h\in\P^{p,p'}_{a,b,e,f}(m_1,m)$ for some $m$.
The bijection then follows.

The expression for $\mwt(h')$ also results from Lemma \ref{BresLem}.
\cqfd
\medskip

Note that the above lemma excludes consideration of the case for which
$\delta^{p,p'}_{a,e}=1$.
In fact, similar results fail in that case.
Nonetheless, it is necessary to tackle the $\delta^{p,p'}_{a,e}=1$ case
for a restricted set of paths in the more general analysis of
\cite{foda-welsh-big}.

\begin{corollary}\label{BijCor}
For $1\le p<p'$, let $1\le a,b<p'$ and $e,f\in\{0,1\}$,
with $\delta^{p,p'}_{a,e}=0$.
Set $a'=a+e+\lfloor ap/p'\rfloor$ and $b'=b+f+\lfloor bp/p'\rfloor$.
Fix $m_0,m_1\ge0$.
Then
\begin{displaymath}
\mchi^{p,p'+p}_{a',b',e,f}(m_0,m_1)=
q^{\frac{1}{4}\left( (m_0-m_1)^2 - \beta^2 \right)}
\sum_{\scriptstyle m\equiv m_0\atop
  \scriptstyle\strut(\mbox{\scriptsize\rm mod}\,2)}
\left[{\frac{1}{2}(m_0+m)\atop m_1}\right]_q
\mchi^{p,p'}_{a,b,e,f}(m_1,m),
\end{displaymath}
where $\beta=\beta^{p,p'}_{a,b,e,f}$.
\end{corollary}

\Proof Apart from the case where $m_1=0$ and $e\ne f$,
this follows immediately from Lemma \ref{BijLem} on setting
$m=2m_1+2k-m_0$,
once it is noted, via Lemma \ref{PartitionGenLem},
that $\left[{k+m_1\atop m_1}\right]_q$ is the
generating function for ${\mathcal Y}(k,m_1)$.

For the case $m_1=0$ and $e\ne f$, both sides are zero
unless $a=b$ and $m_0$ is odd.
In this case, $\P^{p,p'+p}_{a',b',e,f}(m_0,0)$ has precisely
one element $h$ for which
(via the same calculation as in the proof of \ref{SeedLem})
$\mwt(h)=\frac14(m_0^2-1)$.
Thus the two sides are also equal in this case.
\cqfd
\medskip

\section{The $\D$-transform}\label{DTranSec}

The {\it $\D$-transform}
is defined to act on each $h \in \P^{p,p'}_{a,b,e,f}(L)$ to yield
a path $\hat h\in\P^{p'-p,p'}_{a,b,1-e,1-f}(L)$ with exactly the same
sequence of integer heights, i.e., $\hat h_i=h_i$ for $0\le i\le L$.
Note that, by definition, $e(\hat h)=1-e(h)$ and $f(\hat h)=1-f(h)$.

Since the band structure of 
the $(p'-p,p')$-model is obtained from that of the $(p,p')$-model
simply by replacing odd bands by even bands and vice-versa,
then, ignoring the vertex at $i=0$,
each scoring vertex maps to a non-scoring vertex and vice-versa.
That $e(h)$ and $e(\hat h)$ differ implies that the vertex
at $i=0$ is both scoring or both non-scoring in $h$ and $\hat h$.

\begin{lemma}\label{DresLem}
Let $\hat h \in \P^{p'-p,p'}_{a,b,1-e,1-f}(L)$ be obtained from
$h \in \P^{p,p'}_{a,b,e,f}(L)$ by the action of the $\D$-transform.
Then $\pi(\hat h)=1-\pi(h)$.
Moreover, if $m=m(h)$ then:
\begin{displaymath}
\begin{array}{l}
\bullet\quad L(\hat h)=L;\\[2mm]
\bullet\quad m(\hat h)=
  \left\{
    \begin{array}{ll}
  L-m \quad &
            \mbox{if } e+d+\pi(h)\equiv0\,(\mod2),\\[1.5mm]
  L-m+2 \quad & \mbox{if } e+d+\pi(h)\not\equiv0\,(\mod2);
    \end{array} \right.
\\[2mm]
\bullet\quad \mwt(\hat h)=\frac{1}{4}
   \left( L^2 - \alpha(h)^2 \right) - \mwt(h).
\end{array}
\end{displaymath}
\end{lemma}

\Proof
Let $h$ have striking sequence
$\left({a_1\atop b_1}\:{a_2\atop b_2}\:{a_3\atop b_3}\:
 {\cdots\atop\cdots}\:{a_l\atop b_l} \right)^{(e,f,d)}$.
Since, beyond the zeroth vertex, the $\D$-transform
exchanges scoring vertices for non-scoring vertices and vice-versa,
it follows that the striking sequence for $\hat h$ is
$\left({b_1\atop a_1}\:{b_2\atop a_2}\:{b_3\atop a_3}\:
 {\cdots\atop\cdots}\:{b_l\atop a_l} \right)^{(1-e,f,d)}$.
It is immediate that $L(\hat h)=L$, $\pi(\hat h)=1-\pi(h)$,
$e(\hat h)=1-e(h)$ and $d(\hat h)=d(h)$.
Then
$m(\hat h)=(e(\hat h)+d(\hat h)+\pi(\hat h))\,\mod2+\sum_{i=1}^l b_i
=(e+d+\pi(h))\,\mod2+L-\sum_{i=1}^l a_i
=2((e+d+\pi(h))\,\mod2)+L-m(h)$.

Now let $w_i=a_i+b_i$ for $1\le i\le l$.
Then, using Lemma \ref{WtHashLem}, we obtain
\begin{eqnarray*}
\mwt(h)+\mwt(\hat h)&=& 
\sum_{i=1}^l b_i(w_{i-1}+w_{i-3}+\cdots+w_{1+i\bmod2})\\
&&\qquad
+ \sum_{i=1}^l a_i(w_{i-1}+w_{i-3}+\cdots+w_{1+i\bmod2})\\[0.5mm]
&=&
\sum_{i=1}^l w_i(w_{i-1}+w_{i-3}+\cdots+w_{1+i\bmod2})\\[0.5mm]
&=&
(w_1+w_3+w_5+\cdots)(w_2+w_4+w_6+\cdots).
\end{eqnarray*}
The lemma then follows because
$(w_1+w_3+w_5+\cdots)+(w_2+w_4+w_6+\cdots)=L$ and
$(w_1+w_3+w_5+\cdots)-(w_2+w_4+w_6+\cdots)=\pm\alpha(h)$.
\cqfd
\medskip


\begin{lemma}\label{DParamLem}
Let $1\le p<p'$ with $p$ co-prime to $p'$ and $1\le a<p'$.
Then $\lfloor a(p'-p)/p'\rfloor=a-1-\lfloor ap/p'\rfloor$.

If, in addition, $a$ is interfacial in the $(p,p')$-model
and $\delta^{p,p'}_{a,e}=0$ then
$a$ is interfacial in the $(p'-p,p')$-model
and $\delta^{p'-p,p'}_{a,1-e}=0$.
\end{lemma}

\Proof Since $p$ and $p'$ are co-prime, $\lfloor ap/p'\rfloor<ap/p'$.
Hence $\lfloor ap/p'\rfloor+\lfloor a(p'-p)/p'\rfloor=a-1$.

Since the $(p,p')$-model differs from the $(p'-p,p')$-model only in
that corresponding bands are of the opposite parity, $a$ being interfacial
in one model implies that it also is in the other.
The final part then follows immediately.
\cqfd
\medskip

\begin{corollary}\label{DParamCor}
If $1\le p<p'$ with $p$ co-prime to $p'$,
$1\le a,b<p'$ and $e,f\in\{0,1\}$ then
$\alpha^{p'-p,p'}_{a,b}=\alpha^{p,p'}_{a,b}$ and
$\beta^{p'-p,p'}_{a,b,1-e,1-f}+\beta^{p,p'}_{a,b,e,f}=\alpha^{p,p'}_{a,b}$.
\end{corollary}

\Proof Lemma \ref{DParamLem} gives
$\lfloor ap/p'\rfloor+\lfloor a(p'-p)/p'\rfloor=a-1$ and
likewise, $\lfloor bp/p'\rfloor+\lfloor b(p'-p)/p'\rfloor=b-1$.
The required results then follow immediately.
\cqfd
\medskip


\subsection{The $\B\D$-pair}

It will often be convenient to consider the combined action
of a $\D$-transform followed immediately by a $\B$-transform.
Such a pair will naturally be referred to as a $\B\D$-transform
and maps a path $h\in\P^{p'-p,p'}_{a,b,e,f}(L)$ to a path
$h'\in\P^{p,p'+p}_{a',b',1-e,1-f}(L')$, where $a',b',L'$
are determined by our previous results.

In what follows, the $\B\D$-transform will always follow a $\B$-transform.
Thus we restrict consideration to where $2(p'-p)<p'$.

\begin{lemma}\label{BDresLem}
With $p'<2p$, let $h\in{\P}^{p'-p,p'}_{a,b,e,f}(L)$. Let
$h'\in\P^{p,p'+p}_{a',b',1-e,1-f}(L')$ result from the action
of a $\D$-transform on $h$, followed by a $\B(k,\lambda)$-transform.
Then:
\begin{displaymath}
\begin{array}{l}
\bullet\quad L'=
  \left\{
    \begin{array}{ll}
  L+m(h)+2k-2 \quad &
            \mbox{if } \pi(h)=1 \mbox{ and } e=d(h),\\[1.5mm]
  L+m(h)+2k \quad & \mbox{otherwise};
    \end{array} \right.
\\[2mm]
\bullet\quad m(h')=L;\\[2mm]
\bullet\quad \mwt(h')=\frac{1}{4}\left( L^2 + (L'-L)^2
- \alpha^2 - \beta^2 \right) + \wt(\lambda) - \mwt(h),
\end{array}
\end{displaymath}
where $\alpha=\alpha^{p,p'}_{a,b}$ and $\beta=\beta^{p,p'}_{a,b,1-e,1-f}$.
\end{lemma}

\Proof Let $\hat h$ result from the action of the $\D$-transform
on $h$, and let $d=d(h)$, $\pi=\pi(h)$, $\hat e=e(\hat h)$
$\hat d=d(\hat h)$, $\hat\pi=\pi(\hat h)$. Then we immediately
have $\hat d=d$, $\hat e=1-e$, and $\hat\pi=1-\pi$.

In the case where $\pi=0$ and $e\ne d$, we then have,
using Lemmas \ref{BresLem} and \ref{DresLem},
$L'=2L(\hat h)-m(\hat h)+2k+2=2L-(L-m(h)+2)+2k+2=L+m(h)+2k$.

In the case where $\pi=1$ and $e=d$, we then have,
using Lemmas \ref{BresLem} and \ref{DresLem},
$L'=2L(\hat h)-m(\hat h)+2k=2L-(L-m(h)+2)+2k=L+m(h)+2k-2$.

In the other cases, $e+d+\pi\equiv0\,(\mod2)$ and
so $\hat e+\hat d+\hat\pi\equiv0\,(\mod2)$.
Lemmas \ref{BresLem} and \ref{DresLem} yield
$L'=2L(\hat h)-m(\hat h)+2k=2L-(L-m(h))+2k=L+m(h)+2k$.

The expressions for $m(h')$
and $\mwt(h')$
also follow immediately from Lemmas \ref{BresLem} and \ref{DresLem}.
\cqfd
\medskip

We now obtain analogues of Lemma \ref{BijLem} and
Corollary \ref{BijCor} which combine
the $\D$-transform with the $\B$-transform.
As above, we restrict to where $p'<2p$.

\begin{lemma}\label{DijLem}
For $1\le p<p'<2p$, let $1\le a,b<p'$ and $e,f\in\{0,1\}$,
with $\delta^{p'-p,p'}_{a,e}=0$.
Set $a'=a+1-e+\lfloor ap/p'\rfloor$ and  $b'=b+1-f+\lfloor bp/p'\rfloor$.
Fix $m_0,m_1\ge0$.
Then the map $(h,k,\lambda)\mapsto h'$ effected by the action of
a $\D$-transform on $h$ followed by a $\B(k,\lambda)$-transform,
is a bijection between
$\bigcup_{k}
\P^{p'-p,p'}_{a,b,e,f}(m_1,m_0-m_1-2k)\times{\mathcal Y}(k,m_1)$
and $\P^{p,p'+p}_{a',b',1-e,1-f}(m_0,m_1)$.
Moreover,
\begin{displaymath}
\mwt(h')=\frac{1}{4}\left( m_1^2 + (m_0-m_1)^2
                   - \alpha^2 - \beta^2 \right) + \wt(\lambda) - \mwt(h),
\end{displaymath}
where $\alpha=\alpha^{p,p'}_{a,b}$ and $\beta=\beta^{p,p'}_{a,b,1-e,1-f}$.
\end{lemma}

\Proof Given $h\in\P^{p'-p,p'}_{a,b,e,f}(m_1,m)$, let
$\hat h$ result from the action of a $\D$-transform on $h$,
and let $h'$ be the result of a $\B(k,\lambda)$-transform on $\hat h$.

Since $\delta^{p'-p,p'}_{a,e}=0$ so that
$\lfloor (a+(-1)^e)(p'-p)/p')\rfloor=\lfloor a(p'-p)/p'\rfloor$,
it follows that if $\pi(h)=1$ then $e(h)\ne d(h)$.
Then, for $m=m_0-m_1-2k$, we obtain
$h'\in\P^{p,p'+p}_{a',b',1-e,1-f}(m_0,m_1)$
via Lemma \ref{BDresLem}.

Lemma \ref{StartPtLem} gives
$\delta^{p,p'+p}_{a',1-e}=\delta^{p,p'+p}_{b',1-f}=0$.
Lemma \ref{UniqueLem} then shows that
for arbitrary $h'\in\P^{p,p'+p}_{a',b',1-e,1-f}(m_0,m_1)$,
there is a unique triple $(\hat h,k,\lambda)$,
with $\hat h\in\P^{p,p'}_{a,b,1-e,1-f}(m_1,m')$ for some $m'$,
such that
the action of the $\B(k,\lambda)$-transform on $\hat h$ yields $h'$.
Then, via the $\D$-transform, we obtain a unique
$h\in\P^{p'-p,p'}_{a,b,e,f}(m_1,m'')$, for some $m''$.
The bijection then follows.

The expression for $\mwt(h)$ also results from Lemma \ref{BDresLem}.
\cqfd
\medskip

Note that the above lemma excludes the case for which
$\delta^{p'-p,p'}_{a,e}=1$.
Once more, similar results fail in that case.

\begin{corollary}\label{DijCor}
For $1\le p<p'<2p$, let $1\le a,b<p'$ and $e,f\in\{0,1\}$,
with $\delta^{p'-p,p'}_{a,e}=0$.
Set $a'=a+1-e+\lfloor ap/p'\rfloor$ and  $b'=b+1-f+\lfloor bp/p'\rfloor$.
Fix $m_0,m_1\ge0$.
Then
\begin{displaymath}
\begin{array}{l}
\displaystyle
\mchi^{p,p'+p}_{a',b',1-e,1-f}(m_0,m_1;q)=\\[2mm]
\hskip3mm
\displaystyle
q^{\frac{1}{4}\left( m_1^2 + (m_0-m_1)^2 - \alpha^2 - \beta^2 \right)}
\hskip-4mm
\sum_{\scriptstyle m\equiv m_0-m_1\atop
  \scriptstyle\strut(\mbox{\scriptsize\rm mod}\,2)}
\left[{\frac{1}{2}(m_0+m_1-m)\atop m_1}\right]_q
\mchi^{p'-p,p'}_{a,b,e,f}(m_1,m;q^{-1}),\!\!\!
\end{array}
\end{displaymath}
where $\alpha=\alpha^{p,p'}_{a,b}$ and $\beta=\beta^{p,p'}_{a,b,1-e,1-f}$.
\end{corollary}

\Proof Apart from the case where $m_1=0$ and $e\ne f$,
this follows immediately from Lemma \ref{DijLem} on setting
$m=m_0-m_1-2k$,
once it is noted, via Lemma \ref{PartitionGenLem},
that $\left[{k+m_1\atop m_1}\right]_q$ is the
generating function for ${\mathcal Y}(k,m_1)$.

The case $m_1=0$ and $e\ne f$ is dealt with exactly as in the proof
of Corollary \ref{BijCor}.
\cqfd
\medskip

\begin{lemma}\label{BDParamLem}
Let $1\le p<p'<2p$ with $p$ co-prime to $p'$,
$1\le a,b<p'$ and $e,f\in\{0,1\}$ and set
$a'=a+1-e+\lfloor ap/p'\rfloor$ and  $b'=b+1-f+\lfloor bp/p'\rfloor$.
Then
$\lfloor a'p/(p'+p)\rfloor=a-1-\lfloor a(p'-p)/p'\rfloor$
and
$\lfloor b'p/(p'+p)\rfloor=b-1-\lfloor b(p'-p)/p'\rfloor$.
In addition,
$\alpha^{p,p'+p}_{a',b'}=2\alpha^{p'-p,p'}_{a,b}-\beta^{p'-p,p'}_{a,b,e,f}$
and
$\beta^{p,p'+p}_{a',b',1-e,1-f}
  =\alpha^{p'-p,p'}_{a,b}-\beta^{p'-p,p'}_{a,b,e,f}$.
\end{lemma}

\Proof By Lemma \ref{DParamLem} and Corollary \ref{DParamCor},
$\lfloor ap/p'\rfloor=a-1-\lfloor a(p'-p)/p'\rfloor$,
$\lfloor bp/p'\rfloor=b-1-\lfloor b(p'-p)/p'\rfloor$,
$\alpha^{p,p'}_{a,b}=\alpha^{p'-p,p'}_{a,b}$ and
$\beta^{p,p'}_{a,b,1-e,1-f}=\alpha^{p'-p,p'}_{a,b}-\beta^{p'-p,p'}_{a,b,e,f}$.
The current lemma then follows immediately from Lemma \ref{ParamLem}.
\cqfd
\medskip


\section{The structure of the $(p,p')$-model}

\subsection{Continued fractions}\label{ContFSec}

If $p'$ and $p$ are positive co-prime integers and
$$
\frac{p'}{p}=
{c_0+\frac{\displaystyle\strut 1}{\displaystyle c_1 +
\frac{\displaystyle\strut 1}{\displaystyle c_2 +
\frac{\displaystyle\strut 1}{\frac{\lower-5pt\hbox{$\vdots$}}
{\displaystyle\strut c_{n-1} +
\frac{\displaystyle\strut 1}{\displaystyle c_n}}}}}}
$$
with $c_0\ge0$, $c_i\ge1$ for $0<i<n$, and $c_n\ge2$,
then $(c_0,c_1,c_2,\ldots,c_n)$ is said to be the
{\it continued fraction} for $p'/p$.

We refer to $n$ as the {\em height} of $p'/p$.
We set $t=c_0+c_1+\cdots+c_n-2$ and refer to it as the {\em rank} of $p'/p$.
The height and rank of ${\P}^{p,p'}_{a,b,c}(L)$ are then defined
to be equal to those of $p'/p$.


For $0\le k\le n+1$, we also define
\begin{equation}\label{ZoneEq}
t_k=-1+\sum_{i=0}^{k-1} c_i.
\end{equation}
Then $t_{n+1}=t+1$ and $t_n\le t-1$.
We say that the index $j$ with $0\le j\le t_{n+1}$ is in zone $k$
if $t_{k}<j\le t_{k+1}$.
We then write $k=\zeta(j)$.
Note that there are $n+1$ zones and that for 
$0\le k\le n$, zone $k$ contains $c_k$ indices.

\goodbreak
\subsection{The Takahashi and string lengths}\label{TakSec}
\nobreak

Given positive co-prime integers $p$ and $p'$ with $p'/p$ having
rank $t$, define
the set $\{\kappa_i\}_{i=0}^t$ of
{\em Takahashi lengths},
the set $\{\tilde\kappa_i\}_{i=0}^t$ of
{\em truncated Takahashi lengths},
and the set $\{l_i\}_{i=0}^t$
of {\em string lengths} as follows.
First define $y_k$ and $z_k$ for
$-1\le k\le n+1$ by:
\begin{displaymath}
\begin{array}{lcllcl}
y_{-1}&=&0;&z_{-1}&=&1;\\
y_0&=&1;&z_0&=&0;\\
y_k&=&c_{k-1}y_{k-1}+y_{k-2};&z_k&=&c_{k-1}z_{k-1}+z_{k-2},
\quad (1\le k\le n+1).
\end{array}
\end{displaymath}
Now for $t_k<j\le t_{k+1}$ and $0\le k\le n$, set
\begin{eqnarray*}
\kappa_j&=&y_{k-1}+(j-t_k)y_k;\\
\tilde\kappa_j&=&z_{k-1}+(j-t_k)z_k;\\
l_j&=&y_{k-1}+(j-t_k-1)y_k.
\end{eqnarray*}
Note that $\kappa_{j}=l_{j+1}$ unless $j=t_k$ for some $k$,
in which case $\kappa_{t_k}=y_k$ and $l_{t_k+1}=y_{k-1}$.
We define ${\mathcal T}=\{\kappa_i\}_{i=0}^{t-1}$
and ${\mathcal T}^\prime=\{p'-\kappa_i\}_{i=0}^{t-1}$.
(We don't include $\kappa_t$ in the former since it is
present in the latter.)
Then, for $n>0$, $\mathcal T\cap\mathcal T'=\emptyset$.\footnote
{In fact, when $n=0$,
$\mathcal T\cap\mathcal T'=\{2,3,\ldots,p'-2\}$.
Then, if $2\le a\le p'-2$, different fermionic expressions for
$\P^{p,p'}_{a,b,c}(L)$ arise by considering either
$a\in\mathcal T$ or $a\in\mathcal T'$.
The same holds for $2\le b\le p'-2$.
This $n=0$ case was fully examined in \cite{foda-welsh}.}


For example, in the case $p'=38$, $p=11$, for which the continued
fraction is $(3,2,5)$, so that $n=2$,
$(t_1,t_2,t_3)=(2,4,9)$ and $t=8$.
We then obtain:
\begin{displaymath}
\begin{array}{l}
(y_{-1},y_0,y_1,y_2,y_3)=(0,1,3,7,38),\\
(z_{-1},z_0,z_1,z_2,z_3)=(1,0,1,2,11),\\
(\kappa_0,\kappa_1,\kappa_2,\kappa_3,\kappa_4,\kappa_5,
\kappa_6,\kappa_7)=(1,2,3,4,7,10,17,24),\\
(l_1,l_2,l_3,l_4,l_5,l_6,l_7,l_8)=(1,2,1,4,3,10,17,24),\\
(\tkappa_0,\tkappa_1,\tkappa_2,\tkappa_3,\tkappa_4,\tkappa_5,
\tkappa_6,\tkappa_7)=(1,1,1,1,2,3,5,7).
\end{array}
\end{displaymath}

An induction argument readily establishes
that if $1\le k\le n+1$, then $y_kz_{k-1}-y_{k-1}z_k=(-1)^k$,
that $y_k$ is co-prime to $z_k$, and that $y_k/z_k$
has continued fraction $(c_0,c_1,\ldots,c_{k-1})$.
Thus, in particular, $y_{n+1}=p'$ and $z_{n+1}=p$.

\newpage

\goodbreak

\section{Segmenting the model}\label{SegmentSec}
\nobreak

\subsection{Model comparisons}\label{ModComSec}
\nobreak

Here, we relate the parameters associated with the $(p,p')$-model
for which the continued fraction is $(c_0,c_1,\ldots,c_n)$ to those
associated with certain \lq simpler\rq\ models.
In particular, if $c_0>1$, we compare them with those associated
with the $(p,p'-p)$-model and, if $c_0=1$, we compare them with
those associated with the $(p'-p,p')$-model.

In the following two lemmas, the parameters associated with those
simpler models will be primed to distinguish them from those
associated with the $(p,p')$-model.
In particular if $c_0>1$, $(p'-p)/p$ has
continued fraction $(c_0-1,c_1,\ldots,c_n)$, so that in this case,
$t'=t-1$, $n'=n$ and $t_k'=t_k-1$ for $1\le k\le n$.
If $c_0=1$, $p'/(p'-p)$ has continued fraction
$(c_1+1,c_2,\ldots,c_n)$, so that in this case,
$t'=t$, $n'=n-1$ and $t_k'=t_{k+1}$ for $1\le k\le n'$.

\begin{lemma}\label{BmodelLem}
Let $c_0>1$. For $1\le k\le n$ and $0\le j\le t$, let $y_k$, $z_k$,
$\kappa_j$ and $\tilde\kappa_j$ be the parameters associated with
the $(p,p')$-model as defined in Section \ref{TakSec}.
For $1\le k\le n$ and $0\le j\le t'$, let $y_k'$, $z_k'$,
$\kappa_j'$ and $\tilde\kappa_j'$ be the corresponding parameters
for the $(p,p'-p)$-model.
Then:
\begin{itemize}
\item $y_k=y_k'+z_k'\quad (0\le k\le n)$;
\item $z_k=z_k'\quad (0\le k\le n)$;
\item $\kappa_j=\kappa_{j-1}'+\tilde\kappa_{j-1}' \quad (1\le j\le t)$;
\item $\tilde\kappa_j=\tilde\kappa_{j-1}' \quad (1\le j\le t)$.
\end{itemize}
\end{lemma}

\Proof This result is a straightforward consequence of the definitions.
\cqfd
\medskip

\begin{lemma}\label{DmodelLem}
Let $c_0=1$. For $1\le k\le n$ and $0\le j\le t$, let $y_k$, $z_k$,
$\kappa_j$ and $\tilde\kappa_j$ be the parameters associated with
the $(p,p')$-model as defined in Section \ref{TakSec}.
For $1\le k\le n'$ and $0\le j\le t$, let $y_k'$, $z_k'$,
$\kappa_j'$ and $\tilde\kappa_j'$ be the corresponding parameters
for the $(p'-p,p')$-model.
Then:
\begin{itemize}
\item $y_k=y_{k-1}'\quad (1\le k\le n)$;
\item $z_k=y_{k-1}'-z_{k-1}'\quad (1\le k\le n)$;
\item $\kappa_j=\kappa_j' \quad (1\le j\le t)$;
\item $\tilde\kappa_j=\kappa_j'-\tilde\kappa_j' \quad (1\le j\le t)$.
\end{itemize}
\end{lemma}

\Proof Again, this result is a straightforward consequence of the definitions.
\cqfd
\medskip

\begin{lemma}\label{TakBandLem}
If $t_1\le j\le t$ then\footnote{We use the notation
$\delta^{(2)}_{i,j}=1$ if $i\equiv j\,(\mod2)$
and $\delta^{(2)}_{i,j}=0$ if $i\not\equiv j\,(\mod2)$.}
$$
\left\lfloor \frac{\tkappa_{j}p'}{p}\right\rfloor
  =\kappa_{j}-\delta^{(2)}_{\zeta(j),1},
$$
and if $0\le j\le t$ then
$$
\left\lfloor \frac{\kappa_{j}p}{p'}\right\rfloor
  =\tkappa_{j}-\delta^{(2)}_{\zeta(j),0}.
$$
\end{lemma}

\Proof We prove the first of these two results by induction on
the sum of the height and rank of $p'/p$.
Since $\kappa_{t_1}=c_0$ and $\tkappa_{t_1}=1$ and $\zeta(t_1)=0$,
the required result always holds for the case $j=t_1$.
In particular, it certainly holds in the case where the sum of the height
and rank of $p'/p$ is at most 2.

Now assume that the first part holds in the case that sum of height and rank
is $n+t-1$, and consider the case where $p'/p$ has height $n$ and rank $t$.
First assume that $p'>2p$.
For $j\ge t_1$, the induction hypothesis implies that
$
\kappa'_{j-1}-\delta^{(2)}_{\zeta^\prime(j-1),1}
<\tkappa'_{j-1}(p'-p)/p
<\kappa'_{j-1}-\delta^{(2)}_{\zeta^\prime(j-1),1}+1
$,
where the primed quantities pertain to the continued fraction of $(p'-p)/p$.
Using Lemma \ref{BmodelLem} and noting that
$\zeta'(j-1)=\zeta(j)$, readily yields
$
\kappa_{j}-\delta^{(2)}_{\zeta(j),1}
<\tkappa_{j}p'/p
<\kappa_{j}-\delta^{(2)}_{\zeta(j),1}+1
$.
This immediately gives the required result.

In the case $p'<2p$, first let $j\ge t_2$.
The induction hypothesis implies that
$
\kappa'_{j}-\delta^{(2)}_{\zeta^\prime(j),1}
<\tkappa'_{j}p'/(p'-p)
<\kappa'_{j}-\delta^{(2)}_{\zeta^\prime(j),1}+1
$,
where the primed quantities pertain to the continued fraction of $p'/(p'-p)$.
Using Lemma \ref{DmodelLem} and noting that
$\zeta'(j)=\zeta(j)-1$, readily yields
$
\kappa_{j}-\delta^{(2)}_{\zeta(j),1}(p'-p)/p
<\tkappa_{j}p'/p
<\kappa_{j}+(1-\delta^{(2)}_{\zeta(j),1})(p'-p)/p
$.
Since $(p'-p)/p<1$, this implies the required result.

When $p'<2p$, we have $c_0=1$ so that $t_1=0$ and $t_2=c_1$. 
Then $\tkappa_j=j$ for $t_1<j\le t_2$, whereupon
in view of the continued fraction expression for $p'/p$,
we immediately obtain $\lfloor\tkappa_j p'/p\rfloor=j=\kappa_j-1$,
as required.

The first part of the lemma then follows by induction.
For $t_1\le j\le t$, the second part readily follows from the first.
For $0\le j\le t_1\le t$, both sides are clearly equal to 0.
\cqfd
\medskip

If $t_1\le j\le t$, it follows from this result that, with $k$ such that
$t_k<j\le t_{k+1}$, the $\tkappa_j$th odd band in the $(p,p')$-model
lies between heights $\kappa_j-1$ and $\kappa_j$ when $k$ is odd,
and between heights $\kappa_j$ and $\kappa_j+1$ when $k$ is even.
Since there are no adjacent odd bands when $p'>2p$, it follows
that $\kappa_j$ is interfacial when $j\ge t_1$.
On switching the parity of each band, we then obtain in the
case $p'<2p$ that $\kappa_j$ is interfacial when $j\ge t_2$.

\begin{lemma}\label{SegmentLem}
If $1\le p<p'$ and $p$ is co-prime to $p'$,
then for $1\le s\le y_n-2$, the $s$th band of the $(p,p')$-model
is of the same parity as the $s$th band of the $(z_n,y_n)$-model.
\end{lemma}

\Proof We must establish that
$\lfloor sz_{n}/y_{n}\rfloor=\lfloor sz_{n+1}/y_{n+1}\rfloor$
for $1\le s\le y_n-1$.

With $s$ such that $1\le s<y_n$,
let $r=\lfloor sz_{n+1}/y_{n+1}\rfloor$.
Using $y_nz_{n+1}=y_{n+1}z_n+(-1)^n$ then yields:
\begin{displaymath}
ry_n-(-1)^n\frac{s}{y_{n+1}}\le sz_n
<(r+1)y_n-(-1)^n\frac{s}{y_{n+1}}.
\end{displaymath}
Since $1\le s<y_n<y_{n+1}$, the first inequality here
implies that $sz_n/y_n\ge r$.
For the same reasons, and noting that $sz_n/y_n$ is not integral,
the second inequality here implies that $sz_n/y_n<r+1$.
The lemma then follows.
\cqfd
\medskip

\noindent
This lemma shows that the $(z_n,y_n)$-model resides within
the $(p,p')$-model, between heights $1$ and $y_n-1$.
The up-down symmetry of the $(p,p')$-model then also implies that
the $(z_n,y_n)$-model also resides within the $(p,p')$-model,
between heights $p'-y_n+1$ and $p'-1$.

\subsection{Interfacial retention}\label{RetentionSec}

We now show that if $h$ attains an interfacial height,
then the path resulting from the action of a $\B$-transform on
$h$ attains the corresponding interfacial height.

\begin{lemma}\label{Inter1Lem}
Let $h\in\P^{p,p'}_{a,b,e,f}(L)$, and let
$h\in\P^{p,p'+p}_{a',b',e,f}(L')$ result from the action
of a $\B(k,\lambda)$-transform on $h$.
Let $s$ be interfacial in the $(p,p')$-model with $a\ne s\ne b$,
and set $r=\lfloor (s+1)p/p'\rfloor$. Then $s+r$ is interfacial
in the $(p,p'+p)$-model.

If $h_i=s$ for $0\le i\le L$ then $h'_j=s+r$ for some $j$ with
$0\le j\le L'$.
On the other hand, if $h'_j=s+r$ for $0\le j\le L'$ then
$h_i=s$ for some $i$ with $0\le i\le L$.
\end{lemma}

\Proof 
First note that $s$ borders the $r$th odd band in the $(p,p')$-model.
If $s$ is at the lower (resp.~upper) edge
of the $r$th odd band in the $(p,p')$-model then
$s+r$ is at the lower (resp.~upper) edge of the $r$th odd band
in the $(p,p'+p)$-model. In particular, this implies
that $s+r$ is interfacial in the $(p,p'+p)$-model.
Then note that in the $(p,p')$-model, there is at least one even band
between the two odd
bands on either side of $s$ (assume that there is an odd band
immediately above and immediately below the $(p,p')$-model grid
if necessary). Thus there are at least two even bands between the
two odd bands on either side of $s+r$ in the $(p,p'+p)$-model.

Let $h^{(0)}$ result from the action of the $\B_1$-transform on $h$.
The definition of this transform implies that if $h_i=s$ for
some $i$ then $h^{(0)}_j=s+r$ for some $j$ and vice-versa
(when $\delta^{p,p'}_{a,e}=1$ or $\delta^{p,p'}_{b,f}=1$,
this statement relies on $a\ne s\ne b$).

If $h^{(k)}$ results from the action of the $\B_2(k)$-transform
on $h^{(0)}$, then if $h^{(0)}_j=s+r$ for some $j$ then
$h^{(k)}_{j'}=s+r$ for some $j'$ and vice-versa 
(this statement relies on the two odd bands either side of $s+r$
having at least two even bands between them).

If $h'$ results from the action of the $\B_3(\lambda)$-transform
on $h^{(k)}$, then if $h^{(0)}_j=s+r$ for some $j$,
examination of the ten particle moves and edge-moves described in
Section \ref{PartMovesSec}, shows that
$h^{(k)}_{j'}=s+r$ for some $j'$ and vice-versa 
(this statement also relies on the two odd bands either side of $s+r$
having at least two even bands between them).
Combining these results proves the lemma.
\cqfd
\medskip

We also need the analogue of this result for the $\B\D$-transform.

\begin{lemma}\label{Inter2Lem}
Let $h\in\P^{p'-p,p'}_{a,b,e,f}(L)$ and let
$h'\in\P^{p,p'+p}_{a',b',1-e,1-f}(L')$ result from the action
of a $\D$-transform on $h$ followed by a $\B(k,\lambda)$-transform.
Let $s$ be interfacial in the $(p,p')$-model with $a\ne s\ne b$,
and set $r=\lfloor (s+1)p/p'\rfloor$. Then $s+r$ is interfacial
in the $(p,p'+p)$-model.

If $h_i=s$ for $0\le i\le L$ then $h'_j=s+r$ for some $j$ with
$0\le j\le L'$.
On the other hand, if $h'_j=s+r$ for $0\le j\le L'$ then
$h_i=s$ for some $i$ with $0\le i\le L$.
\end{lemma}

\Proof
This follows immediately from the above result after noting
that if $s$ is interfacial in the $(p'-p,p')$-model then it
is also in the $(p,p')$-model.
\cqfd
\medskip

A set $\mathcal S$ is said to be interfacial in the $(p,p')$-model
if each $s\in\mathcal S$ is interfacial in the $(p,p')$-model.
We now define $\P^{p,p'}_{a,b,e,f}(L,m)\{\mathcal S\}$
to be the subset of $\P^{p,p'}_{a,b,e,f}(L,m)$
comprising those paths $h$ for which for each $s\in\mathcal S$,
there exists $i$ with $0\le i\le L$ such that $h_i=s$.
The generating function for this set is
\begin{displaymath}
\mchi^{p,p'}_{a,b,e,f}(L;q)\{\mathcal S\}
=\sum_{h\in\P^{p,p'}_{a,b,e,f}(L)\{\mathcal S\}} q^{\mwt(h)}.
\end{displaymath}
Of course, $\P^{p,p'}_{a,b,e,f}(L,m)\{\emptyset\}=
\P^{p,p'}_{a,b,e,f}(L,m)$. 

Given $\mathcal S$ as above, we now define
$\mathcal S'=\{s+\lfloor(s+1)p/p'\rfloor:s\in\mathcal S\}$.

\begin{corollary}\label{IntBijCor}
For $1\le p<p'$, let $1\le a,b<p'$ and $e,f\in\{0,1\}$,
with $\delta^{p,p'}_{a,e}=0$.
Let $\mathcal S$ be interfacial in the $(p,p')$-model with
$a\ne s\ne b$ for all $s\in\mathcal S$.
Set $a'=a+e+\lfloor ap/p'\rfloor$ and $b'=b+f+\lfloor bp/p'\rfloor$.
Fix $m_0,m_1\ge0$.
Then
\begin{displaymath}
\begin{array}{l}
\displaystyle
\mchi^{p,p'+p}_{a',b',e,f}(m_0,m_1)\{\mathcal S'\}\\[2mm]
\hskip20mm
\displaystyle
=q^{\frac{1}{4}\left( (m_0-m_1)^2 - \beta^2 \right)}
\sum_{\scriptstyle m\equiv m_0\atop
  \scriptstyle\strut(\mbox{\scriptsize\rm mod}\,2)}
\left[{\frac{1}{2}(m_0+m)\atop m_1}\right]_q
\mchi^{p,p'}_{a,b,e,f}(m_1,m)\{\mathcal S\},
\end{array}
\end{displaymath}
where $\beta=\beta^{p,p'}_{a,b,e,f}$.
\end{corollary}

\Proof Combining Lemmas \ref{BijLem} and \ref{Inter1Lem}
implies that the map 
$(h,k,\lambda)\mapsto h'$ effected by the action of
a $\B(k,\lambda)$-transform on $h$, is a bijection between
$\bigcup_{k}
\P^{p,p'}_{a,b,e,f}(m_1,2k+2m_1-m_0)\{\mathcal S\}\times{\mathcal Y}(k,m_1)$
and $\P^{p,p'+p}_{a',b',e,f}(m_0,m_1)\{\mathcal S'\}$.
The result then follows as in the proof of Corollary \ref{BijCor}.
\cqfd
\medskip

\begin{corollary}\label{IntDijCor}
For $1\le p<p'<2p$, let $1\le a,b<p'$ and $e,f\in\{0,1\}$,
with $\delta^{p'-p,p'}_{a,e}=0$.
Let $\mathcal S$ be interfacial in the $(p,p')$-model with
$a\ne s\ne b$ for all $s\in\mathcal S$.
Set $a'=a+1-e+\lfloor ap/p'\rfloor$ and  $b'=b+1-f+\lfloor bp/p'\rfloor$.
Fix $m_0,m_1\ge0$.
Then
\begin{displaymath}
\begin{array}{l}
\displaystyle
\mchi^{p,p'+p}_{a',b',1-e,1-f}(m_0,m_1;q)\{\mathcal S'\}\\[4mm]
\hskip10mm
\displaystyle
=q^{\frac{1}{4}\left( m_1^2 + (m_0-m_1)^2 - \alpha^2 - \beta^2 \right)}\\
\displaystyle
\hskip24mm\times
\sum_{\scriptstyle m\equiv m_0-m_1\atop
  \scriptstyle\strut(\mbox{\scriptsize\rm mod}\,2)}
\hskip-1mm
\left[{\frac{1}{2}(m_0+m_1-m)\atop m_1}\right]_q
\mchi^{p'-p,p'}_{a,b,e,f}(m_1,m;q^{-1})\{\mathcal S\},
\end{array}
\end{displaymath}
where $\alpha=\alpha^{p,p'}_{a,b}$ and $\beta=\beta^{p,p'}_{a,b,1-e,1-f}$.
\end{corollary}

\Proof Combining Lemmas \ref{DijLem} and \ref{Inter2Lem}
implies that the map 
$(h,k,\lambda)\mapsto h'$ effected by the action of
a $\D$-transform on $h$ immediately followed by a $\B(k,\lambda)$-transform,
is a bijection between
$\bigcup_{k}
\P^{p'-p,p'}_{a,b,e,f}(m_1,m_0-m_1-2k)\{\mathcal S\}\times{\mathcal Y}(k,m_1)$
and $\P^{p,p'+p}_{a',b',1-e,1-f}(m_0,m_1)\{\mathcal S'\}$.
The result then follows as in the proof of Corollary \ref{DijCor}.
\cqfd
\medskip


\section{Extending and truncating paths}\label{ExtAttSec}

\subsection{Extending paths}\label{ExtendSec}
In this section, we specify a process by which a path
$h\in\P^{p,p'}_{a,b,e,f}(L)$ may be extended by
a single unit to its left, or by a single unit to its right.
One extension may follow the other to yield a path of length $L+2$.

Path extension on the left is restricted to where
$\delta^{p,p'}_{a,e}=0$ so that the pre-segment of $h$ lies in the even band.

We obtain $h'$ by defining
$h'_0=a'=a+(-1)^{e}$ and $h'_{i}=h_{i-1}$ for $1\le i\le L+1$.
In particular, $\pi(h')=0$.
We also define $e(h')=e'=1-e$, so that then
$h'\in\P^{p,p'}_{a',b,e',f}(L+1)$.

\begin{figure}[ht]
\centerline{\epsfig{file=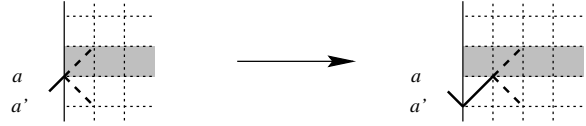}}
\caption{Extending on the left.}
\label{TypicalExt1Fig}
\medskip
\end{figure}

This extending process is depicted in Fig.~\ref{TypicalExt1Fig}.

\begin{lemma}\label{Extend1Lem}
Let $h\in\P^{p,p'}_{a,b,e,f}(L)$, where $\delta^{p,p'}_{a,e}=0$.
Let $h'\in\P^{p,p'}_{a',b,e',f}(L')$ be
obtained from $h$ by the above process of path extension.
If $\Delta=a'-a$ then $\Delta=(-1)^e=-(-1)^{e'}$,
and
\begin{displaymath}
\begin{array}{l}
\bullet\quad L'=L+1;\\[1mm]
\bullet\quad m(h')=m(h);\\[1mm]
\bullet\quad \mwt(h')=\mwt(h)+ \frac{1}{2}(L-m(h)+\Delta\beta(h)).
\end{array}
\end{displaymath}
Furthermore, $\alpha^{p,p'}_{a',b}=\alpha^{p,p'}_{a,b}-\Delta$
and $\beta^{p,p'}_{a',b,e',f}=\beta^{p,p'}_{a,b,e,f}-\Delta$.
\end{lemma}

\Proof That $\Delta=(-1)^e=-(-1)^{e'}$ is immediate from the definition.
Let $h$ have striking sequence
$\left({a_1\atop b_1}\:{a_2\atop b_2}\:{a_3\atop b_3}\:
 {\cdots\atop\cdots}\:{a_l\atop b_l} \right)^{(e,f,d)}$.

If $e=d$, we are restricted to the case $\pi(h)=0$,
since $\delta^{p,p'}_{a,e}=0$.
The striking sequence of $h'$ is then
$\left({0\atop1}\:{a_1\atop b_1}\:{a_2\atop b_2}\:
 {\cdots\atop\cdots}\:{a_l\atop b_l} \right)^{(e',f,e')}$.
Thereupon, since $\pi(h')=0$, we obtain $m(h')=m(h)$.
In this case we immediately obtain,
via Lemma \ref{WtHashLem}, that $\mwt(h')=\mwt(h)+(b_1+b_3+\cdots)$.
Thereupon, since $\Delta=(-1)^e=(-1)^d$,
$\beta(h)=(-1)^d((b_1+b_3+\cdots)-(b_2+b_4+\cdots))$
and $m(h)=(a_1+a_2+a_3+\cdots)$, we obtain
$\mwt(h')=\mwt(h)+(L-m(h)+\Delta\beta(h))/2$.

If $e\ne d$, the striking sequence of $h'$ is
$\left({a_1+1-\pi\atop b_1+\pi}\:{a_2\atop b_2}\:
 {\cdots\atop\cdots}\:{a_l\atop b_l} \right)^{(e',f,e')}$.
Then $m(h')=1-\pi+\sum_{i=1}^la_i$
which equals $m(h)=(e+d+\pi)\,\mod2+\sum_{i=1}^la_i$
for both $\pi=0$ and $\pi=1$.
Here Lemma \ref{WtHashLem} implies that $\mwt(h')=\mwt(h)+(b_2+b_4+\cdots)$.
Thereupon, since $\Delta=(-1)^e=-(-1)^d$,
$\beta(h)=(-1)^d((b_1+b_3+\cdots)-(1-\pi+b_2+b_4+\cdots))$
and $m(h)=(1-\pi+a_1+a_2+a_3+\cdots)$, we also obtain
$\mwt(h')=\mwt(h)+(L-m(h)+\Delta\beta(h))/2$.

That $\alpha^{p,p'}_{a',b}=\alpha^{p,p'}_{a,b}-\Delta$ is immediate.
Since $\pi(h')=0$ then $\lfloor a'p/p'\rfloor=\lfloor ap/p'\rfloor$.
That $\beta^{p,p'}_{a',b,e',f}=\beta^{p,p'}_{a,b,e,f}-\Delta$ now follows.
\cqfd
\medskip

In the following lemma, we consider the special case when
$2p<p'<3p$ so that the first and second bands of the $(p,p')$-model
are even and odd respectively.
We then only consider path extension into the first
or the $(p'-2)$th band of the $(p,p')$-model.

\begin{lemma}\label{ExtGen1Lem}
Let $2<2p<p'<3p$ and either $a=2$ and $e=1$, or $a=p'-2$ and $e=0$.
Then $a$ is interfacial in the $(p,p')$-model.
Let $\mathcal S$ be interfacial in the $(p,p')$-model, and
set $\Delta=(-1)^e$, $a'=a+\Delta$ and $e'=1-e$.
Then:
\begin{displaymath}
\mchi^{p,p'}_{a',b,e',f}(L,m)\{\mathcal S\cup\{a\}\}
=
q^{\frac{1}{2}(L-1-m+\Delta\beta)}
\mchi^{p,p'}_{a,b,e,f}(L-1,m)\{\mathcal S\},
\end{displaymath}
where $\beta=\beta^{p,p'}_{a,b,e,f}$.

In addition, $\alpha^{p,p'}_{a',b}=\alpha^{p,p'}_{a,b}-\Delta$,
$\beta^{p,p'}_{a',b,e',f}=\beta^{p,p'}_{a,b,e,f}-\Delta$.
\end{lemma}

\Proof Since $2p<p'<3p$, it follows that
$0=\lfloor 2p/p'\rfloor\ne\lfloor 3p/p'\rfloor$ wherepon
$2$ and $p'-2$ are both interfacial in the $(p,p')$-model,
and $\delta^{p,p'}_{a,e}=0$.

Let $h\in\P^{p,p'}_{a,b,e,f}(L-1,m)\{\mathcal S\}$.
Extend $h$ on the left to obtain $h'$ with $h'_0=a'=a+\Delta$.
Clearly, $h'$ attains $a$.
Then, Lemma \ref{Extend1Lem} implies that
$h'\in\P^{p,p'}_{a',b,e',f}(L,m)\{\mathcal S\cup\{a\}\}$.

Conversely, any such $h'$ arises from
some $h\in\P^{p,p'}_{a,b,e,f}(L-1,m)\{\mathcal S\}$ in this way
since either $h'_0=1$ and $e'=0$, or $h'_0=p'-1$ and $e'=1$.
The result then follows from the expression for $\mwt(h')$
given in Lemma \ref{Extend1Lem}, and $\beta(h)=\beta^{p,p'}_{a,b,e,f}$
from Lemma \ref{BetaConstLem}.

The final statement also follows from Lemma \ref{Extend1Lem}.
\cqfd
\medskip

For $h\in\P^{p,p'}_{a,b,e,f}(L)$, we now define path
extension to the right in a similar way.
Here we restrict path extension to the cases where $\delta^{p,p'}_{b,f}=0$
so that the post-segment of $h$ lies in the even band.

We obtain $h'$ by defining
$h'_{i}=h_{i}$ for $0\le i\le L$ and $h'_{L+1}=b'=b+(-1)^{f}$ and 
We also define $f(h')=f'=1-f$, so that then
$h'\in\P^{p,p'}_{a,b',e,f'}(L+1)$.

\begin{figure}[ht]
\centerline{\epsfig{file=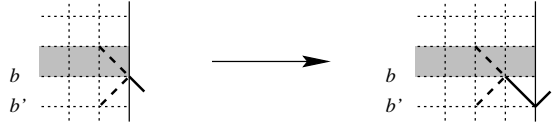}}
\caption{Extending on the right.}
\label{TypicalExt2Fig}
\medskip
\end{figure}

This extending process is depicted in Fig.~\ref{TypicalExt2Fig}.

\begin{lemma}\label{Extend2Lem}
Let $h\in\P^{p,p'}_{a,b,e,f}(L)$, where $\delta^{p,p'}_{b,f}=0$.
Let $h'\in\P^{p,p'}_{a,b',e,f'}(L')$ be
obtained from $h$ by the above process of path extension.
If $\Delta=b'-b$ then $\Delta=(-1)^f=-(-1)^{f'}$,
and
\begin{displaymath}
\begin{array}{l}
\bullet\quad L'=L+1;\\[1mm]
\bullet\quad m(h')=m(h);\\[1mm]
\bullet\quad \mwt(h')=\mwt(h)+ \frac{1}{2}(L-\Delta\alpha(h)).
\end{array}
\end{displaymath}
Furthermore, $\alpha^{p,p'}_{a,b'}=\alpha^{p,p'}_{a,b}+\Delta$
and $\beta^{p,p'}_{a,b',e,f'}=\beta^{p,p'}_{a,b,e,f}+\Delta$.
\end{lemma}

\Proof That $\Delta=(-1)^f=-(-1)^{f'}$ is immediate from the definition.
Let $h$ have striking sequence
$\left({a_1\atop b_1}\:{a_2\atop b_2}\:{a_3\atop b_3}\:
 {\cdots\atop\cdots}\:{a_l\atop b_l} \right)^{(e,f,d)}$.
It is easily checked that the $L$th vertex of $h'$ is scoring
if and only if the $L$th vertex of $h$ is scoring.

Then, if the extending segment is in the same direction as the $L$th
segment, $h'$ has striking sequence
$\left({a_1\atop b_1}\:{a_2\atop b_2}\:{a_3\atop b_3}\:
 {\cdots\atop\cdots}\:{a_l\atop b_l+1} \right)^{(e,f',d)}$
and $\Delta=-(-1)^{d+l}$.
That $m(h')=m(h)$ is immediate.

When the extending segment is in the direction opposite to that of
the $L$th segment, $h'$ has striking sequence
$\left({a_1\atop b_1}\:{a_2\atop b_2}\:
 {\cdots\atop\cdots}\:{a_l\atop b_l}\:{0\atop 1} \right)^{(e,f',d)}$
and $\Delta=(-1)^{d+l}$.
We immediately obtain $m(h')=m(h)$ in this case.

For $1\le i\le l$, let $w_i=a_i+b_i$.
We find $\alpha(h)=-(-1)^{d+l}((w_l+w_{l-2}\cdots)-(w_{l-1}+w_{l-3}+\cdots))$.
In the first case above, Lemma \ref{WtHashLem} gives
$\mwt(h')=\mwt(h)+(w_{l-1}+w_{l-3}+w_{l-5}+\cdots)$, whereupon
we obtain $\mwt(h')=\mwt(h)+\frac12(L(h)-\Delta\alpha(h))$.
In the second case above, Lemma \ref{WtHashLem} gives
$\mwt(h')=\mwt(h)+(w_{l}+w_{l-2}+w_{l-4}+\cdots)$, and we again
obtain $\mwt(h')=\mwt(h)+\frac12(L(h)-\Delta\alpha(h))$.

That $\alpha^{p,p'}_{a,b'}=\alpha^{p,p'}_{a,b}+\Delta$ is immediate.
That $\beta^{p,p'}_{a,b',e,f'}=\beta^{p,p'}_{a,b,e,f}+\Delta$ now follows
because $\lfloor bp/p'\rfloor=\lfloor b'p/p'\rfloor$.
\cqfd
\medskip

\begin{lemma}\label{ExtGen2Lem}
Let $2<2p<p'<3p$ and either $b=2$ and $f=1$, or $b=p'-2$ and $f=0$.
Then $b$ is interfacial in the $(p,p')$-model.
Let $\mathcal S$ be interfacial in the $(p,p')$-model, and
set $\Delta=(-1)^f$, $b'=b+\Delta$ and $f'=1-f$.
Then:
\begin{displaymath}
\mchi^{p,p'}_{a,b',e,f'}(L,m)\{\mathcal S\cup\{b\}\}
=
q^{\frac{1}{2}(L-1-\Delta\alpha)}
\mchi^{p,p'}_{a,b,e,f}(L-1,m)\{\mathcal S\},
\end{displaymath}
where $\alpha=\alpha^{p,p'}_{a,b}$.

In addition, $\alpha^{p,p'}_{a,b'}=\alpha^{p,p'}_{a,b}+\Delta$ and
$\beta^{p,p'}_{a,b',e,f'}=\beta^{p,p'}_{a,b,e,f}+\Delta$.
\end{lemma}

\Proof
\Proof Since $2p<p'<3p$, it follows that
$0=\lfloor 2p/p'\rfloor\ne\lfloor 3p/p'\rfloor$ wherepon
$2$ and $p'-2$ are both interfacial in the $(p,p')$-model,
and $\delta^{p,p'}_{b,f}=0$.

Let $h\in\P^{p,p'}_{a,b,e,f}(L-1,m)\{\mathcal S\}$.
Extend this path on the right to obtain $h'$ with $h'_L=b'=b+\Delta$.
Clearly, $h'$ attains height $b$.
Then, via Lemma \ref{Extend2Lem},
$h'\in\P^{p,p'}_{a,b',e,f'}(L,m)\{\mathcal S\cup\{b\}\}$.
Conversely, any such $h'$ arises in this way from some
$h\in\P^{p,p'}_{a,b,e,f}(L-1,m)\{\mathcal S\}$,
since either $h'_L=1$ and $f'=1$, or $h'_L=p'-1$ and $f'=0$.
The required result then follows from the expression for $\mwt(h')$
given in Lemma \ref{Extend2Lem}, and $\alpha(h)=\alpha^{p,p'}_{a,b}$
from Lemma \ref{BetaConstLem}.

The final statement follows from Lemma \ref{Extend2Lem}.
\cqfd
\medskip

\subsection{Truncating paths}\label{AttenSec}
In this section, we specify a process by which a path
$h\in\P^{p,p'}_{a,b,e,f}(L)$, for $L>0$ may be shortened by removing
just the leftmost (first) segment, or by removing just the rightmost
($L$th) segment.
Consequently, the new path $h'$ is of length $L'=L-1$.
One shortening may follow the other to yield a path of length $L-2$.

In fact, we will only use these shortening processes when $p'>2p$,
so that in particular, the 1st and the $(p'-2)$th bands of the
$(p,p')$-model are even.

Shortening on the left side will occur only when $a=1$ or $a=p'-1$
so that the removed segment is in an even band, and will occur
when the 0th vertex is scoring.

\begin{lemma}\label{AttenGen1Lem}
Let $p'>2p$ and either $a=1$ and $e=0$, or $a=p'-1$ and $e=1$.
Let $\mathcal S$ be interfacial in the $(p,p')$-model, with
$a\notin\mathcal S$.
Define $\Delta=-(-1)^e$, $e'=1-e$ and $a'=a-\Delta$.
Then
\begin{displaymath}
\begin{array}{l}
\displaystyle
\mchi^{p,p'}_{a',b,e',f}(L,m)\{\mathcal S\}
=
q^{-\frac{1}{2}(L+1-m+\Delta\beta)}
\mchi^{p,p'}_{a,b,e,f}(L+1,m)\{\mathcal S\},
\end{array}
\end{displaymath}
where $\beta=\beta^{p,p'}_{a,b,e,f}$.

In addition, $\alpha^{p,p'}_{a',b}=\alpha^{p,p'}_{a,b}+\Delta$,
and $\beta^{p,p'}_{a',b,e',f}=\beta^{p,p'}_{a,b,e,f}+\Delta$.
\end{lemma}

\Proof Let $h\in\P^{p,p'}_{a,b,e,f}(L+1,m)\{\mathcal S\}$, and note
that necessarily $h_1=a'$.
Let $h'\in\P^{p,p'}_{a',b,e',f}(L,m)\{\mathcal S\}$ be defined
by $h'_i=h_{i+1}$ for $0\le i\le L$.
The lemma then follows on noting that $\delta^{p,p'}_{a',e'}=0$
and using Lemma \ref{Extend1Lem} after switching the roles of
$h$ and $h'$ there.
\cqfd
\medskip

Shortening on the right side will occur only when $b=1$ or $b=p'-1$
so that the removed segment is in an even band, and will occur
when the $L$th vertex is scoring.

\begin{lemma}\label{AttenGen2Lem}
Let $p'>2p$ and either $b=1$ and $f=0$, or $b=p'-1$ and $f=1$.
Let $\mathcal S$ be interfacial in the $(p,p')$-model, with
$b\notin\mathcal S$.
Define $\Delta=-(-1)^f$, $f'=1-f$ and $b'=b-\Delta$.
Then
\begin{displaymath}
\mchi^{p,p'}_{a,b',e,f'}(L,m)\{\mathcal S\}
=
q^{-\frac{1}{2}(L+1-\Delta\alpha)}
\mchi^{p,p'}_{a,b,e,f}(L+1,m)\{\mathcal S\},
\end{displaymath}
where $\alpha=\alpha^{p,p'}_{a,b}$.

In addition, $\alpha^{p,p'}_{a,b'}=\alpha^{p,p'}_{a,b}-\Delta$,
and $\beta^{p,p'}_{a,b',e,f'}=\beta^{p,p'}_{a,b,e,f}-\Delta$.
\end{lemma}

\Proof Let $h\in\P^{p,p'}_{a,b,e,f}(L+1,m)\{\mathcal S\}$, and note
that necessarily $h_L=b'$.
Let $h'\in\P^{p,p'}_{a,b',e,f'}(L,m)\{\mathcal S\}$ be defined
by $h'_i=h_{i}$ for $0\le i\le L$.
The lemma then follows on noting that $\delta^{p,p'}_{b',f'}=0$
and using Lemma \ref{Extend2Lem} after switching the roles of
$h$ and $h'$ there.
\cqfd
\medskip


\goodbreak
\section{Fermionic expressions}\label{FermSec}
\nobreak
\subsection{Results}\label{ResultsSec}
\nobreak

In this section, we fix co-prime $p$ and $p'$, and fix
$a,b\in\mathcal T\cup\mathcal T'$, with $1\le a,b<p'$.
We make use of the definitions of \ref{ContFSec} and \ref{TakSec}.
For certain $c$, we present two fermionic expressions for
$\P^{p,p'}_{a,b,c}(L)$.
The value of $c$ depends on $b$ and, for $p'>2p$, is given by:
\begin{equation}\label{CEq}
c=
\left\{ \begin{array}{ll}
    2
       &\mbox{if } b=1;\\[1.5mm]
    b-1
       &\mbox{if } 1<b\le t_1;\\[1.5mm]
    p'-2
       &\mbox{if } b=p'-1;\\[1.5mm]
    b+1
       &\mbox{if } p'-t_1\le b<p'-1;\\[1.5mm]
    b\pm1
       &\mbox{otherwise.}
         \end{array} \right.
\end{equation}
For $p'<2p$, change $t_1$ to $t_2$ in this definition.

The statement of these fermionic expressions requires the following notation.
For convenience, set $a^L=a$ and $a^R=b$.
Now, for $A\in\{L,R\}$, define $\sigma^A$ such that:
\begin{equation}\label{SigmaEq}
\kappa_{\sigma^A}=\quad
\left\{
  \begin{array}{ll}
  \displaystyle
  a^A\quad
  &\mbox{if }a^A\in\mathcal T;\\[1mm]
  \displaystyle
  p'-a^A\quad
  &\mbox{if }a^A\in\mathcal T^\prime.
  \end{array}
\right.
\end{equation}
For $0\le j\le t$, define%
\footnote{In this paper, all vectors $\boldQ$, $\boldm$, $\hat{\boldm}$,
$\boldn$, $\boldu$, $\boldDelta$ and $\bolde$
should be considered as column vectors. However, for typographical
convenience, we shall express their components in row vector form.}
$\bolde_j=(e_1,e_2,\ldots,e_t)$ with $e_i=\delta_{ij}$.
Then define
\begin{equation}
{\boldu}^A
=\bolde_{\sigma^A}-\sum_{k:\sigma^A\le t_k<t} \bolde_{t_k}
+
\left\{
  \begin{array}{ll}
  \displaystyle
  0\quad
  &\mbox{if }a^A\in\mathcal T;\\[1mm]
  \displaystyle
  \bolde_t\quad
  &\mbox{if }a^A\in\mathcal T^\prime,
  \end{array}
\right.
\end{equation}
and
\begin{equation}
{\boldDelta}^A
=
\left\{
  \begin{array}{ll}
  \displaystyle
  -\bolde_{\sigma^A}+\sum_{k:\sigma^A\le t_k<t} \bolde_{t_k}
  \quad
  &\mbox{if }a^A\in\mathcal T;\\[5mm]
  \displaystyle
  -\bolde_t
  +\bolde_{\sigma^A}-\sum_{k:\sigma^A\le t_k<t} \bolde_{t_k}
  \quad
  &\mbox{if }a^A\in\mathcal T^\prime.
  \end{array}
\right.
\end{equation}

We define the matrix $\boldC$ to be the $t\times t$ tri-diagonal
matrix with entries $\boldC_{ij}$ for $0\le i,j\le t-1$ where,
when the indices are in this range,

\begin{equation}\label{CDefEq}
\begin{array}{cccl}
\boldC_{j,j-1}=-1,
&\boldC_{j,j}=1,
&\boldC_{j,j+1}=\phantom{-}1,
&\hbox{if $j=t_k,\quad k=1,2,\ldots,n$;}\\
\boldC_{j,j-1}=-1,
&\boldC_{j,j}=2,
&\boldC_{j,j+1}=-1,
&\hbox{$0\le j<t$ otherwise.}
\end{array}
\end{equation}

It is also useful to define $\hboldC$ to be the $t\times t$
upper-triangular matrix with entries $\hboldC_{ij}={\boldC}_{ij}$,
as above, with $1\le i\le t$ and $0\le j\le t-1$.

For example, in the case $p=9$ and $p'=31$, where the continued fraction
of $p'/p$ is $(3,2,4)$ and $t_1=2$, $t_2=4$ and $t_3=8$, we have:

\begin{displaymath}
\boldC=
\left(
\vcenter{\halign{&$\scriptstyle\hfil\; #$\cr
2  &-1 &.\,&.\,&.\,&.\,&.\,\cr
-1 & 2 &-1 &.\,&.\,&.\,&.\,\cr
.\,&-1 & 1 & 1 &.\,&.\,&.\,\cr
.\,&.\,&-1 & 2 &-1 &.\,&.\,\cr
.\,&.\,&.\,&-1 & 1 & 1 &.\,\cr
.\,&.\,&.\,&.\,&-1 & 2 &-1 \cr
.\,&.\,&.\,&.\,&.\,&-1 & 2 \cr}}
\right),
\qquad
\hboldC=
\left(
\vcenter{\halign{&$\scriptstyle\hfil\; #$\cr
-1 &2  &-1 &.\,&.\,&.\,&.\,\cr
.\,&-1 & 1 & 1 &.\,&.\,&.\,\cr
.\,&.\,&-1 & 2 &-1 &.\,&.\,\cr
.\,&.\,&.\,&-1 & 1 & 1 &.\,\cr
.\,&.\,&.\,&.\,&-1 & 2 &-1 \cr
.\,&.\,&.\,&.\,&.\,&-1 & 2 \cr
.\,&.\,&.\,&.\,&.\,&.\,&-1 \cr}}
\right).
\end{displaymath}

Since $\hboldC$ is upper-triangular, its inverse is readily obtained.
Given a $t$-dimensional vector $\boldu$, we then define
$Q_i\in\{0,1\}$ for $0\le i<t$, by%
\footnote{For $\boldv=(v_1,v_2,\ldots,v_t)$,
we define $\boldv\,\mod2=(v_1\,\mod2,v_2\,\mod2,\ldots,v_t\,\mod2,)$.}
\begin{equation}\label{ParityDef}
(Q_0,Q_1,Q_2,\ldots,Q_{t-1})^T= \hboldC{}^{-1}\boldu \;\mod2.
\end{equation}
We thus define
the {\it parity vector}
$\boldQ(\boldu)=(Q_1,Q_2,\ldots,Q_{t-1})$.

Now, given a $t$-dimensional vector 
$\boldu=(u_1,u_2,\ldots,u_t)$, define the $(t-1)$-dimensional vector
$\boldu^{\flat}=(u^{\flat}_1,u^{\flat}_2,\ldots,u^{\flat}_{t-1})$ by:
\begin{equation}
u^{\flat}_j=
\left\{
  \begin{array}{cl}
       0 &\quad\mbox{if } t_k<j\le t_{k+1},\; k\equiv 0\,(\mod 2);\\[0.5mm]
     u_j &\quad\mbox{if } t_k<j\le t_{k+1},\; k\not\equiv 0\,(\mod 2),
  \end{array} \right.
\end{equation}
and the $(t-1)$-dimensional vector
$\boldu^{\sharp}=(u^{\sharp}_1,u^{\sharp}_2,\ldots,u^{\sharp}_{t-1})$ by:
\begin{equation}
u^{\sharp}_j=
\left\{
  \begin{array}{cl}
     u_j &\quad\mbox{if } t_k<j\le t_{k+1},\; k\equiv 0\,(\mod 2);\\[0.5mm]
       0 &\quad\mbox{if } t_k<j\le t_{k+1},\; k\not\equiv 0\,(\mod 2).
  \end{array} \right.
\end{equation}
Then, of course, $(\boldu)_j=(\boldu^{\flat}+\boldu^{\sharp})_j$
for $1\le j<t$.
For convenience, we sometimes write $\boldu_\flat$ and
$\boldu_\sharp$ for $\boldu^\flat$ and $\boldu^\sharp$ respectively.

Finally, we define a value $\gamma$ that depends on
$\boldDelta^L$ and $\boldDelta^R$.
This value is obtained by iteratively generating the sequences
$(\beta_{t},\beta_{t-1},\ldots,\beta_0)$,
$(\alpha_{t},\alpha_{t-1},\ldots,\alpha_0)$,
and $(\gamma_{t},\gamma_{t-1},\ldots,\gamma_0)$ as follows.
Let $\alpha_t=\beta_t=\gamma_t=0$.
Now, for $j=t,t-1,\ldots,1$, obtain $\alpha_{j-1}$, $\beta_{j-1}$,
and $\gamma_{j-1}$
from $\alpha_j$, $\beta_j$, and $\gamma_j$ in the following
three stages.
Firstly, obtain:
\begin{equation}\label{Const1Eq}
(\beta'_{j-1},\gamma'_{j-1})=
(\beta_j+(\boldDelta^L)_j-(\boldDelta^R)_j,\gamma_j+2\alpha_j(\boldDelta^R)_j).
\end{equation}
Then obtain:
\begin{equation}\label{Const2Eq}
(\alpha''_{j-1},\gamma''_{j-1})=
(\alpha_j+\beta'_{j-1},\gamma'_{j-1}-(\beta'_{j-1})^2).
\end{equation}
Finally, set
\begin{equation}\label{Const3Eq}
\begin{array}{l}
(\alpha_{j-1},\beta_{j-1},\gamma_{j-1})\\[2mm]
\quad=
\left\{
  \begin{array}{cl}
     (\alpha''_{j-1},\alpha''_{j-1}-\beta'_{j-1},
              -(\alpha''_{j-1})^2-\gamma''_{j-1}) &
          \mbox{if } j=t_k+1,\;  1\le k\le n;\\
     (\alpha''_{j-1},\beta'_{j-1},\gamma''_{j-1}) &
          \mbox{otherwise}.
  \end{array} \right.
\end{array}
\end{equation}
We then set $\gamma=\gamma_0$.


\begin{theorem}\label{Ferm1Thrm}
If $a,b\in\mathcal T\cup\mathcal T'$, define everything as above.
Then:
\begin{displaymath}
\begin{array}{l}
\displaystyle
\hskip-2mm
\ochi^{p,p'}_{a,b,c}(L)=\\[1.5mm]
\displaystyle
\hskip5mm
  \sum_{\sboldm\equiv\sboldQ(\sboldu^L+\sboldu^R)}
  \hskip-5mm
  q^{\frac{1}{4}\hat{\sboldm}^T\sboldC\hat{\sboldm}-\frac{1}{4} L^2
   -\frac{1}{2}(\sboldu^L_\flat+\sboldu^R_\sharp)\cdot\sboldm
  +\frac{1}{4}\gamma}
  \prod_{j=1}^{t-1}
  \left[
  {m_j-\frac{1}{2}(\hboldC\hat{\boldm}\!-\!\boldu^L\!-\!\boldu^R)_j\atop m_j}
  \right]_q\hskip-5mm\phantom{.}\\[8mm]
\displaystyle
\hskip30mm
+\left\{ \begin{array}{ll}
    \chi^{z_n,y_n}_{a,b,c}(L)
       &\mbox{if } a<y_n \mbox{ and } b<y_n;\\[1.5mm]
    \chi^{z_n,y_n}_{p'-a,p'-b,p'-c}(L)
       &\mbox{if } a>p'-y_n \mbox{ and } b>p'-y_n;\\[1.5mm]
   0, &\mbox{otherwise}.
         \end{array} \right.
\end{array}
\end{displaymath}
With $\boldQ(\boldu^L+\boldu^R)=(Q_1,Q_2,\ldots,Q_{t-1})$,
the summation here is over all vectors $\boldm=(m_1,m_2,\ldots,m_{t-1})$
such that $m_j\in\Z_{\ge0}$ and $m_j\equiv Q_j\,(\mod2)$ for $1\le j<t$.
Then $\hat{\boldm}=(L,m_1,m_2,\ldots,m_{t-1})$.
\end{theorem}

The second fermionic expression for $\ochi^{p,p'}_{a,b,c}(L)$
that we present, involves the modified form
$\left[ {A \atop B} \right]^\prime_q$
of the Gaussian polynomial defined in (\ref{ModGaussian}).


\begin{theorem}\label{Ferm2Thrm}
If $a,b\in\mathcal T\cup\mathcal T'$, define everything as above.
Then, if $L\ge0$:
\begin{displaymath}
\begin{array}{l}
\displaystyle
\hskip-2mm
\ochi^{p,p'}_{a,b,c}(L)=\\[1.5mm]
\displaystyle
\hskip3mm
  \sum_{\sboldm\equiv\sboldQ(\sboldu^L+\sboldu^R)}
  \hskip-5mm
  q^{\frac{1}{4}\hat{\sboldm}^T\sboldC\hat{\sboldm}-\frac{1}{4} L^2
   -\frac{1}{2}(\sboldu^L_\flat+\sboldu^R_\sharp)\cdot\sboldm
  +\frac{1}{4}\gamma}
  \prod_{j=1}^{t-1}
  \left[
  {m_j-\frac{1}{2}(\hboldC\hat{\boldm}\!-\!\boldu^L\!-\!\boldu^R)_j\atop m_j}
  \right]^\prime_q\!.
\end{array}
\end{displaymath}
With $\boldQ(\boldu^L+\boldu^R)=(Q_1,Q_2,\ldots,Q_{t-1})$,
the summation here is over all vectors $\boldm=(m_1,m_2,\ldots,m_{t-1})$
such that $m_j\in\Z_{\ge0}$ and $m_j\equiv Q_j\,(\mod2)$ for $1\le j<t$.
Then $\hat{\boldm}=(L,m_1,m_2,\ldots,m_{t-1})$.
\end{theorem}

\subsection{Carrying out the induction}\label{CarrySec}

With $p$ and $p'$ fixed, employ the definitions of Section \ref{ContFSec}.
Then, for $0\le i\le t$,
let $k(i)$ be such that $t_{k(i)}\le i<t_{k(i)+1}$ (i.e.\ $k(i)=\zeta(i+1)$),
and define $p_i$ and $p_i^{\prime}$ to be the positive co-prime
integers for which $p_i^{\prime}/p_i$ has continued
fraction $(t_{k(i)+1}+1-i,c_{k(i)+1},\ldots,c_n)$.
Thus $p_i^{\prime}/p_i$ has rank $t-i$.
As in Section \ref{TakSec}, we obtain
Takahashi lengths $\{\kappa^{(i)}_j\}_{j=0}^{t-i}$ and
truncated Takahashi lengths $\{\tilde\kappa^{(i)}_j\}_{j=0}^{t-i}$
for $p_i^{\prime}/p_i$.

\begin{lemma}\label{IndTakLem} Let $1\le i\le t$. If $i\ne t_{k(i)}$
then:
\begin{displaymath}
\begin{array}{lll}
p^{(i-1)\prime}=
 &p^{(i)\prime}+p^{(i)};\\[0.5mm]
p^{(i-1)}=
 &p^{(i)};\\[0.5mm]
\kappa^{(i-1)}_j=
 &\kappa^{(i)}_{j-1}+\tkappa^{(i)}_{j-1}\qquad
 &(1\le j\le t^{(i-1)});\\[0.5mm]
\tkappa^{(i-1)}_j=
 &\tkappa^{(i)}_{j-1}\qquad
 &(1\le j\le t^{(i-1)}).
\end{array}
\end{displaymath}
If $i=t_{k(i)}$ then:
\begin{displaymath}
\begin{array}{lll}
p^{(i-1)\prime}=
 &2p^{(i)\prime}-p^{(i)};\\[0.5mm]
p^{(i-1)}=
 &p^{(i)\prime}-p^{(i)};\\[0.5mm]
\kappa^{(i-1)}_j=
 &2\kappa^{(i)}_{j-1}-\tkappa^{(i)}_{j-1}\qquad
 &(2\le j\le t^{(i-1)});\\[0.5mm]
\tkappa^{(i-1)}_j=
 &\tkappa^{(i)}_{j-1}-\tkappa^{(i)}_{j-1}\qquad
 &(2\le j\le t^{(i-1)}).
\end{array}
\end{displaymath}
\end{lemma}

\Proof If $i\ne t_{k(i)}$ then $k(i-1)=k(i)$.
Then $p^{(i)\prime}/p^{(i)}$ and
$p^{(i-1)\prime}/p^{(i-1)}$ have continued fractions
$(t_{k(i)}+1-i,c_{k(i)+1},\ldots,c_n)$ and
$(t_{k(i)}+2-i,c_{k(i)+1},\ldots,c_n)$ respectively.
That $p^{(i-1)\prime}=p^{(i)\prime}+p^{(i)}$ and $p^{(i-1)}=p^{(i)}$
follows immediately.
The expressions for $\kappa^{(i-1)}_j$ and $\tkappa^{(i-1)}_j$
then follow from Lemma \ref{BmodelLem}.

If $i= t_{k(i)}$ then $k(i-1)=k(i)-1$.
Then $p^{(i)\prime}/p^{(i)}$ and
$p^{(i-1)\prime}/p^{(i-1)}$ have continued fractions
$(c_{k(i)},c_{k(i)+1},\ldots,c_n)$ and
$(2,c_{k(i)},c_{k(i)+1},\ldots,c_n)$ respectively.
That $p^{(i-1)\prime}=2p^{(i)\prime}-p^{(i)}$
and $p^{(i-1)}=p^{(i)\prime}-p^{(i)}$
follows immediately.
The expressions for $\kappa^{(i-1)}_j$ and $\tkappa^{(i-1)}_j$
then follow from combining Lemma \ref{DmodelLem} with Lemma \ref{BmodelLem}.
\cqfd
\medskip

As above, take $A\in\{R,L\}$.
If $a^A\in\mathcal T$, set
\begin{displaymath}
\begin{array}{ll}
a^A_i=&\left\{
  \begin{array}{ll}
    1 & \mbox{if } \sigma^A\le i<t;\\
    \kappa^{(i)}_{\sigma^A-i} & \mbox{if } 0\le i\le\sigma^A,
  \end{array}
\right.\\[5mm]
a^{A\prime}_i=&\left\{
  \begin{array}{ll}
    1+\delta_{i,t_{k(i)}} & \mbox{if } \sigma^A\le i<t;\\
    \kappa^{(i-1)}_{\sigma^A-i+1} & \mbox{if } 0\le i<\sigma^A,
  \end{array}
\right.
\end{array}
\end{displaymath}
and if $a^A\in\mathcal T'$, set
\begin{displaymath}
\begin{array}{ll}
a^A_i=&\left\{
  \begin{array}{ll}
    p_i'-1 & \mbox{if } \sigma^A\le i<t;\\
    p_i'-\kappa^{(i)}_{\sigma^A-i} & \mbox{if } 0\le i\le\sigma^A.
  \end{array}
\right.\\[5mm]
a^{A\prime}_i=&\left\{
  \begin{array}{ll}
    p_i'-1-\delta_{i,t_{k(i)}} & \mbox{if } \sigma^A\le i<t;\\
    p_i'-\kappa^{(i-1)}_{\sigma^A-i+1} & \mbox{if } 0\le i<\sigma^A.
  \end{array}
\right.
\end{array}
\end{displaymath}

In addition, define $k^A$ to be such that
$t_{k^A}<\sigma^A\le t_{k^A+1}$.
Then, if $a^A\in\mathcal T$, set
\begin{displaymath}
e^A_i=\left\{
  \begin{array}{ll}
    0 & \mbox{if } \sigma^A\le i<t;\\
    \delta^{(2)}_{k,k^A} & \mbox{if } 0\le i<\sigma^A,
  \end{array}
\right.
\end{displaymath}
and if $a^A\in\mathcal T'$, set
\begin{displaymath}
e^A_i=\left\{
  \begin{array}{ll}
    1 & \mbox{if } \sigma^A\le i<t;\\
    1-\delta^{(2)}_{k,k^A} & \mbox{if } 0\le i<\sigma^A.
  \end{array}
\right.
\end{displaymath}

\begin{lemma}\label{A_AprimeLem}
Let $1\le i<t$.
Then for $A\in\{L,R\}$:
\begin{displaymath}
a^{A\prime}_i=
\left\{
  \begin{array}{ll}
     \displaystyle a^A_i+\left\lfloor\frac{a^A_ip_i}{p'_i}\right\rfloor+e^A_i
       & \mbox{if } i\ne t_{k(i)};  \\[3.5mm]
     \displaystyle 2a^A_i-\left\lfloor\frac{a^A_ip_i}{p'_i}\right\rfloor-e^A_i
       & \mbox{if } i= t_{k(i)}.
  \end{array}
\right.
\end{displaymath}
\end{lemma}

\Proof For $p'_i/p_i$, in view of the continued fraction specified
above, the analogues of the quantities defined in (\ref{ZoneEq})
are $t'_j=t_{k(i)+j}-i$ for $1\le j\le n-k(j)+1$.
For $i<\sigma^A$, the various cases are then readily proved
using Lemmas \ref{TakBandLem} and \ref{IndTakLem}.
For $i\ge\sigma^A$, the results follow immediately.
\cqfd
\medskip


For each $t$-dimensional vector
$\boldu=(u_1,u_2,\ldots,u_t)$, define the $(t-1)$-dimensional vector
$\boldu^{(\flat,k)}=(u^{(\flat,k)}_1,u^{(\flat,k)}_2,\ldots,
u^{(\flat,k)}_{t-1})$ by
\begin{equation}
u^{(\flat,k)}_j=
\left\{
  \begin{array}{cl}
       0 &\quad\mbox{if } t_{k'}<j\le t_{k'+1},
                                \; k'\equiv k\;(\mod 2);\\[0.5mm]
      u_j &\quad\mbox{if } t_{k'}<j\le t_{k'+1},
                                \; k'\not\equiv k\;(\mod 2),
  \end{array} \right.
\end{equation}
and the $(t-1)$-dimensional vector
$\boldu^{(\sharp,k)}=(u^{(\sharp,k)}_1,u^{(\sharp,k)}_2,\ldots,
u^{(\sharp,k)}_{t-1})$ by
\begin{equation}
u^{(\sharp,k)}_j=
\left\{
  \begin{array}{cl}
     u_j &\quad\mbox{if } t_{k'}<j\le t_{k'+1},
                                \; k'\equiv k\;(\mod 2);\\[0.5mm]
      0 &\quad\mbox{if } t_{k'}<j\le t_{k'+1},
                                \; k'\not\equiv k\;(\mod 2),
  \end{array} \right.
\end{equation}
For convenience, we sometimes write
$\boldu_{(\flat,k)}$ instead of $\boldu^{(\flat,k)}$, and
$\boldu_{(\sharp,k)}$ instead of $\boldu^{(\sharp,k)}$.

Now for $0\le i\le t-2$, define:
\begin{equation}\label{FDef1Eq}
\begin{array}{l}
\displaystyle
F^{(i)}_{a,b}(\boldu^{L},\boldu^{R},m_i,m_{i+1};q)=\\[3mm]
\quad
\displaystyle
\sum
\left(
  q^{\frac{1}{4}\hat{\sboldm}^{(i+1)T}\sboldC\hat{\sboldm}^{(i+1)}
   +\frac{1}{4} m_i^2
   -\frac{1}{2} m_im_{i+1}
   -\frac{1}{2}(\sboldu^{L}_{(\flat,k(i))}+\sboldu^{R}_{(\sharp,k(i))})\cdot
   \sboldm^{(i)}
  +\frac{1}{4}\gamma''_i}\hskip-20mm
  \phantom{\prod_{j=i+1}^{t-1}}
\right.
\\[3mm]
\displaystyle
\hskip45mm
\left.
  \prod_{j=i+1}^{t-1}
  \left[
  {m_j-\frac{1}{2}(\hboldC\hat{\boldm}^{(i)}
                       \!-\!\boldu^{L}\!-\!\boldu^{R})_j\atop m_j}
  \right]_q\right),\hskip-20mm
\end{array}
\end{equation}
where the sum here is taken over all vectors
$(m_{i+2},m_{i+3},\ldots,m_{t-1})\equiv(Q_{i+2},Q_{i+3},\ldots,Q_{t-1})$,
where $(Q_1,Q_2,\ldots,Q_{t-1})=\boldQ(\boldu^{L}+\boldu^{R})$.
The $(t-1)$-dimensional vector
$\boldm^{(i)}=(0,0,\ldots,0,m_{i+1},m_{i+2},m_{i+3},\ldots,m_{t-1})$
has its first $i$ components equal to zero.
The $t$-dimensional vector
$\hat{\boldm}^{(i)}=(0,0,\ldots,0,m_{i},m_{i+1},m_{i+2},\ldots,m_{t-1})$
has its first $i$ components equal to zero.

We also define:
\begin{equation}\label{FDef2Eq}
F^{(t-1)}_{a,b}(\boldu^{L},\boldu^{R},m_{t-1},m_{t};q)
=q^{\frac{1}{4}m_{t-1}^2+\frac{1}{4}\gamma''_i}
\delta_{m_t,0}.
\end{equation}
For convenience, we set $Q_t=0$.

Since $\left[{m+n\atop m}\right]_{q^{-1}}=q^{-mn}\left[{m+n\atop m}\right]_q$,
it follows that for $0\le i\le t-2$:
\begin{equation}\label{Inverse1Eq}
\begin{array}{l}
\displaystyle
F^{(i)}_{a,b}(\boldu^{L},\boldu^{R},m_i,m_{i+1};q^{-1})=\\[3mm]
\quad
\displaystyle
\sum
\left(
  q^{\frac{1}{4}\hat{\sboldm}^{(i+1)T}\sboldC\hat{\sboldm}^{(i+1)}
   -\frac{1}{4} m_i^2
   -\frac{1}{2}(\sboldu^{L}_{(\flat,k(i)-1)}+\sboldu^{R}_{(\sharp,k(i)-1)})
   \cdot \sboldm^{(i)}
  -\frac{1}{4}\gamma''_i}
  \phantom{\prod_{j=i+1}^{t-1}}
\right.
\\[3mm]
\displaystyle
\hskip48mm
\left.
  \prod_{j=i+1}^{t-1}
  \left[
  {m_j-\frac{1}{2}(\hboldC\hat{\boldm}^{(i)}
                       \!-\!\boldu^{L}\!-\!\boldu^{R})_j\atop m_j}
  \right]_q
\right),
\end{array}
\end{equation}
where the sum here is taken over all vectors
$(m_{i+2},m_{i+3},\ldots,m_{t-1})\equiv(Q_{i+2},Q_{i+3},\ldots,Q_{t-1})$,
as above.
Of course, we also have:
\begin{equation}\label{Inverse2Eq}
F^{(t-1)}_{a,b}(\boldu^{L},\boldu^{R},m_{t-1},m_{t};q^{-1})
=q^{-\frac{1}{4}m_{t-1}^2-\frac{1}{4}\gamma''_i}
\delta_{m_t,0}.
\end{equation}

\begin{lemma}\label{CoreIndLem}
Let $0\le i<t$, $m_i\equiv Q_i$ and $m_{i+1}\equiv Q_{i+1}$.
If
\begin{displaymath}
\mathcal S^{(i)}=
\left\{ \begin{array}{ll}
    \left\{ \kappa^{(i)}_{t_n-i} \right\}
       &\mbox{if } i<t_n,\sigma^L<t_n,\sigma^R<t_n
        \mbox{ and } a,b\in\mathcal T;\\[1.5mm]
    \left\{ p'_i-\kappa^{(i)}_{t_n-i} \right\}
       &\mbox{if } i<t_n,\sigma^L<t_n,\sigma^R<t_n
        \mbox{ and } a,b\in\mathcal T';\\[1.5mm]
    \emptyset
       &\mbox{otherwise},
         \end{array} \right.
\end{displaymath}
then:
\begin{equation}\label{IndEq}
\mchi^{p_i,p_i^\prime}_{a^L_i,a^R_i,e^L_i,e^R_i}
              (m_i,m_{i+1})\left\{\mathcal S^{(i)}\right\}
=
F^{(i)}_{a,b}(\boldu^{L},\boldu^{R},m_i,m_{i+1}).
\end{equation}
In addition,
$\alpha^{p_i,p_i^\prime}_{a^L_i,a^R_i}=\alpha_i''$ and
$\beta^{p_i,p_i^\prime}_{a^L_i,a^R_i,e^L_i,e^R_i}=\beta_i'$.
\end{lemma}

\Proof
For $i=t-1$, we have $p_i^\prime=3$, $p_i=1$,
and if $a^A\in\mathcal T$ then
$a^A_i=1$, $e^A_i=0$, $(\boldDelta^A)_t=0$;
and if $a^A\in\mathcal T'$ then
$a^A_i=2$, $e^A_i=1$, $(\boldDelta^A)_t=-1$.
Furthermore, we have $i\ge t_n$.
Via (\ref{Const1Eq}) and (\ref{Const2Eq}), we obtain
$\alpha_{t-1}''=\beta_{t-1}'=(\boldDelta^L)_t-(\boldDelta^R)_t$
and
$\gamma_{t-1}''=-((\boldDelta^L)_t-(\boldDelta^R)_t)^2$.
For $i=t-1$, the first statement of our induction proposition
is now seen to hold via  Lemma \ref{SeedLem}.
The definitions of $\alpha^{p_i,p_i^\prime}_{a^L_i,a^R_i}$ and
$\beta^{p_i,p_i^\prime}_{a^L_i,a^R_i,e^L_i,e^R_i}$ then
yield the two final statements.


Now assume the result holds for a particular $i$ with
$1\le i\le t-1$.
As above, let $k(i)$ be such that $t_{k(i)}\le i<t_{{k(i)}+1}$.
First consider the case $t_{k(i)}<i<t_{{k(i)}+1}$.
Equation (\ref{Const3Eq}) gives $\alpha_i=\alpha_i''$,
$\beta_i=\beta_i'$ and $\gamma_i=\gamma_i''$.
Let $m_{i-1}\equiv Q_{i-1}$.
On setting $M=m_{i-1}+u^L_i+u^R_i$, equation
(\ref{ParityDef}) implies that $M\equiv Q_{i+1}$.
Then, use of the induction hypothesis,
Lemmas \ref{BijCor} or Lemma \ref{IntBijCor} as appropriate,
and Lemmas \ref{IndTakLem} and \ref{A_AprimeLem} yields:
\begin{equation}\label{StepEasyEq}
\begin{array}{l}
\displaystyle
\mchi^{p_{i-1},p_{i-1}^\prime}_{a^{L\prime}_{i},a^{R\prime}_{i},
                                       e^L_{i},e^R_{i}} (M,m_{i})
\left\{\mathcal S^{(i-1)}\right\}
\\[2mm]
\hskip5mm
\displaystyle
=\sum_{m_{i+1}\equiv Q_{i+1}}
q^{\frac{1}{4}(M-m_i)^2-\frac{1}{4}\beta^{2}_i}
\left[{\frac{1}{2}(M+m_{i+1})}\atop m_i\right]_q
F^{(i)}_{a,b}(\boldu^{L},\boldu^{R},m_i,m_{i+1}).
\end{array}
\end{equation}
Here, Lemma \ref{ParamLem} also gives
$\alpha^{p_{i-1},p_{i-1}^\prime}_{a^{L\prime}_{i},a^{R\prime}_{i}}
=\alpha_i+\beta_i$, and
$\beta^{p_{i-1},p_{i-1}^\prime}_{a^{L\prime}_{i},a^{R\prime}_{i},
                        e^L_{i},e^R_{i}}=\beta_i$.

That $\left\{\mathcal S^{(i-1)}\right\}$ appears
on the leftside here is because, via Lemma \ref{TakBandLem},
if $i<t_n$ then
$\kappa^{(i)}_{t_n-i}$ is interfacial in the $(p_i,p_i')$-model,
and borders the $\tkappa^{(i)}_{t_n-i}$th odd band,
and then $\kappa^{(i)}_{t_n-i}+\tkappa^{(i)}_{t_n-i}=
\kappa^{(i-1)}_{t_n-i+1}$, by Lemma \ref{IndTakLem}, and finally noting
that $i\ne t_n$ so that if $i\ge t_n$ then $i-1\ge t_n$.
(The up-down symmetry of the $(p,p')$-model implies that
if $i<t_n$ then
$p'_i-\kappa^{(i)}_{t_n-i}$ is interfacial in the $(p_i,p_i')$-model,
and borders the $(p_i-\tkappa^{(i)}_{t_n-i})$th odd band.
Then we use $(p'_i-\kappa^{(i)}_{t_n-i})+(p_i-\tkappa^{(i)}_{t_n-i})=
p'_{i-1}-\kappa^{(i-1)}_{t_n-i+1}$, from Lemma \ref{IndTakLem}.)

Since $M=m_{i-1}+u^L_i+u^R_i$, on noting that $t_k<i<t_{k+1}$, we have:
\begin{displaymath}
M+m_{i+1}=2m_i-(\hboldC\hat{\boldm}^{(i-1)}-\boldu^{L}-\boldu^{R})_i,
\end{displaymath}
and
\begin{displaymath}
\begin{array}{l}
\displaystyle
\hat{\boldm}^{(i+1)T}\boldC\hat{\boldm}^{(i+1)}+m_i^2-2m_im_{i+1}+(M-m_i)^2\\
\qquad=
\hat{\boldm}^{(i)T}\boldC\hat{\boldm}^{(i)}+M^2-2Mm_i\\
\qquad=
\hat{\boldm}^{(i)T}\boldC\hat{\boldm}^{(i)}+m_{i-1}^2-2m_im_{i-1}\\
\hskip30mm
+2(m_{i-1}-m_i)(u^L_i+u^R_i)+(u^L_i+u^R_i)^2.
\end{array}
\end{displaymath}
(In the case $i=t-1$, we require this expression after substituting
$m_t=0$.) 
Thence,
\begin{displaymath}
\begin{array}{l}
\displaystyle
\hskip-1mm
\mchi^{p_{i-1},p_{i-1}^\prime}_{a^{L\prime}_{i},a^{R\prime}_{i},
                       e^L_{i},e^R_{i}} (m_{i-1}+u^L_i+u^R_i,m_{i})
\left\{\mathcal S^{(i-1)}\right\}
\\[4mm]
\displaystyle
=\sum
\left(
q^{\frac{1}{4}\hat{\sboldm}^{(i)T}\sboldC\hat{\sboldm}^{(i)}
  -\frac{1}{2}m_im_{i-1}
  +\frac{1}{2}(m_{i-1}-m_i)(u^L_i+u^R_i)
  -\frac{1}{2}(\sboldu^{L}_{(\flat,k(i))}+\sboldu^{R}_{(\sharp,k(i))})\cdot
    \sboldm^{(i)}}
\vrule height20pt depth2pt width0pt
\hskip-10mm
\right.
\\[3mm]
\hskip15mm
\displaystyle
\left.
q^{
  \frac{1}{4}m_{i-1}^2
  +\frac{1}{4}(u^L_i+u^R_i)^2
  +\frac{1}{4}\gamma_i-\frac{1}{4}\beta_i^2}
  \prod_{j=i}^{t-1}
  \left[
  {m_j-\frac{1}{2}(\hboldC\hat{\boldm}^{(i-1)}
                       \!-\!\boldu^{L}\!-\!\boldu^{R})_j\atop m_j}
  \right]_q\right),
\end{array}
\end{displaymath}
where the sum is over all
$(m_{i+1},m_{i+2},\ldots,m_{t-1})\equiv(Q_{i+1},Q_{i+2},\ldots,Q_{t-1})$.


If $i=\sigma^R$ then $u^R_i=1$. 
In this case, by definition, we have either
$a^{R\prime}_i=1$, $e^R_i=0$, $a^R_{i-1}=2$ and $e^R_{i-1}=1$, or
$a^{R\prime}_i=p'_{i-1}-1$, $e^R_i=1$, $a^R_{i-1}=p'_{i-1}-2$ and
$e^R_{i-1}=0$.
It is easily checked that $a^{R\prime}_i\notin\mathcal S^{(i-1)}$.
Then, use of Lemma \ref{AttenGen2Lem} yields:

\begin{equation}\label{RightEasyEq}
\begin{array}{l}
\displaystyle
\hskip-1mm
\mchi^{p_{i-1},p_{i-1}^\prime}_{a^{L\prime}_{i},a^{R\prime}_{i-1},
                   e^L_{i},e^R_{i-1}} (m_{i-1}+u^L_i,m_{i})
                   \left\{\mathcal S^{(i-1)}\right\}
\\[3mm]
\hskip10mm
\displaystyle
=
q^{-\frac12u^R_i(m_{i-1}+u^L_i+u^R_i)
   +\frac12(\sboldDelta^R)_i(\alpha_i+\beta_i)}
\\[2mm]
\hskip35mm
\displaystyle
\mchi^{p_{i-1},p_{i-1}^\prime}_{a^{L\prime}_{i},a^{R\prime}_{i},
                   e^L_{i},e^R_{i}} (m_{i-1}+u^L_i+u^R_i,m_{i})\!
                   \left\{\mathcal S^{(i-1)}\right\}\!.
\end{array}
\end{equation}
If $i\ne\sigma^R$ then (noting that $i\ne t_{k}$) $u^R_i=(\boldDelta^R)_i=0$,
$e^R_{i-1}=e^R_{i}$ and $a^R_{i-1}=a^{R\prime}_i$.
The preceding expression thus also holds in this case.

We also immediately obtain
\begin{displaymath}
\begin{array}{ll}
\alpha^{p_{i-1},p_{i-1}^\prime}_{a^{L\prime}_{i},a^{R}_{i-1}}
&=\alpha_i+\beta_i-(\boldDelta^R)_i;\\[2mm]
\beta^{p_{i-1},p_{i-1}^\prime}_{a^{L\prime}_{i},a^{R}_{i-1},
                                                 e^L_{i},e^R_{i-1}}
&=\beta_i-(\boldDelta^R)_i.
\end{array}
\end{displaymath}


If $i=\sigma^L$ then $u^L_i=1$. 
In this case, by definition, we have either
$a^{L\prime}_i=1$, $e^L_i=0$, $a^L_{i-1}=2$ and $e^L_{i-1}=1$, or
$a^{L\prime}_i=p'_{i-1}-1$, $e^L_i=1$, $a^L_{i-1}=p'_{i-1}-2$ and
$e^L_{i-1}=0$.
It is easily checked that $a^{L\prime}_i\notin\mathcal S^{(i-1)}$.
Then, use of Lemma \ref{AttenGen1Lem} yields:

\begin{equation}\label{LeftEasyEq}
\begin{array}{l}
\displaystyle
\hskip-1mm
\mchi^{p_{i-1},p_{i-1}^\prime}_{a^L_{i-1},a^R_{i-1},e^L_{i-1},e^R_{i-1}}
                  (m_{i-1},m_{i})
                   \left\{\mathcal S^{(i-1)}\right\}
\\[3mm]
\hskip10mm
\displaystyle
=
q^{-\frac12u^L_i(m_{i-1}-m_i+u^L_i)
      -\frac12(\sboldDelta^L)_i(\beta_i-(\sboldDelta^R)_i) }
\\[2mm]
\hskip35mm
\displaystyle
\mchi^{p_{i-1},p_{i-1}^\prime}_{a^{L\prime}_{i},a^{R}_{i-1},
                   e^L_{i},e^R_{i-1}} (m_{i-1}+u^L_i,m_{i})\!
                   \left\{\mathcal S^{(i-1)}\right\}\!.
\end{array}
\end{equation}

If $i\ne\sigma^R$ then (noting that $i\ne t_{k}$) $u^L_i=(\boldDelta^L)_i=0$,
$e^L_{i-1}=e^L_{i}$ and $a^L_{i-1}=a^{L\prime}_i$.
The preceding expression thus also holds in this case.

We also obtain:
\begin{displaymath}
\begin{array}{ll}
\alpha^{p_{i-1},p_{i-1}^\prime}_{a^L_{i-1},a^R_{i-1}}
&=\alpha_i+\beta_i
          -(\boldDelta^R)_{i}
          +(\boldDelta^L)_{i};\\[2mm]
\beta^{p_{i-1},p_{i-1}^\prime}_{a^L_{i-1},a^R_{i-1},e^L_{i-1},e^R_{i-1}}
&=\beta_i-(\boldDelta^R)_{i}
         +(\boldDelta^L)_{i}.
\end{array}
\end{displaymath}

Combining all the above, and using the expression for $\gamma_{i-1}''$
given by (\ref{Const1Eq}) and (\ref{Const2Eq}), yields:

\begin{displaymath}
\begin{array}{l}
\displaystyle
\mchi^{p_{i-1},p_{i-1}^\prime}_{a^L_{i-1},a^R_{i-1},e^L_{i-1},e^R_{i-1}}
                   (m_{i-1},m_{i})
                   \left\{\mathcal S^{(i-1)}\right\}
\\[1.5mm]
\hskip5mm
\displaystyle
=\sum\left(
  q^{\frac{1}{4}\hat{\sboldm}^{(i)T}\sboldC\hat{\sboldm}^{(i)}
   +\frac{1}{4} m_{i-1}^2
   -\frac{1}{2} m_{i-1}m_{i}
   -\frac{1}{2}(\sboldu^{L}_{(\flat,k(i))}+\sboldu^{R}_{(\sharp,k(i))})
   \cdot \sboldm^{(i-1)}
  +\frac{1}{4}\gamma''_{i-1}}
\vrule height20pt depth2pt width0pt
\right.\\[5mm]
\hskip55mm
\displaystyle
\left.
  \prod_{j=i}^{t-1}
  \left[
  {m_j-\frac{1}{2}(\hboldC\hat{\boldm}^{(i-1)}
                       \!-\!\boldu^{L}\!-\!\boldu^{R})_j\atop m_j}
  \right]_q\right)\\[3mm]
\hskip5mm
\displaystyle
=F^{(i-1)}_{a,b}(\boldu^L,\boldu^R,m_{i-1},m_i),
\end{array}
\end{displaymath}
which is the required result when $i\ne t_k$, since $k(i)=k(i-1)$.

In this $i\ne t_k$ case, making use of (\ref{Const1Eq}),
(\ref{Const2Eq}), we also immediately obtain:
\begin{displaymath}
\begin{array}{ll}
\alpha^{p_{i-1},p_{i-1}^\prime}_{a^L_{i-1},a^R_{i-1}}
&=\alpha_{i-1}'';\\[2mm]
\beta^{p_{i-1},p_{i-1}^\prime}_{a^L_{i-1},a^R_{i-1},e^L_{i-1},e^R_{i-1}}
&=\beta_{i-1}'.
\end{array}
\end{displaymath}


Now consider the case for which $i=t_k$.
Equation (\ref{Const3Eq}) gives $\alpha_i=\alpha_i''$,
$\beta_i=\alpha_i''-\beta_i'$ and $\gamma_i=-\alpha_i^2-\gamma_i''$.
Corollary \ref{DParamLem} gives
$\alpha^{p_{i}^\prime-p_{i},p_{i}^\prime}_{a^L_{i},a^R_{i}}
=\alpha_i$ and
$\beta^{p_{i}^\prime-p_{i},p_{i}^\prime}_{a^L_{i},a^R_{i},1-e^L_{i},1-e^R_{i}}
=\beta_i$.
Let $m_{i-1}\equiv Q_{i-1}$.
On setting $M=m_{i-1}+u^L_i+u^R_i$, equation 
(\ref{ParityDef}) implies that $M\equiv Q_{i+1}$.
Then, use of the induction hypothesis,
Lemmas \ref{DijCor} or Lemma \ref{IntDijCor} as appropriate,
and Lemmas \ref{IndTakLem} and \ref{A_AprimeLem} yields:
\begin{equation}\label{StepHardEq}
\begin{array}{l}
\displaystyle
\mchi^{p_{i},p_{i}^\prime}_{a^{L\prime}_{i},a^{R\prime}_{i},
                   1-e^L_{i},1-e^R_{i}} (M,m_{i};q)
                   \left\{\mathcal S^{(i)\prime}\right\}
\\[3mm]
\hskip10mm
\displaystyle
=\sum_{m_{i+1}\equiv Q_{i+1}}
\left(
q^{\frac{1}{4}(m_i^2+(M-m_i)^2-\alpha_i^{2}-\beta_i^{2})}
\left[{\frac{1}{2}(M+m_i-m_{i+1})}\atop m_i\right]_q
\right.\\[3mm]
\hskip60mm
\displaystyle
\left.
F^{(i)}_{a,b}(\boldu^{L},\boldu^{R},m_i,m_{i+1};q^{-1})
\vrule height20pt depth2pt width0pt
\right),
\end{array}
\end{equation}
where
\begin{displaymath}
\mathcal S^{(i)\prime}=
\left\{ \begin{array}{ll}
    \left\{ \kappa^{(i-1)}_{t_n-i+1} \right\}
       &\mbox{if } i<t_n,\sigma^L<t_n,\sigma^R<t_n
        \mbox{ and } a,b\in\mathcal T;\\[1.5mm]
    \left\{ p'_i-\kappa^{(i-1)}_{t_n-i+1} \right\}
       &\mbox{if } i<t_n,\sigma^L<t_n,\sigma^R<t_n
        \mbox{ and } a,b\in\mathcal T';\\[1.5mm]
    \emptyset
       &\mbox{otherwise},
         \end{array} \right.
\end{displaymath}
using a similar argument to that in the $i\ne t_{k(i)}$ case.
Here, Lemma \ref{ParamLem} also gives
$\alpha^{p_{i-1},p_{i-1}^\prime}_{a^{L\prime}_{i},a^{R\prime}_{i}}
=\alpha_i+\beta_i$, and
$\beta^{p_{i-1},p_{i-1}^\prime}_{a^{L\prime}_{i},a^{R\prime}_{i},
1-e^L_{i},1-e^R_{i}}=\beta_i$.

Now set $M=m_{i-1}+u^L_i+u^R_i$, whence on noting that $i=t_k$,
\begin{displaymath}
M+m_i-m_{i+1}=2m_i-(\hboldC\hat{\boldm}^{(i-1)}-\boldu^{L}-\boldu^{R})_i
\end{displaymath}
(in the case $i=t-1$, we require this expression after substituting
$m_t=0$),
and
\begin{displaymath}
\begin{array}{l}
\displaystyle
\hat{\boldm}^{(i+1)T}\boldC\hat{\boldm}^{(i+1)}-m_i^2+m_i^2+(M-m_i)^2\\
\qquad=
\hat{\boldm}^{(i)T}\boldC\hat{\boldm}^{(i)}+M^2-2Mm_i\\
\qquad=
\hat{\boldm}^{(i)T}\boldC\hat{\boldm}^{(i)}+m_{i-1}^2-2m_im_{i-1}\\
\hskip30mm
+2(m_{i-1}-m_i)(u^L_i+u^R_i)+(u^L_i+u^R_i)^2.
\end{array}
\end{displaymath}
Use of (\ref{Inverse1Eq}) or (\ref{Inverse2Eq}) then gives:
\begin{displaymath}
\begin{array}{l}
\displaystyle
\mchi^{p_{i-1},p_{i-1}^\prime}_{a^{L\prime}_{i},a^{R\prime}_{i},
                       1-e^L_{i},1-e^R_{i}} (m_{i-1}+u^L_i+u^R_i,m_{i})
                       \left\{\mathcal S^{(i)\prime}\right\}
\\[3mm]
\displaystyle
=\sum
\left(
q^{\frac{1}{4}\hat{\sboldm}^{(i)T}\sboldC\hat{\sboldm}^{(i)}
  -\frac{1}{2}m_im_{i-1}
  +\frac{1}{2}(m_{i-1}-m_i)(u^L_i+u^R_i)
  -\frac{1}{2}(\sboldu^{L}_{(\flat,k(i)-1)}+\sboldu^{R}_{(\sharp,k(i)-1)})
    \cdot\sboldm^{(i)}}
\vrule height20pt depth2pt width0pt
\right.
\\[3mm]
\hskip15mm
\displaystyle
\left.
q^{
  \frac{1}{4}m_{i-1}^2
  +\frac{1}{4}(u^L_i+u^R_i)^2
  +\frac{1}{4}\gamma_i-\frac{1}{4}\beta_i^2}
  \prod_{j=i}^{t-1}
  \left[
  {m_j-\frac{1}{2}(\hboldC\hat{\boldm}^{(i-1)}
                       \!-\!\boldu^{L}\!-\!\boldu^{R})_j\atop m_j}
  \right]_q
\right)
\!,
\end{array}
\end{displaymath}
where the sum is over all
$(m_{i+1},m_{i+2},\ldots,m_{t-1})\equiv(Q_{i+1},Q_{i+2},\ldots,Q_{t-1})$.


Now set $\mathcal S^{(i)R}=\mathcal S^{(i)\prime}\cup a^{R\prime}_i$
if $i>\sigma^R$ and $\mathcal S^{(i)R}=\mathcal S^{(i)\prime}$
otherwise.

Since $i=t_k$, it follows that $u^R_i=-1$ if $i>\sigma^R$.
In this case, by definition, we have either
$a^{R\prime}_i=2$, $1-e^R_i=1$, $a^R_{i-1}=1$ and $e^R_{i-1}=0$, or
$a^{R\prime}_i=p'_{i-1}-2$, $1-e^R_i=0$, $a^R_{i-1}=p'_{i-1}-1$ and
$e^R_{i-1}=1$.
Then Lemma \ref{ExtGen2Lem} yields:
\begin{equation}\label{RightHardEq}
\begin{array}{l}
\displaystyle
\mchi^{p_{i-1},p_{i-1}^\prime}_{a^{L\prime}_{i},a^{R}_{i-1},
                       1-e^L_{i},e^R_{i-1}} (m_{i-1}+u^L_i,m_{i})
                       \left\{\mathcal S^{(i)R}\right\}
\\[3mm]
\hskip10mm
\displaystyle
=
q^{-\frac12u^R_i(m_{i-1}+u^L_i+u^R_i)
   +\frac12(\sboldDelta^R)_i(\alpha_i+\beta_i)}
\\[2mm]
\hskip35mm
\displaystyle
\mchi^{p_{i-1},p_{i-1}^\prime}_{a^{L\prime}_{i},a^{R\prime}_{i},
                       1-e^L_{i},1-e^R_{i}} (m_{i-1}+u^L_i+u^R_i,m_{i})
                       \left\{\mathcal S^{(i)\prime}\right\}.
\end{array}
\end{equation}

In addition, the same expression clearly also holds in the case
$i\le\sigma^R$, for which
$u^R_i=(\boldDelta^R)_i=0$,
$e^R_{i-1}=1-e^R_{i}$ and $a^R_{i-1}=a^{R\prime}_i$.
(In the $i=\sigma^R$ case, note that $k(i-1)=k(i)-1=k^R(i)$.)

Lemma \ref{ExtGen2Lem} also implies that:
\begin{displaymath}
\begin{array}{ll}
\alpha^{p_{i-1},p_{i-1}^\prime}_{a^{L\prime}_{i},a^{R}_{i-1}}
&=\alpha_i+\beta_i-(\boldDelta^R)_i;\\[2mm]
\beta^{p_{i-1},p_{i-1}^\prime}_{a^{L\prime}_{i},a^{R}_{i-1},
                                                 e^L_{i},e^R_{i-1}}
&=\beta_i-(\boldDelta^R)_i.
\end{array}
\end{displaymath}


Now set $\mathcal S^{(i)L}=\mathcal S^{(i)R}\cup a^{L\prime}_i$
if $i>\sigma^L$ and $\mathcal S^{(i)L}=\mathcal S^{(i)R}$
otherwise.

Since $i=t_k$, it follows that $u^L_i=-1$ if $i>\sigma^L$.
In this case, by definition, we have either
$a^{L\prime}_i=2$, $1-e^L_i=1$, $a^L_{i-1}=1$ and $e^L_{i-1}=0$, or
$a^{L\prime}_i=p'_{i-1}-2$, $1-e^L_i=0$, $a^L_{i-1}=p'_{i-1}-1$ and
$e^R_{i-1}=1$.
Then Lemma \ref{ExtGen1Lem} yields:
\begin{equation}\label{LeftHardEq}
\begin{array}{l}
\displaystyle
\mchi^{p_{i-1},p_{i-1}^\prime}_{a^L_{i-1},a^R_{i-1},e^L_{i-1},e^R_{i-1}}
                       (m_{i-1},m_{i})
                       \left\{\mathcal S^{(i)L}\right\}
\\[3mm]
\hskip10mm
\displaystyle
=
q^{-\frac12u^L_i(m_{i-1}-m_i+u^L_i)
      -\frac12(\sboldDelta^L)_i(\beta_i-(\sboldDelta^R)_i) }
\\[2mm]
\hskip35mm
\displaystyle
\mchi^{p_{i-1},p_{i-1}^\prime}_{a^{L\prime}_{i},a^{R}_{i},
                       1-e^L_{i},e^R_{i-1}} (m_{i-1}+u^L_i,m_{i})
                       \left\{\mathcal S^{(i)R}\right\}.
\end{array}
\end{equation}

In addition, the same expression clearly also holds in the case
$i\le\sigma^L$, for which
$u^L_i=(\boldDelta^L)_i=0$,
$e^L_{i-1}=1-e^L_{i}$ and $a^L_{i-1}=a^{L\prime}_i$.
(In the $i=\sigma^L$ case, note that $k(i-1)=k(i)-1=k^L(i)$.)

Lemma \ref{ExtGen1Lem} also implies that:
\begin{displaymath}
\begin{array}{ll}
\alpha^{p_{i-1},p_{i-1}^\prime}_{a^L_{i-1},a^R_{i-1}}
&=\alpha_i+\beta_i
          -(\boldDelta^R)_{i}
          +(\boldDelta^L)_{i};\\[2mm]
\beta^{p_{i-1},p_{i-1}^\prime}_{a^L_{i-1},a^R_{i-1},e^L_{i-1},e^R_{i-1}}
&=\beta_i-(\boldDelta^R)_{i}
         +(\boldDelta^L)_{i}.
\end{array}
\end{displaymath}

Combining all the above cases for $i=t_k$ yields:
\begin{displaymath}
\begin{array}{l}
\displaystyle
\mchi^{p_{i-1},p_{i-1}^\prime}_{a^L_{i-1},a^R_{i-1},e^L_{i-1},e^R_{i-1}}
                      (m_{i-1},m_{i})
                       \left\{\mathcal S^{(i)L}\right\}.
\\[1.5mm]
\hskip5mm
\displaystyle
=\sum\left(
  q^{\frac{1}{4}\hat{\sboldm}^{(i)T}\sboldC\hat{\sboldm}^{(i)}
   +\frac{1}{4} m_{i-1}^2
   -\frac{1}{2} m_{i-1}m_{i}
   -\frac{1}{2}(\sboldu^{L}_{(\flat,k(i)-1)}+\sboldu^{R}_{(\sharp,k(i)-1)})
   \cdot \sboldm^{(i-1)}
  +\frac{1}{4}\gamma''_{i-1}}
\vrule height20pt depth2pt width0pt
\right.\\[5mm]
\hskip50mm
\displaystyle
\left.
  \prod_{j=i}^{t-1}
  \left[
  {m_j-\frac{1}{2}(\hboldC\hat{\boldm}^{(i-1)}
                       \!-\!\boldu^{L}\!-\!\boldu^{R})_j\atop m_j}
  \right]_q\right)\\[3mm]
\hskip5mm
\displaystyle
=F^{(i-1)}_{a,b}(\boldu^L,\boldu^R,m_{i-1},m_i).
\end{array}
\end{displaymath}
Once it is established that
\begin{displaymath}
\P^{p_{i-1},p_{i-1}^\prime}_{a^L_{i-1},a^R_{i-1},e^L_{i-1},e^R_{i-1}}\!
                      (m_{i-1},m_{i})
                       \left\{\mathcal S^{(i)L}\right\}
=\P^{p_{i-1},p_{i-1}^\prime}_{a^L_{i-1},a^R_{i-1},e^L_{i-1},e^R_{i-1}}\!
                      (m_{i-1},m_{i})
                       \left\{\mathcal S^{(i-1)}\right\}.
\end{displaymath}
we obtain the required result when $i=t_k$, since $k(i)=k(i-1)+1$.

If $i=t_n$ then $\{\mathcal S^{(i)L}\}=\{\mathcal S^{(i-1)}\}$
immediately.
Now let $i<t_n$.
For $A\in\{L,R\}$, if $\sigma^A_i=-1$ then necessarily
$\sigma^A_{t_n}=-1$.
In the case that $a^A\in\mathcal T$, this implies that
$\{2,\kappa^{(i)}_{t_n-i}\}\subset\mathcal S^{(i)L}$
and $\kappa^{(i)}_{t_n-i}\in \mathcal S^{(i-1)}$.
Since $a^A_{i-1}=1$, we may drop the element $2$ from
$\mathcal S^{(i)L}$ with no effect.
Similar reasoning holds for $a^A\in\mathcal T'$
whereupon the claim is established.

In this $i=t_k$ case, making use of (\ref{Const1Eq}),
(\ref{Const2Eq}), we also immediately obtain:
\begin{displaymath}
\begin{array}{ll}
\alpha^{p_{i-1},p_{i-1}^\prime}_{a^L_{i-1},a^R_{i-1}}
&=\alpha_{i-1}'';\\[2mm]
\beta^{p_{i-1},p_{i-1}^\prime}_{a^L_{i-1},a^R_{i-1},e^L_{i-1},e^R_{i-1}}
&=\beta_{i-1}'.
\end{array}
\end{displaymath}

\noindent
The lemma then follows by induction.
\cqfd
\medskip

Before performing a sum over $m_1$, we require the following result.

\begin{lemma}\label{ParityLem}
For $0\le j\le t$,
\begin{displaymath}
\begin{array}{ll}
\alpha_j''&\:\equiv Q_j\;(\mod2);\\[0.5mm]
\beta_j'&\:\equiv Q_j-Q_{j+1}\;(\mod2).
\end{array}
\end{displaymath}
\end{lemma}

\Proof Since $\alpha_t''=0$, $\beta_t'=0$ and $Q_t=Q_{t+1}=0$,
this result is manifest for $j=t$.

We now proceed by downward induction. Thus assume the result holds
for a particular $j>0$.
When $j\ne t_{k(j)}$, equations (\ref{Const3Eq}) and (\ref{Const1Eq})
imply that $\beta_{j-1}'=\beta_j'+(\boldu^L)_j-(\boldu^R)_j$.
Equation (\ref{ParityDef}) implies that
$Q_{j-1}\equiv Q_{j+1}-(\boldu^L)_j-(\boldu^R)_j$.
Thus the induction hypothesis immediately gives
$\beta_{j-1}'\equiv Q_{j-1}-Q_{j}$ in this case.

When $j= t_{k(j)}$, equations (\ref{Const3Eq}) and (\ref{Const1Eq})
imply that $\beta_{j-1}'=\alpha_j''-\beta_j'+(\boldu^L)_j-(\boldu^R)_j$.
Equation (\ref{ParityDef}) implies that
$Q_{j-1}\equiv Q_{j}+Q_{j+1}-(\boldu^L)_j-(\boldu^R)_j$.
Thus the induction hypothesis also gives
$\beta_{j-1}'\equiv Q_{j-1}-Q_{j}$ in this case.

In both cases, equations (\ref{Const3Eq}), (\ref{Const1Eq}) and
(\ref{Const2Eq}) give $\alpha_{j-1}''=\alpha_j''+\beta_{j-1}''$,
whence the induction hypothesis immediately gives
$\alpha_j''\equiv Q_{j-1}$ as required.
\cqfd
\medskip

Define:
\begin{displaymath}
\begin{array}{l}
\displaystyle
F_{a,b}(\boldu^L,\boldu^R,L;q)
\\[1mm]
\displaystyle
\hskip2mm
=
\hskip-2mm
\sum_{\sboldm\equiv\sboldQ(\sboldu^L+\sboldu^R)}
  \hskip-5mm
  q^{\frac{1}{4}\hat{\sboldm}^T\sboldC\hat{\sboldm}-\frac{1}{4} L^2
   -\frac{1}{2}(\sboldu^L_\flat+\sboldu^R_\sharp)\cdot\sboldm
  +\frac{1}{4}\gamma}
  \prod_{j=1}^{t-1}
  \left[
  {m_j-\frac{1}{2}(\hboldC\hat{\boldm}\!-\!\boldu^L\!-\!\boldu^R)_j\atop m_j}
  \right]_q\!\!.
\end{array}
\end{displaymath}
The summation here is over all vectors $\boldm=(m_1,m_2,\ldots,m_{t-1})$
such that $m_j\in\Z_{\ge0}$ and $m_j\equiv Q_j\,(\mod2)$ for $1\le j<t$.
Then, $\hat{\boldm}=(m_0,m_1,m_2,\ldots,m_{t-1})$.

On defining
\begin{displaymath}
\mathcal S=
\left\{ \begin{array}{ll}
    \left\{ \kappa_i \right\}
       &\mbox{if } \sigma^L<t_n,\sigma^R<t_n
        \mbox{ and } a,b\in\mathcal T;\\[1.5mm]
    \left\{ p'_i-\kappa_{i} \right\}
       &\mbox{if } \sigma^L<t_n,\sigma^R<t_n
        \mbox{ and } a,b\in\mathcal T';\\[1.5mm]
    \emptyset
       &\mbox{otherwise},
         \end{array} \right.
\end{displaymath}
we then obtain:

\begin{lemma}\label{CoreEFLem}
Let $p'>2p$.
If $L\equiv\alpha^{p,p'}_{a,b}$ then
\begin{displaymath}
\mchi^{p,p'}_{a,b,e^L_0,e^R_0}(L)\left\{\mathcal S\right\}
=F_{a,b}(\boldu^L,\boldu^R,L).
\end{displaymath}
In addition, $\delta^{p,p'}_{b,e^R_0}=0$.
\end{lemma}

\Proof Lemma \ref{ParityLem} implies that $L\equiv Q_0$.
Lemma \ref{ResPathGenLem} requires the sum over all
$m_1\equiv L+\beta^{p,p'}_{a,b,e,f}$ of the $i=0$ case of
Lemma \ref{CoreIndLem}.
This is applicable since for such $m_1$, Lemma \ref{ParityLem} implies
that $m_1\equiv Q_1$.

The lemma follows after noting that in the $p'>2p$ case,
$\hat{\boldm}^{(1)T}C\hat{\boldm}^{(1)}+L^2-2Lm_{1}
=\hat{\boldm}^{T}C\hat{\boldm}-L^2$
and $\gamma''_0=\gamma$.
\cqfd
\medskip

We now transfer this result to the original weighting function of
(\ref{WtDef}).
To do this we require the value of $c$ given by (\ref{CEq}).
Then, defining
$\ochi^{p,p'}_{a,b,c}(L)\left\{\mathcal S\right\}$
in the way analogous to
$\mchi^{p,p'}_{a,b,e,f}(L)\left\{\mathcal S\right\}$,
we obtain:

\begin{lemma}\label{CoreABCLem}
If $L\equiv\alpha^{p,p'}_{a,b}\,(\mod2)$ then
\begin{displaymath}
\ochi^{p,p'}_{a,b,c}(L)\left\{\mathcal S\right\}
=F_{a,b}(\boldu^L,\boldu^R,L).
\end{displaymath}
\end{lemma}

\Proof
For the moment, assume that $p'>2p$.
Consider $h\in\P^{p,p'}_{a,b,e,f}(L)$ and $h'\in\P^{p,p'}_{a,b,c'}(L)$
given by $h'_i=h_i$ for $0\le i\le L$.
If $\delta^{p,p'}_{b,f}=0$ and $c'=b+(-1)^f$ then,
as noted in Section \ref{WingSec}, $\mwt(h)=\owt(h')$.
Consequently,
$\mchi^{p,p'}_{a,b,e,f}(L)\left\{\mathcal S\right\}
=\ochi^{p,p'}_{a,b,c'}(L)\left\{\mathcal S\right\}$.
However, if $b$ is interfacial then the same is true for $c'=b\pm1$.
As noted at the end of Section \ref{ModComSec},
$b$ is interfacial if $\sigma^R\ge t_1$.
Otherwise, the current lemma follows from noting that for the
$c$ defined above, $c=b+(-1)^{e^R_0}$.

Now given $h\in\P^{p,p'}_{a,b,c}(L)$, define
$\hat h\in\P^{p'-p,p'}_{a,b,c}(L)$ by
by $\hat h_i=h_i$ for $0\le i\le L$.
As in Lemma \ref{DresLem},
$\wt(\hat h)=\frac14(L^2-\alpha^2)-\wt(h)$,
where $\alpha=\alpha^{p,p'}_{a,b}$.
Therefore
$\ochi^{p,p'}_{a,b,c}(L)\left\{\mathcal S\right\}=
q^{\frac14(L^2-\alpha^2)}
\ochi^{p,p'}_{a,b,c}(L;q^{-1})\left\{\mathcal S\right\}$.
Since $\alpha^{p,p'}_{a,b}=\alpha''_0$ by Lemma \ref{CoreIndLem},
and $\gamma_0=-(\alpha''_0)^2-\gamma_0''$ by (\ref{Const3Eq}),
the $p'<2p$ case follows from the $p'>2p$ case obtained
above after using 
$\left[{m+n\atop m}\right]_{q^{-1}}=q^{-mn}\left[{m+n\atop m}\right]_q$,
and noting the change in the definition of $\boldC$.
\cqfd
\bigskip

\noindent {\it Proof of Theorem \ref{Ferm1Thrm}: }
First consider the case where $a<y_n$ and $b<y_n$.
Then necessarily $a,b\in\mathcal T$.
Since $y_n=\kappa_{t_n}$, we have $\sigma^L<t_n$ and $\sigma^R<t_n$.
Thereupon, $\mathcal S=\{y_n\}$.
Let $h\in\P^{p,p'}_{a,b,c}(L)\backslash\P^{p,p'}_{a,b,c}(L)\{y_n\}$.
Then $1\le h_i<y_n$ for $0\le i\le L$.
Since, by Lemma \ref{SegmentLem}, the lowermost $y_n-2$ bands of
the $(p,p')$-model have exactly the same parities as the
corresponding bands of the $(z_n,y_n)$-model,
we see that if $h'\in\P^{z_n,y_n}_{a,b,c}(L)$ is defined by
$h'_i=h_i$ for $0\le i\le L$ then $\owt(h')=\owt(h)$.
Since all of $\P^{z_n,y_n}_{a,b,c}(L)$ arises in this way,
we have
$\ochi^{p,p'}_{a,b,c}(L)=
\ochi^{p,p'}_{a,b,c}(L)\{y_n\}+\ochi^{z_n,y_n}_{a,b,c}(L)$.
This proves the first case of Theorem \ref{Ferm1Thrm}.

The second case arises if $a>p'-y_n$ and $b>p'-y_n$.
Here, necessarily $a,b\in\mathcal T'$,
whence again $\sigma^L<t_n$ and $\sigma^R<t_n$.
The argument proceeds as above, noting that both the
$(p,p')$- and $(z_n,y_n)$-models are up-down symmetric.

The other cases are immediate since $\mathcal S=\emptyset$.
\cqfd
\medskip


\subsection{The $\boldm\boldn$-system}\label{MNsysSec}

Each term in the fermionic expressions given by
Theorem \ref{Ferm1Thrm} or Theorem \ref{Ferm2Thrm} corresponds
to a vector $\boldm=(m_1,m_2,\ldots,m_{t-1})$ where
$\boldm\equiv\boldQ(\boldu^L+\boldu^R)$.
As usual, we set $\hat{\boldm}=(L,m_1,m_2,\ldots,m_{t-1})$.
Now, for each $\boldm$, define a vector
$\boldn=(n_1,n_2,\ldots,n_t)$ by
\begin{equation}\label{nDefEq}
{\boldn}=\frac12(-\hat{\boldC}\hat{\boldm}+\boldu).
\end{equation}
In view of (\ref{ParityDef}), we see that $n_j\in\Z$ for $1\le j\le t$.
Then since
\begin{equation}
\frac{1}{2}(\boldC\hat{\boldm}\!-\!\boldu^L\!-\!\boldu^R)_j=-n_j,
\end{equation}
in those terms that provide a non-zero
contribution to the fermionic expression of
Theorem \ref{Ferm1Thrm}, $n_j\ge0$ for $1\le j\le t$.

On examining the proof of Lemma \ref{CoreIndLem},
we see that $n_i$ is the number of particles added at the
$i$th induction step to pass from
$\P^{p_i,p_i^\prime}_{a^L_i,a^R_i,e^L_i,e^R_i}
              (m_i,m_{i+1})\left\{\mathcal S^{(i)}\right\}$
to
$\P^{p_{i-1},p_{i-1}^\prime}_{a^L_{i-1},a^R_{i-1},e^L_{i-1},e^R_{i-1}}
              (m_{i-1},m_{i})\left\{\mathcal S^{({i-1})}\right\}$.

The set of equations that link the two vectors
$\hat{\boldm}$ and $\boldn$ is known as the $\boldm\boldn$-system.
On account of (\ref{CDefEq}), the equations are more explicitly given by,
for $1\le j\le t$:
\begin{eqnarray}
&&m_{j-1}-m_{j+1}=m_{j}+2n_{j}-u_{j}\qquad
\hbox{if $j=t_k,\quad k=1,2,\ldots,n$;}
\label{MNEq1}\\
&&m_{j-1}+m_{j+1}=2m_j+2n_j-u_{j}\quad\hbox{ otherwise,}
\label{MNEq2}
\end{eqnarray}
where we set $m_{t}=m_{t+1}=0$.





Using these two expressions, and setting $m_0=L$, it may be
readily shown that:
\begin{equation}\label{PartProbEq}
\sum_{i=1}^{t} l_in_i=\frac{1}{2}\left(L+\sum_{i=1}^t l_iu_i\right).
\end{equation}
Thereupon, the summands in the expression for
$F_{a,b}(\boldu^L,\boldu^R,L)$ given in Theorem \ref{Ferm1Thrm}
correspond to solutions of (\ref{PartProbEq})
with each $n_i$ a non-negative integer.

\subsection{The second fermionic form}\label{BadFermSec}

The proof of Theorem \ref{Ferm2Thrm} follows the same lines
as that of Theorem \ref{Ferm1Thrm}.
We will not give the full description, but indicate how the proof
of Lemma \ref{CoreIndLem} is affected by the use of the modified Gaussians.
We first define $F^{(i)\prime}_{a,b}(\boldu^L,\boldu^R,m_i,m_{i+1};q)$
for $0\le i<t$ in the same way as
$F^{(i)}_{a,b}(\boldu^L,\boldu^R,m_i,m_{i+1};q)$ in (\ref{FDef1Eq})
and (\ref{FDef2Eq}), except employing the
modified Gaussians instead of the classical Gaussians.
Note that this modified form of the Gaussian differs from the form
defined in (\ref{Gaussian}) if and only if $A<0$ and $B\ge0$.
In this case, $\left[ {A \atop B} \right]=0$.
In addition, since $\left[{m+n\atop m}\right]'_{q^{-1}}
=q^{-mn}\left[{m+n\atop m}\right]'_q$,
it follows that the analogues of (\ref{Inverse1Eq}) and (\ref{Inverse2Eq})
hold.

\begin{lemma}\label{CoreInd2Lem}
Let $0\le i<t$, $m_i\equiv Q_i$ and $m_{i+1}\equiv Q_{i+1}$.
If $m_i\ge0$ then:
\begin{equation}\label{Ind2Eq}
\mchi^{p_i,p_i^\prime}_{a^L_i,a^R_i,e^L_i,e^R_i}
              (m_i,m_{i+1})
=
F^{(i)}_{a,b}(\boldu^{L},\boldu^{R},m_i,m_{i+1}).
\end{equation}
In addition,
$\alpha^{p_i,p_i^\prime}_{a^L_i,a^R_i}=\alpha_i''$ and
$\beta^{p_i,p_i^\prime}_{a^L_i,a^R_i,e^L_i,e^R_i}=\beta_i'$.
\end{lemma}

\Proof
The proof proceeds much as in the proof of \ref{CoreIndLem}.
However, we must certainly check that using the modified
Gaussians does not introduce unwanted terms.

Consider the $i\ne t_{k(i)}$ case. Combining the analogues
of (\ref{StepEasyEq}), (\ref{RightEasyEq}) and (\ref{LeftEasyEq})
yields:

\begin{displaymath}
\begin{array}{l}
\displaystyle
\hskip-1mm
\mchi^{p_{i-1},p_{i-1}^\prime}_{a^L_{i-1},a^R_{i-1},e^L_{i-1},e^R_{i-1}}
                  (m_{i-1},m_{i})
\\[3mm]
\displaystyle
=
\hskip-3mm
\sum_{\scriptstyle m_{i+1}\equiv Q_{i+1}\atop
       \scriptstyle 0\le m_{i+1}\le m_{i}+1}
\hskip-5mm
q^{\frac12\left(
   m_iu^L_i
   -m_{i-1}(u^L_i+u^R_i)
   -u^L_iu^R_i-2
   +\beta_i((\sboldDelta^R)_i-(\sboldDelta^L)_i)
   +\alpha_i(\sboldDelta^R)_i
   +(\sboldDelta^L)_i(\sboldDelta^R)_i\right)}
\\[2mm]
\hskip23mm
\displaystyle
\times\quad
q^{\frac{1}{4}(M-m_i)^2-\frac{1}{4}\beta^{2}_i}
\left[{\frac{1}{2}(M+m_{i+1})}\atop m_i\right]_q
F^{(i)\prime}_{a,b}(\boldu^{L},\boldu^{R},m_i,m_{i+1}),
\end{array}
\end{displaymath}
where $M=m_{i-1}+u^L_i+u^R_i$.
Since $m_{i-1},m_{i+1}\ge0$, and $u^L_i,u^R_i\ge0$ (because $i\ne t_{k(i)}$),
we have
\begin{equation}\label{ModifiedEq}
\left[{\frac{1}{2}(m_{i-1}+m_{i+1}+u^L_i+u^R_i)}\atop m_i\right]'_q
=\left[{\frac{1}{2}(m_{i-1}+m_{i+1}+u^L_i+u^R_i)}\atop m_i\right]_q.
\end{equation}
The induction step for $i\ne t_{k(i)}$ then proceeds exactly
as in the proof of Lemma \ref{CoreIndLem}.

For the $i=t_{k(i)}$ case, combining the analogues
of (\ref{StepHardEq}), (\ref{RightHardEq}) and (\ref{LeftHardEq})
yields:

\begin{equation}\label{AnalogueEq}
\begin{array}{l}
\displaystyle
\hskip-1mm
\mchi^{p_{i-1},p_{i-1}^\prime}_{a^L_{i-1},a^R_{i-1},e^L_{i-1},e^R_{i-1}}
                  (m_{i-1},m_{i})
                   \left\{\tilde{\mathcal S}\right\}
\\[3mm]
\hskip7mm
\displaystyle
=\sum_{\scriptstyle m_{i+1}\equiv Q_{i+1}\atop
       \scriptstyle 0\le m_{i+1}\le m_{i}+1}
\hskip-6mm
q^{\frac12\left(
   m_iu^L_i
   -m_{i-1}(u^L_i+u^R_i)
   -u^L_iu^R_i
   +(\sboldDelta^L)_i(\sboldDelta^R)_i\right)-1}
\\[3mm]
\hskip29mm
\displaystyle
\times\:
q^{\frac12\left(
   \beta_i((\sboldDelta^R)_i-(\sboldDelta^L)_i)
   +\alpha_i(\sboldDelta^R)_i\right)
   +\frac{1}{4}\left(m_i^2+(M-m_i)^2-\alpha^2_i-\beta^{2}_i\right)}
\\[3mm]
\hskip29mm
\displaystyle
\times\:
\left[{\frac{1}{2}(M\!+\!m_i\!-\!m_{i+1})}\atop m_i\right]_q
F^{(i)\prime}_{a,b}(\boldu^{L},\boldu^{R},m_i,m_{i+1};q^{-1}),
\end{array}
\end{equation}
where $M=m_{i-1}+u^L_i+u^R_i$, and $2\in\tilde{\mathcal S}$
if and only if either $a^L_i=1$ or $a^R_i=1$;
$p'-2\in\tilde{\mathcal S}$
if and only if either $a^L_i=p'-1$ or $a^R_i=p'-1$;
and $\tilde{\mathcal S}$ contains no other values.

We must check that (\ref{AnalogueEq}) holds
if the Gaussian is replaced by its modified form, and the
\lq$\{\tilde{\mathcal S}\}$\rq\ is removed.

If $u^L_i=u^R_i=0$ then $\tilde{\mathcal S}=\emptyset$.
In addition $m_{i+1}\le m_i+1$ implies that:
\begin{equation}\label{Modified2Eq}
\left[{\frac{1}{2}(m_{i-1}\!+\!m_i\!-\!m_{i+1}\!+\!u^L_i\!+\!u^R_i)}
      \atop m_i\right]'_q
=\left[{\frac{1}{2}(m_{i-1}\!+\!m_i\!-\!m_{i+1}\!+\!u^L_i\!+\!u^R_i)}
      \atop m_i\right]_q.
\end{equation}
Thereupon, the induction step for this subcase of
$i=t_{k(i)}$ follows as in the proof of Lemma \ref{CoreIndLem}.

Now consider $u^L_i\ne u^R_i$.
We tackle the case $u^L_i=0$ and $u^R_i=-1$
(the case $u^L_i=-1$ and $u^R_i=0$ is similar).
This implies that $\sigma^L\ge t_{k(i)}$ and
$\sigma^R<t_{k(i)}$.
Then either $a^R_{i-1}=1$ and $\tilde{\mathcal S}=\{2\}$, or
$a^R_{i-1}=p'-1$ and $\tilde{\mathcal S}=\{p'-2\}$.
In addition, $2\le a^L_{i-1}\le p'-2$.
We immediately see that
\begin{equation}\label{GenEqualEq}
\mchi^{p_{i-1},p_{i-1}^\prime}_{a^L_{i-1},a^R_{i-1},e^L_{i-1},e^R_{i-1}}
                  (m_{i-1},m_{i})
                   \left\{\tilde{\mathcal S}\right\}
=\mchi^{p_{i-1},p_{i-1}^\prime}_{a^L_{i-1},a^R_{i-1},e^L_{i-1},e^R_{i-1}}
                  (m_{i-1},m_{i}).
\end{equation}
On the other hand, since $m_{i+1}\le m_i+1$,
(\ref{Modified2Eq}) is valid here unless $m_{i-1}=m_i=0$
and $m_{i+1}=1$.
Now $\sigma^L\ge t_{k(i)}$ implies that if $a^L_i=a^R_i$ then
$\sigma^L=t_{k(i)}$ and $e^L_i=e^R_i$ whereupon
$F^{(i)\prime}_{a,b}(\boldu^{L},\boldu^{R},0,1;q^{-1})=0$.
In this case, since $a^L_{i-1}\ne a^R_{i-1}$, then
$\mchi^{p_{i-1},p_{i-1}^\prime}_{a^L_{i-1},a^R_{i-1},e^L_{i-1},e^R_{i-1}}
                  (0,0)=0$.
Therefore, the induction step holds in this $u^L_i\ne u^R_i$ case.

Now consider $u^L_i=u^R_i=-1$,
so that $\sigma^L<t_{k(i)}$ and $\sigma^R<t_{k(i)}$.
If $a^A\in\mathcal T$ then $a^A_{i-1}=1$, and if
$a^A\in\mathcal T'$ then $a^A_{i-1}=p'-1$.
Thereupon, (\ref{GenEqualEq}) holds unless $m_{i-1}=m_i=0$
and either both $a,b\in\mathcal T$ or both $a,b\in\mathcal T'$.
In these cases,
\begin{equation}\label{GenEqual2Eq}
\begin{array}{l}
\displaystyle
\mchi^{p_{i-1},p_{i-1}^\prime}_{a^L_{i-1},a^R_{i-1},e^L_{i-1},e^R_{i-1}}
                  (0,0)
                   \left\{\tilde{\mathcal S}\right\}
=0;\\[3mm]
\mchi^{p_{i-1},p_{i-1}^\prime}_{a^L_{i-1},a^R_{i-1},e^L_{i-1},e^R_{i-1}}
                  (0,0)
=1,
\end{array}
\end{equation}
by direct enumeration.
On the other hand, (\ref{Modified2Eq}) is valid here unless
$m_{i-1}+m_i-m_{i+1}=0$, and $m_i=0$.
If $m_{i-1}=m_i=0$ then since $\left[{-1\atop 0}\right]'=1$,
and $\alpha_i=\beta_i=0$, the required analogue of
(\ref{AnalogueEq}) holds in this case.
If $m_{i-1}=1$ and $m_i=0$ then both sides of the analogue
of (\ref{AnalogueEq}) are easily seen to be zero.

The induction step is now complete, whence the lemma follows.
\cqfd
\medskip

Note that, at the $i$th step in the induction, an extra term
arises due to the modified Gaussian only if $i=t_{k(i)}$,
$\sigma^L<i$, $\sigma^R<i$ and either both
$a,b\in\mathcal T$ or both $a,b\in\mathcal T'$.
In this case, consider the term 
$F^{(i)\prime}_{a,b}(\boldu^{L},\boldu^{R},m_i,m_{i+1};q^{-1})$,
in (\ref{AnalogueEq}) which enumerates the elements of
$\P^{p_i,p'_i}_{a^L_i,a^R_i,e^L_i,e^R_i}(m_i,m_{i+1})$.
In the case where the extra term arises,
$m_i=m_{i+1}=0$ and either
both $a^L_i=a^R_i=1$ and $e^L_i=e^R_i=0$, or
both $a^L_i=a^R_i=p'-1$ and $e^L_i=e^R_i=1$.
Thus there is precisely one path $\tilde h$ of zero length.

\begin{figure}[ht]
\centerline{\epsfig{file=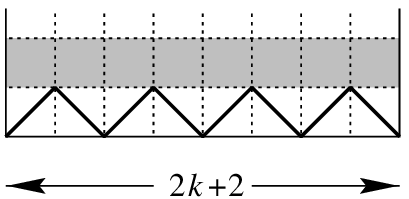}}
\caption{}
\label{HumpyFig}
\end{figure}

Equation (\ref{AnalogueEq}) encapsulates the action of a $\D$-transform
followed by a $\B(k,\lambda)$-transform on $\tilde h$, followed by
extending the result on both sides (since $u^L_i=u^R_i=-1$).
We thus obtain a path of length $m_{i-1}=2k+2$ in the
$(p_{i-1},p'_{i-1})$-model.
This path has the form given in Fig.~\ref{HumpyFig}.
\noindent
That this path contains $n_i=k$ particles, is also encoded in
(\ref{MNEq1}).

When the classical Gaussians are employed, equation (\ref{AnalogueEq})
thus fails to account for the case of a zero length path.
Use of the modified Gaussian remedies this, by permitting
the case $n_i=-1$. This may be viewed as an annihilation of
the $k=0$ case of Fig.~\ref{HumpyFig}, which although appearing
to be a particle (c.f.~Lemma \ref{UniqueLem}), arises through
solely the action of the $B_1$-transform followed by path extension.

\bigskip
\noindent
{\em Acknowledgments:}
We would like to thank Professor Y.~Pugai for collaboration on an 
earlier stage of this work, and on related works, and for many useful 
discussions. His contributions to this work are gratefully acknowledged.  
We also wish to thank Professors A.~Berkovich, B.~McCoy and A.~Schilling 
for many informative discussions on \cite{bms}. Finally, we wish to thank 
Professors M.~Kashiwara and T.~Miwa for the invitation to attend 
{\it `Physical Combinatorics'} where a preliminary version of this work 
was presented, and for their excellent hospitality.


\end{document}